\documentclass[preprint, 11pt]{elsarticle}

\usepackage[T1]{fontenc}
\usepackage[utf8]{inputenc}
\usepackage[hmargin=0.9in]{geometry}
\usepackage[babel,german=quotes]{csquotes}
\usepackage{subcaption}
\usepackage{graphicx}
\usepackage[colorlinks=true,citecolor=blue,urlcolor=blue]{hyperref}
\usepackage{caption}
\usepackage{amsmath}
\usepackage{amsthm}
\usepackage{amssymb}
\usepackage{bm, amsfonts}
\usepackage{xcolor,soul,framed}
\usepackage{fixmath}
\usepackage{tikz}
\usepackage{bm}
\usepackage{makecell}
\usepackage{titlesec}
\usepackage{multirow}
\usepackage{indentfirst}

\colorlet{shadecolor}{green}

\titleformat{\paragraph}
{\normalfont\normalsize\itshape}{\theparagraph}{1em}{}
\titlespacing*{\paragraph}
{0pt}{3.25ex plus 1ex minus .2ex}{1.5ex plus .2ex}

\newdefinition{remark}{Remark}

\makeatletter
\def\ps@pprintTitle{%
\let\@oddhead\@empty
\let\@evenhead\@empty
\def\@oddfoot{
\footnotesize\itshape
\ifx\@journal\@empty Elsevier
\else\@journal\fi
\hfill\today
}%
\let\@evenfoot\@oddfoot}
\makeatother

\begin{document}

\begin{frontmatter}
\title{Dissipative WENO stabilization of high-order
discontinuous Galerkin methods for hyperbolic problems}
\author{Joshua Vedral}
\ead{joshua.vedral@math.tu-dortmund.de}

\address{Institute of Applied Mathematics (LS III), TU Dortmund University\\ Vogelpothsweg 87,
	D-44227 Dortmund, Germany}

\journal{To be submitted}

\begin{abstract}
We present a new approach to stabilizing high-order Runge-Kutta discontinuous Galerkin (RKDG) schemes using weighted essentially non-oscillatory (WENO) reconstructions in the context of hyperbolic conservation laws. In contrast to RKDG schemes that overwrite finite element solutions with WENO reconstructions, our approach employs the reconstruction-based smoothness sensor presented by Kuzmin and Vedral (J. Comput. Phys. 487:112153, 2023) to control the amount of added numerical dissipation. Incorporating a dissipation-based WENO stabilization term into a discontinuous Galerkin (DG) discretization, the proposed methodology achieves high-order accuracy while effectively capturing discontinuities in the solution. As such, our approach offers an attractive alternative to WENO-based slope limiters for DG schemes. The reconstruction procedure that we use performs Hermite interpolation on stencils composed of a mesh cell and its neighboring cells. The amount of numerical dissipation is determined by the relative differences between the partial derivatives of reconstructed candidate polynomials and those of the underlying finite element approximation. The employed smoothness sensor takes all derivatives into account to properly assess the local smoothness of a high-order DG solution. Numerical experiments demonstrate the ability of our scheme to capture discontinuities sharply. Optimal convergence rates are obtained for all polynomial degrees.
\end{abstract}


\end{frontmatter}
%%%
%%% Introduction
%%%
\section{Introduction}
\label{sec:intro}
It is well known that the presence of shocks, discontinuities and steep gradients in exact solutions of hyperbolic conservation laws makes the design of accurate and robust numerical approximations a challenging task. Standard finite volume (FV) and finite element (FE) methods produce numerical solutions which converge to wrong weak solutions and introduce either spurious oscillations or excessive dissipation near shocks. To address this issue, stabilized high-resolution numerical methods were introduced to capture sharp features and discontinuities accurately. Many representatives of such schemes are based on the \textit{essentially non-oscillatory} (ENO) paradigm introduced by Harten et al. \cite{harten1987}. The key idea is to select the smoothest polynomial from a set of candidate polynomials which are constructed using information from neighboring elements. By employing a weighted combination of all polynomial approximations, Liu et al. \cite{liu1994} introduced the framework of weighted ENO (WENO) schemes which allows for better representation of steep gradients, offering more accurate solutions in the presence of shocks and discontinuities. The seminal paper by Jiang and Shu \cite{jiang1996} introduced a fifth-order WENO scheme using a new measurement of the solution smoothness. 

More recently, numerous variants of the original WENO methodology have been developed. A~significant finding in \cite{henrick2005} revealed that the WENO scheme introduced in \cite{jiang1996} may encounter accuracy issues when the solution comprises critical points. For instance, the fifth-order scheme presented in \cite{jiang1996} degrades to third order in the vicinity of smooth extrema. To address this concern, a mapping function was proposed to correct the weights, resulting in the WENO-M schemes \cite{henrick2005}. Moreover, central WENO (CWENO) schemes \cite{levy1999,levy2000,qiu2002} were introduced by splitting the reconstruction stencil into a central stencil and a WENO stencil. The contributions of both stencils are blended using suitable weights. The so-called WENO-Z schemes \cite{acker2016,castro2011,wang2018} reduce the dissipative error by using a linear combination of existing low-order smoothness indicators to construct a new high-order one.

Qiu and Shu \cite{qiu2005} applied the WENO framework to Runge-Kutta discontinuous Galerkin (RKDG) methods by introducing a limiter based on WENO reconstructions.  They further investigated the use of Hermite WENO (HWENO) schemes as limiters for RKDG schemes in one \cite{qiu2004} and two \cite{qiu2005b} dimensions. An extension to unstructured grids was presented by Zhu and Qiu \cite{zhu2009}. To account for the local properties of DG schemes, Zhong and Shu \cite{zhong2013} developed a limiter which uses information only from immediate neighboring elements. Their idea of reconstructing the entire polynomial instead of reconstructing point values was adopted by Zhu et al. for structured \cite{zhu2016} and unstructured \cite{zhu2017} grids. The simplified limiter presented by Zhang and Shu \cite{zhang2010} makes it possible to construct uniformly high-order accurate FV and DG schemes satisfying a strict maximum principle for scalar conservation laws on rectangular meshes. Their approach was later generalized to positivity preserving schemes for the compressible Euler equations of gas dynamics on rectangular \cite{zhang2010c} and triangular \cite{zhang2012} meshes. An overview of existing DG-WENO methods can be found, e.g., in \cite{shu2016,zhang2011}. 

The use of WENO reconstructions as limiters is clearly not an option for continuous Galerkin (CG) methods. To our best knowledge, the first successful attempt to construct WENO-based stabilization in the CG setting was recently made in \cite{kuzmin2023a}. Instead of overwriting polynomial solutions with WENO reconstructions, the approach proposed in \cite{kuzmin2023a} uses
a reconstruction-based smoothness sensor to adaptively blend the numerical viscosity operators of high- and low-order stabilization terms.

In this work, we extend the methodology developed in \cite{kuzmin2023a}
to DG methods and hyperbolic systems. Building upon the standard DG scheme, we introduce a nonlinear stabilization technique that incorporates low-order dissipation in the vicinity of discontinuities. The smoothness sensor developed in \cite{kuzmin2023a} is employed to control the amount of added diffusion. The resulting scheme achieves high-order accuracy while effectively capturing discontinuities in the solution. Moreover, the stabilized space discretization corresponds to a consistent weighted residual formulation of the problem at hand.

In the next section, we provide a review of standard RKDG methods. Section \ref{sec:stab} introduces the dissipation-based nonlinear stabilization technique. The WENO-based smoothness sensor is presented in Section \ref{sec:sensor}. To assess the performance of the proposed methodology, we perform a series of numerical experiments in Section \ref{sec:numex}. The paper ends with concluding remarks in Section \ref{sec:concl}.

%%%
%%% RKDG
%%%
\section{Standard RKDG methods}
\label{sec:rkdg}
Let $U(\mathbf{x},t)\in \mathbb{R}^m$, $m\in \mathbb{N}$ be the vector of conserved quantities at the space location $\mathbf{x}\in \bar{\Omega}$ and time instant $t\geq 0$. We consider the initial value problem
\begin{subequations}
	\begin{alignat}{2}
	\frac{\partial U}{\partial t}+\nabla\cdot\mathbf{F}(U)&=0 \quad &&\text{in }\Omega \times \mathbb{R}_+, \label{eq:pde}\\
	U(\cdot,0)&=U_0 \quad &&\text{in } \Omega,
	\end{alignat}\label{eq:ivp}%))
\end{subequations}
where $\Omega\in \mathbb{R}^d$, $d\in\{1,2,3\}$ is a Lipschitz domain, $\mathbf{F}(U)\in \mathbb{R}^{m\times d}$ is the flux function and $U_0\in\mathbb{R}^m$ is the initial datum. In the scalar ($m=1$) case, boundary conditions are imposed weakly at the inlet $\partial \Omega_{-}=\{\mathbf{x}\in \partial \Omega:\mathbf{v}\cdot\mathbf{n}<0\}$ of the domain. Here, $\mathbf{v}(\mathbf{x},t)=\mathbf{F}'(U(\mathbf{x},t))$ is a velocity field and $\mathbf{n}$ denotes the unit outward normal to the boundary of the domain. For hyperbolic systems ($m>1$), appropriate boundary conditions (see, e.g., \cite{guaily2013}) need to be imposed.

Let $\mathcal{T}_h=\{K_1,\ldots,K_{E_h}\}$ be a decomposition of the domain $\Omega$ into non-overlapping elements $K_e$, $e=1,...,E_h$. To discretize \eqref{eq:ivp} in space, we use the standard discontinuous Galerkin finite element space  
\begin{equation}
\mathbb{V}_p(\mathcal{T}_h)=\{v_h\in L^2(\Omega):v_h|_{K_e}\in \mathbb{V}_p(K_e)\,\forall K_e\in\mathcal{T}_h\},
\end{equation}
where $\mathbb{V}_p(K_e)\in \{\mathbb{P}_p(K_e), \mathbb{Q}_p(K_e)\}$ is the space spanned by polynomials of degree $p$. We seek an approximate solution $U_h^e\in \mathbb{V}_p(K_e)$ of the form 
\begin{equation}
U_h^e(\mathbf{x},t):=U_h(\mathbf{x},t)|_{K_e}=\sum_{j=1}^{N_e}U_j(t)\varphi_j(\mathbf{x}), \quad \mathbf{x}\in K_e, \quad t\in [0,T],
\label{eq:fesol}
\end{equation}
where $\varphi_j$, $j=1,\ldots,N_e$ are finite element basis functions. We remark that our methodology does not rely on a particular choice of basis functions. Popular choices include, e.g., Lagrange, Bernstein or Legendre-Gauss-Lobatto polynomials. 

Inserting \eqref{eq:fesol} into \eqref{eq:pde}, multiplying by a sufficiently smooth test function $W_h$ and applying integration by parts yields the local weak formulation  
\begin{equation}
\int_{K_e}W_h\cdot\frac{\partial U_h^e}{\partial t}\,\mathrm{d}\mathbf{x}-\int_{K_e}\nabla W_h\cdot\mathbf{F}(U_h^e)\,\mathrm{d}\mathbf{x}+\int_{\partial K_e}W_h\cdot\mathbf{F}(U_h^e)\cdot\mathbf{n}_e\,\mathrm{ds}=0.
\label{eq:weak1}
\end{equation}
The DG solution $U_h$ is generally discontinuous across element interfaces. The one-sided limits 
\begin{equation}
U_h^{e,\pm}(\mathbf{x},t):=\lim_{\varepsilon\to\pm 0}U_h(\mathbf{x}+\varepsilon\mathbf{n}_e,t), \quad \mathbf{x}\in\partial K_e
\end{equation}
are used to construct an approximation  $\mathcal{H}(U_h^{e,-},U_h^{e,+},\mathbf{n}_e)$ to $\mathbf{F}(U_h^e)\cdot \mathbf{n}_e$. A numerical flux function $\mathcal{H}(U_L,U_R,\mathbf{n})$ must provide consistency
\begin{equation}
\mathcal{H}(U,U,\mathbf{n})=\mathbf{F}(U)\cdot \mathbf{n}
\end{equation}
and conservation
\begin{equation}
\mathcal{H}(U,V,\mathbf{n})=-\mathcal{H}(V,U,-\mathbf{n}).
\end{equation}
We employ solely the local Lax-Friedrichs numerical flux
\begin{equation}
\mathcal{H}_{LLF} (U_h^{e,-},U_h^{e,+},\mathbf{n}_e)=\mathbf{n}_e\cdot \frac{\mathbf{F}(U_h^{e,-})+\mathbf{F}(U_h^{e,+})}{2}-\frac{1}{2}s^{+}(U_h^{e,+}-U_h^{e,-}),
\end{equation}
where $s^{+}$ is the maximum wave speed of the 1D Riemann problem in the normal direction $\mathbf{n}_e$, i.e.,
\begin{equation}
s^{+}=\max_{\omega \in [0,1]}|\mathbf{F}'(\omega U_h^{e,-}+(1-\omega)U_h^{e,+})\cdot \mathbf{n}_e|.
\end{equation}
\begin{remark}
	Many stabilized RKDG schemes employ upwind fluxes, the (local) Lax-Friedrichs flux or HLL-type fluxes. While the choice of the numerical flux function is important for piecewise constant approximations, it has little influence on the quality of high-order DG methods \cite{hajduk2022diss}.
\end{remark}
By inserting a numerical flux $\mathcal{H} (U_h^{e,-},U_h^{e,+},\mathbf{n}_e)$ into \eqref{eq:weak1} and replacing $W_h$ with $\varphi_i^e$, we obtain a system of $N_e$ semi-discrete equations
\begin{equation}
\int_{K_e}\varphi_i^e\cdot\frac{\partial U_h^e}{\partial t}\,\mathrm{d}\mathbf{x}=\int_{K_e}\nabla \varphi_i^e\cdot\mathbf{F}(U_h^e)\,\mathrm{d}\mathbf{x}-\int_{\partial K_e}\varphi_i^e\cdot\mathcal{H}(U_h^{e,-},U_h^{e,+},\mathbf{n}_e)\,\mathrm{ds},
\label{eq:weak2}
\end{equation}
which can be written in matrix form as
\begin{equation}
M\frac{\mathrm{d}U}{\mathrm{dt}}=R(U),
\end{equation}
where $M=(m_{ij})$ denotes the mass matrix, $U$ is the vector of degrees of freedom and $R(U)$ is the residual vector.

The numerical solution is evolved in time using the optimal third-order strong stability preserving (SSP) Runge-Kutta method \cite{gottlieb2001}
\begin{align}
\begin{split}
U^{(1)}&=U^n+\Delta t M^{-1}R(U^n),\\
U^{(2)}&=\frac{3}{4}U^n+\frac{1}{4}U^{(1)}+\frac{1}{4}\Delta t M^{-1}R(U^{(1)}),\\
U^{(n+1)}&=\frac{1}{3}U^n+\frac{2}{3}U^{(2)}+\frac{2}{3}\Delta tM^{-1}R(U^{(2)}),
\end{split}
\label{eq:ssp}
\end{align}
unless stated otherwise.

%%% 
%%% Stabilization
%%%
\section{Dissipation-based nonlinear stabilization}
\label{sec:stab}
The RKDG scheme maintains numerical stability in applications to initial value problems that involve smooth solutions or weak shocks. While the semi-discrete scheme satisfies a cell entropy inequality \cite{jiang1994} and, therefore, provides $L^2$ stability in the scalar case \cite{cockburn1998a}, it still tends to generate significant oscillations in the presence of strong discontinuities. 

To address this issue, we introduce artificial diffusion near shocks which effectively dampens the spurious oscillations, ensures stability and improves the accuracy of the RKDG scheme. By incorporating a local stabilization operator $s_h^e(U_h^e,\varphi_i^e)$ into the weak formulation \eqref{eq:weak2}, the general form of a stabilized scheme can be written as 
\begin{equation}
\int_{K_e}\varphi_i^e\cdot\frac{\partial U_h^e}{\partial t}\,\mathrm{d}\mathbf{x}+s_h^e(U_h^e,\varphi_i^e)=\int_{K_e}\nabla \varphi_i^e\cdot\mathbf{F}(U_h^e)\,\mathrm{d}\mathbf{x}-\int_{\partial K_e}\varphi_i^e\cdot\mathcal{H}(U_h^{e,-},U_h^{e,+},\mathbf{n}_e)\,\mathrm{ds}.
\end{equation}
Since piecewise-constant (DG-$\mathbb{P}_0$) approximations naturally introduce sufficient diffusion to ensure numerical stability, a straightforward and practical approach is to add isotropic artificial diffusion throughout the computational domain. In this context, the corresponding low-order nonlinear stabilization operator reads
\begin{equation}
s_h^{e,L}(U_h^e,\varphi_i^e)=\nu_e\int_{K_e}\nabla \varphi_i^e\cdot \nabla U_h^e\,\mathrm{d}\mathbf{x},
\label{eq:lostab}
\end{equation}
where $\nu_e$ is the viscosity parameter given by
\begin{equation}
\nu_e=\frac{\lambda_e h_e}{2p}.
\label{eq:visc}
\end{equation}
Here, $\lambda_e=\|\mathbf{F}'(U_h^e)\|_{L^{\infty}(K_e)}$ denotes the maximum wave speed and $h_e$ is the local mesh size.

To apply stabilization only in regions with discontinuities, a shock detector $\gamma_e\in [0,1]$ is introduced to control the amount of added diffusion. In our numerical experiments, we use the adaptive nonlinear stabilization operator 
\begin{equation}
s_h^{e,A}(U_h^e,\varphi_i^e)=\gamma_e\nu_e\int_{K_e}\nabla \varphi_i^e\cdot \nabla U_h^e\,\mathrm{d}\mathbf{x}.
\label{eq:astab}
\end{equation}
Until now, the employed stabilization technique is completely independent of the WENO framework. The accuracy and stability of the resulting scheme depend solely on the definition of the smoothness indicator $\gamma_e$, which we give in the next section.
\begin{remark}
	In contrast to the postprocessing nature of standard WENO-based RKDG schemes, in which the numerical solution is overwritten by WENO reconstructions, our stabilization technique is directly integrated into the semi-discrete problem. This incorporation is akin to standard shock-capturing schemes, which makes our stabilization technique well suited for further modifications, such as flux limiting \cite{kuzmin2021a,kuzmin2021,rueda2023}.	
\end{remark}
\begin{remark}
	The stabilized continuous Galerkin (CG) version presented in \cite{kuzmin2023a} adaptively blends the numerical viscosity operators of low- and high-order stabilization terms. However, the purpose of employing high-order linear stabilization in \cite{kuzmin2023a} was primarily to obtain optimal convergence rates for problems with smooth exact solutions. Given the inherent properties of DG methods, our scheme does not require the incorporation of high-order linear stabilization. Instead, it relies solely on the low-order nonlinear stabilization operator \eqref{eq:astab} equipped with an appropriate shock detector $\gamma_e$.
\end{remark}

%%%
%%% Smoothness sensor
%%%
\section{WENO-based smoothness sensor}
\label{sec:sensor}
A well-designed shock-capturing method for DG discretizations of hyperbolic problems should suppress spurious oscillations in a manner that preserves high-order accuracy in smooth regions. To meet this requirement, the use of nonlinear stabilization is commonly restricted to cells that are identified as `troubled' by smoothness sensors. Alongside traditional shock detectors that rely on measures such as entropy \cite{guermond2011,krivodonova2004,lv2016}, total variation and slope \cite{ducros1999,harten1978,harten1983a,hendricks2018,jameson1981,ren2003} or residual \cite{marras2018,stiernstrom2021} of a finite element approximation, remarkable progress has been made in the development of smoothness sensors based on the WENO methodology \cite{hill2004,li2020,movahed2013,visbal2005,wang2023,zhao2019,zhao2020}. Building upon these advances, we present a smoothness sensor that aligns with the principles of WENO-based shock detection. 

Let $U_h^{e,*}$ denote a WENO reconstruction. We equip our stabilization operator \eqref{eq:astab} with the smoothness sensor \cite{kuzmin2023a}
\begin{equation}
\gamma_e = \min\Bigg(1, \frac{\|U_h^e-U_h^{e,*}\|_e}{\|U_h^e\|_e}\Bigg)^q,
\label{eq:sensor}
\end{equation}
where $q\geq 1$ allows for tuning the sensitivity of the relative difference between $U_h^e$ and $U_h^{e,*}$. To accurately evaluate the smoothness of the numerical solution, we use the scaled Sobolev semi-norm \cite{friedrich1998,jiang1996}
\begin{equation}
\|v\|_e=\Bigg(\sum_{1\leq |\mathbf{k}|\leq p}h_e^{2|\mathbf{k}|-d}\int_{K_e}|D^{\mathbf{k}}v|^2\,\mathrm{d}\mathbf{x}\Bigg)^{1/2} \quad \forall v\in H^p(K_e),
\end{equation}
where $\mathbf{k}=(k_1,\ldots,k_d)$ is the multiindex of the partial derivative
\begin{equation}
D^{\mathbf{k}}v=\frac{\partial^{|\mathbf{k}|}v}{\partial x_1^{k_1}\cdots \partial x_d^{k_d}}, \quad |\mathbf{k}|=k_1+\ldots+k_d.
\end{equation}
Jiang and Shu \cite{jiang1996} and Friedrich \cite{friedrich1998} measured the smoothness of candidate polynomials using derivative-based metrics of the form $\|\cdot\|_e^q$, $q \in \{1,2\}$. In this context, our smoothness sensor shares similarities with existing WENO-based shock sensors \cite{hill2004,movahed2013,visbal2005,zhao2020}.
\begin{remark}
	We modify the smoothness sensor \eqref{eq:sensor} by setting it to zero if the relative difference between the numerical solution and the WENO reconstruction is smaller than a prescribed tolerance. Thus, we effectively avoid introducing unnecessary numerical dissipation in smooth regions, thereby enhancing overall accuracy. A generalized version of the smoothness indicator \eqref{eq:sensor} is given by 
	\begin{align}
	\gamma_{e}^b = \begin{cases}
	\gamma_e\quad &\text{if } \frac{\|U_h^e-U_h^{e,*}\|_e}{\|U_h^e\|_e}\geq b, \\
	0 \quad &\text{otherwise},
	\end{cases}
	\label{eq:gensensor}
	\end{align}
	where $b\in[0,1]$ is the threshold for activating numerical dissipation. We recover the smoothness sensor \eqref{eq:sensor} by using $b=0$.
\end{remark}

We complete the description of our smoothness sensor by discussing the WENO reconstruction polynomial $U_h^{e,*}$. Exploiting the local nature of standard DG schemes, we follow the Hermite WENO (HWENO) approach introduced by Qiu and Shu \cite{qiu2004}. Various extensions of the HWENO method can be found in the literature (see, e.g., \cite{liu2015,luo2007,luo2012,qiu2004,zhang2023,zhu2009,zhu2017}). By incorporating information about the derivatives, HWENO schemes allow for a more accurate representation of the solution, especially in the presence of discontinuities and sharp gradients. The weights assigned to candidate polynomials are determined based on smoothness indicators, with larger weights given to less oscillatory polynomials. This weighting strategy effectively reduces numerical oscillations and improves the overall accuracy of the solution. An adaptive HWENO reconstruction can be written as
\begin{equation}
U_h^{e,*}=\sum_{l=0}^{m_e}\omega_l^eU_{h,l}^e \in \mathbb{P}_p(K_e),
\end{equation}
where $\omega_l^e$ and $U_{h,l}^e$, $l=0,\ldots,m_e$ denote nonlinear weights and Hermite candidate polynomials, respectively. 

We remark that the main contribution of our work is not the construction of
$U_h^{e,*}$ but the way in which we use $U_h^{e,*}$  to stabilize the underlying RKDG scheme. Therefore, we omit the details of the HWENO reconstruction procedure employed in this work. For a detailed description of this procedure, we refer the reader to \cite[Sec.~5]{kuzmin2023a}.

%%%
%%% Numerical examples
%%%
\section{Numerical examples}
\label{sec:numex}
\setcounter{secnumdepth}{4}
\setcounter{tocdepth}{4}
To assess the properties of the WENO-based stabilization procedure, we conduct a series of numerical experiments involving both linear and nonlinear scalar test problems, as well as the compressible Euler equations of gas dynamics. Our objective is to demonstrate the optimal convergence behavior of our method for smooth problems, its superb shock-capturing capabilities, and its ability to converge to correct entropy solutions in the nonlinear case. Throughout the  numerical experiments, we utilize Lagrange basis functions of polynomial order $p\in\{1, 2, 3\}$. We set $q=1$ in \eqref{eq:sensor} and deviate slightly from the default setting of the linear weights  in \cite{zhu2017} by setting $\tilde{\omega}_0^e=1-m_e\cdot10^{-2}$, $\tilde{\omega}_l^e=10^{-2}$, $l=1,\ldots,m_e$. Numerical solutions are evolved in time using the third-order explicit SSP Runge-Kutta scheme \eqref{eq:ssp}, unless stated otherwise. For clarity, we assign the following labels to the methods under investigation:
\medskip

\begin{tabular}{c|l}
	DG& discontinuous Galerkin method without any stabilization;\\
	LO& DG + linear low-order stabilization, i.e., using $s_h^{e,L}$ defined  by \eqref{eq:lostab};\\
	WENO& DG + nonlinear stabilization, i.e., using $s_h^{e,A}$ defined by \eqref{eq:astab}.
\end{tabular}
\medskip

We measure the experimental order of convergence (EOC) using \cite{leveque1996,lohmann2017} 
\begin{equation}
p^{EOC} = \log\Bigg(\frac{E_1(h_2)}{E_1(h_1)}\Bigg)\log\Bigg(\frac{h_2}{h_1}\Bigg)^{-1},
\end{equation}
where $E_1(h_1)$ and $E_1(h_2)$ are the $L^1$ errors on two meshes of mesh size $h_1$ and $h_2$, respectively. The implementation of all schemes is based on the open-source C++ library MFEM \cite{anderson2021,mfem}. For visualization purposes, numerical solutions to two-dimensional problems are $L^2$-projected into the same-order space of continuous functions. The corresponding results are visualized using the open source C++ software GLVis \cite{glvis}.
\subsection{One-dimensional linear advection with constant velocity}
To begin, we consider the one-dimensional linear transport problem
\begin{equation}
\frac{\partial u}{\partial t}+v\frac{\partial u}{\partial x}=0 \quad \text{in } \Omega= (0,1)
\label{eq:linadv}
\end{equation}
with constant velocity $v=1$ and periodic boundary conditions. To assess the convergence properties of our scheme, we evolve the smooth initial condition 
\begin{equation}
u_0(x)=\cos(2\pi(x-0.5))
\end{equation}
up to the final time $t=1.0$ using SSP Runge-Kutta methods of order $p+1$. Table \ref{tab:convlinadv} shows the $L^1$ errors and EOCs for Lagrange finite elements of order $p\in\{1,2,3\}$. Both the DG and the WENO scheme deliver the optimal $L^1$ convergence rates, i.e., $\text{EOC}\approx p+1$. Notably, we observe that the error of the WENO scheme converges towards the Galerkin error as the mesh is refined.  In contrast, the LO scheme achieves only first-order accuracy because numerical dissipation is introduced throughout the domain.

% Linear advection
\begin{table}[h!]
	\begin{subtable}{\textwidth}
		\centering
		\begin{tabular}{ccccccc}
			\hline
			&\multicolumn{2}{c}{DG}&\multicolumn{2}{c}{LO}&\multicolumn{2}{c}{WENO}\\
			\hline
			$E_h$ & $\|u_h-u_{\text{exact}}\|_{L^1}$ & EOC & $\|u_h-u_{\text{exact}}\|_{L^1}$ & EOC & $\|u_h-u_{\text{exact}}\|_{L^1}$ & EOC\\
			\hline
			16		 &   7.28e-03 & --       & 2.97e-01 & --       & 7.36e-02 &	--       \\
			32		 &   1.67e-03 &	2.13 & 1.70e-01 & 0.80 & 1.94e-02 & 1.93 \\
			64		 &   3.97e-04 & 2.07 & 9.14e-02 & 0.90 & 2.89e-03 & 2.75 \\
			128	   &   9.69e-05 & 2.04 & 4.75e-02 & 0.95 & 1.72e-04 & 4.07 \\
			256	   &   2.39e-05 & 2.02 & 2.42e-02 & 0.97 & 2.39e-05 & 2.84 \\
			512	   &   5.95e-06 & 2.01 & 1.22e-02 & 0.99 & 5.95e-06 & 2.01 \\
			1024 &   1.48e-06 & 2.00 & 6.13e-03 & 0.99 & 1.48e-06 & 2.00 \\
			2048 &   3.70e-07 & 2.00 & 3.07e-03 & 1.00 & 3.70e-07 & 2.00 \\
			\hline
		\end{tabular}
		\caption{$p=1$}
	\end{subtable}
	\begin{subtable}{\textwidth}
		\centering
		\begin{tabular}{ccccccc}
			\hline
			&\multicolumn{2}{c}{DG}&\multicolumn{2}{c}{LO}&\multicolumn{2}{c}{WENO}\\
			\hline
			$E_h$ & $\|u_h-u_{\text{exact}}\|_{L^1}$ & EOC & $\|u_h-u_{\text{exact}}\|_{L^1}$ & EOC & $\|u_h-u_{\text{exact}}\|_{L^1}$ & EOC\\
			\hline
			16      & 1.59e-04 & --       & 2.41e-01 & --       & 1.59e-04 & --       \\
			32      & 1.95e-05 & 3.02 & 1.35e-01 & 0.84 & 1.95e-05 & 3.03 \\
			64      & 2.43e-06 & 3.01 & 7.14e-02 & 0.92 & 2.43e-06 & 3.01 \\
			128   & 3.03e-07 & 3.00 & 3.67e-02 & 0.96 & 3.03e-07 & 3.00 \\
			256   & 3.78e-08 & 3.00 & 1.86e-02 & 0.98 & 3.78e-08 & 3.00 \\
			512   & 4.73e-09 & 3.00 & 9.39e-03 & 0.99 &	4.73e-09 & 3.00 \\
			1024 & 5.91e-10 & 3.00 & 4.71e-03 & 0.99 & 5.91e-10 & 3.00 \\
			2048 & 7.40e-11 & 3.00 & 2.36e-03 & 1.00 & 7.40e-11 & 3.00 \\
			\hline
		\end{tabular}
		\caption{$p=2$}
	\end{subtable}
	\begin{subtable}{\textwidth}
		\centering
		\begin{tabular}{ccccccc}
			\hline
			&\multicolumn{2}{c}{DG}&\multicolumn{2}{c}{LO}&\multicolumn{2}{c}{WENO}\\
			\hline
			$E_h$ & $\|u_h-u_{\text{exact}}\|_{L^1}$ & EOC & $\|u_h-u_{\text{exact}}\|_{L^1}$ & EOC & $\|u_h-u_{\text{exact}}\|_{L^1}$ & EOC\\
			\hline
			16    & 3.86e-06 & --       & 1.84e-01 & --       & 3.86e-06 & --       \\
			32    & 2.40e-07 & 4.01 & 1.00e-01 & 0.88 & 2.40e-07 & 4.01 \\
			64    & 1.50e-08 & 4.00 & 5.21e-02 & 0.94 & 1.50e-08 & 4.00 \\
			128 & 9.35e-10 & 4.00 & 2.66e-02 & 0.97 & 9.35e-10 & 4.00 \\
			256 & 5.85e-11 & 4.00 & 1.34e-02 & 0.98 & 5.85e-11 & 4.00 \\
			512 & 3.83e-12 & 3.93 & 6.76e-03 & 0.99 & 3.83e-12 & 3.93 \\
			\hline
		\end{tabular}
		\caption{$p=3$}
	\end{subtable}
	\caption{One-dimensional linear advection with constant velocity, grid convergence history for finite elements of degree $p\in\{1,2,3\}$.}
	\label{tab:convlinadv}
\end{table}

To test the robustness of our scheme, we solve \eqref{eq:linadv} using the initial condition 
\begin{equation}
u_0(x) = \begin{cases} 
1 & \mbox{if} \ 0.15 \le x \le 0.45, \\
\Big[\cos\Big(\frac{10\pi}{3}(x-0.7)\Big)\Big]^2 & \mbox{if}  \ 0.55 < x < 0.85, \\
0 & \mbox{otherwise},
\end{cases}
\label{num:advinittwo}
\end{equation}
consisting of a discontinuous and an infinitely differentiable smooth profile.

\begin{figure}[!htb]
	\centering
	\begin{subfigure}[b]{\linewidth}
		\centering
		\begin{tikzpicture}
		\draw[rounded corners] (0, 0) rectangle (12, 0.5) node[pos=.5]{};
		\draw[very thick, color={rgb:red,0;green,0.4470;blue,0.7410}] (6.75,0.25)--(7.25,0.25);
		\node at (2,0.25) (a) {$p=1$};
		\node at (5,0.25) (a) {$p=2$};
		\node at (8,0.25) (a) {$p=3$};
		\node at (11,0.25) (a) {Exact};
		\draw[very thick,color={rgb:red,0.9290;green,0.6940;blue,0.1250}] (0.75,0.25)--(1.25,0.25);
		\draw[very thick,color={rgb:red,0.8500;green,0.3250;blue,0.0980}] (3.75,0.25)--(4.25,0.25);
		\draw[very thick] (9.75,0.25)--(10.25,0.25);
		\end{tikzpicture}
		\vspace*{0.25cm}
	\end{subfigure}
	\begin{subfigure}[b]{.32\linewidth}
		\includegraphics[width=\linewidth]{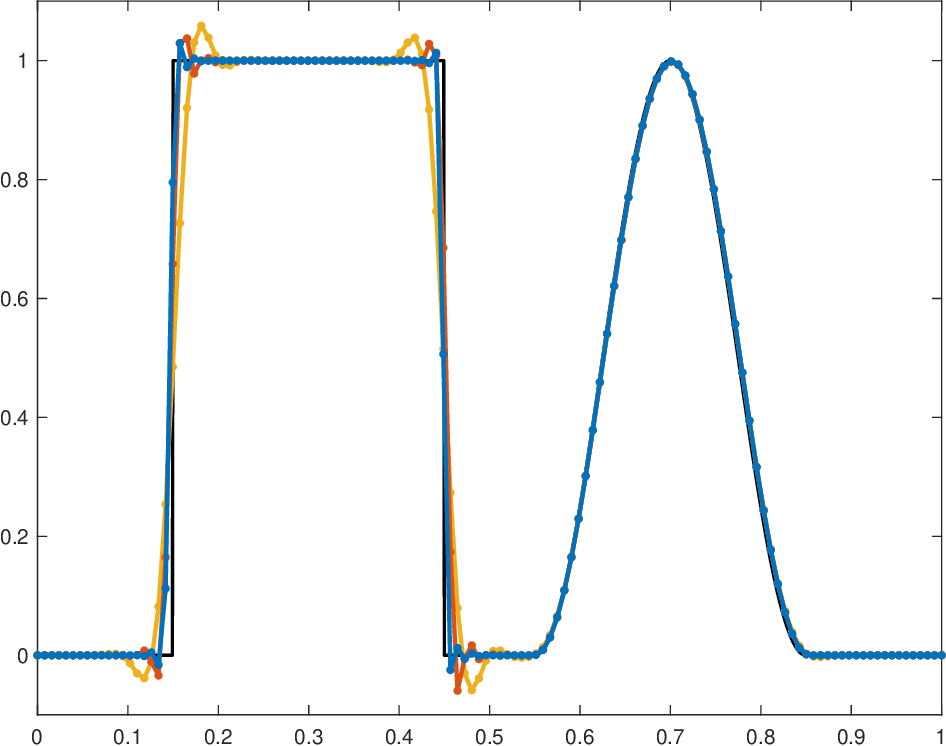}
		\caption{DG}
		\label{fig:advdg}
	\end{subfigure}
	\begin{subfigure}[b]{.32\linewidth}
		\includegraphics[width=\linewidth]{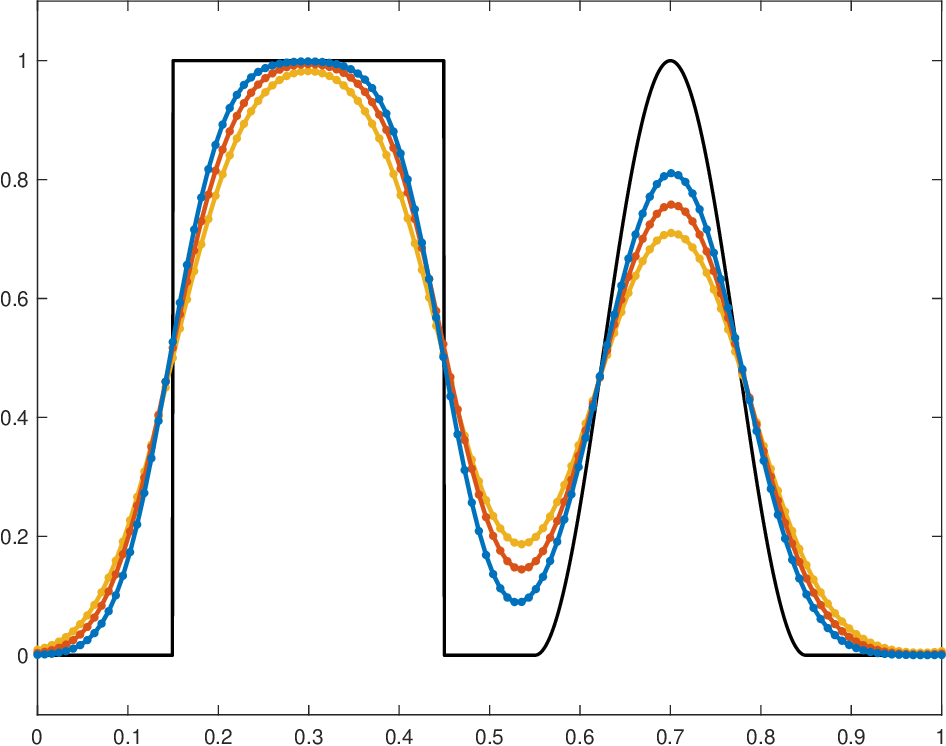}
		\caption{LO}
		\label{fig:advlo}
	\end{subfigure}
	\begin{subfigure}[b]{.32\linewidth}
		\includegraphics[width=\linewidth]{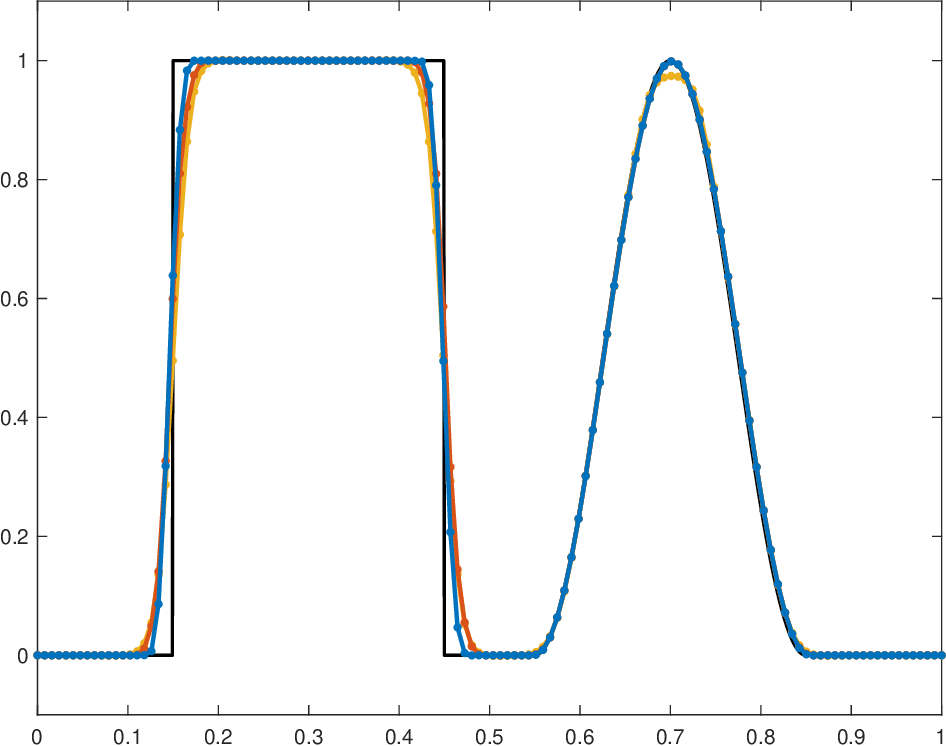}
		\caption{WENO}
		\label{fig:advweno}
	\end{subfigure}
	\caption{One-dimensional linear advection with constant velocity, numerical solutions at $t=1.0$ obtained using $E_h=128$ and $p\in\{1,2,3\}$.}
	\label{fig:adv}
\end{figure}

We run the numerical simulations up to the final time $t=1.0$ using $E_h=128$ elements. The results are shown in Figs \ref{fig:advdg}-\ref{fig:advweno}. As expected, the DG solutions display spurious oscillations near the discontinuities, while the solutions obtained with the LO method suffer from excessive numerical dissipation, resulting in reduced accuracy. The WENO solutions provide a sharp representation of the discontinuities without any apparent over- or undershoots. However, it is worth noting that for linear finite elements a peak clipping effect can be observed in the smooth profile.

%To investigate the long-term behavior of our scheme, we extend the simulation up to $t=100.0$, again employing $E_h=128$ elements. The corresponding results are presented in Figs \ref{fig:ladvdg}-\ref{fig:ladvweno}. The DG solution continues to exhibit spurious oscillations while the LO method produces even more dissipative results due to additional dissipation over time. Remarkably, the WENO scheme maintains its accuracy for long-term simulations.

% Linear advection: t = 1.0

%% Linear advection: t = 100.0
%\begin{figure}[!htb]
%	\centering
%	\begin{subfigure}[b]{\linewidth}
%		\centering
%		\begin{tikzpicture}
%		\draw[rounded corners] (0, 0) rectangle (12, 0.5) node[pos=.5]{};
%		\draw[very thick, color={rgb:red,0;green,0.4470;blue,0.7410}] (6.75,0.25)--(7.25,0.25);
%		\node at (2,0.25) (a) {$p=1$};
%		\node at (5,0.25) (a) {$p=2$};
%		\node at (8,0.25) (a) {$p=3$};
%		\node at (11,0.25) (a) {Exact};
%		\draw[very thick,color={rgb:red,0.9290;green,0.6940;blue,0.1250}] (0.75,0.25)--(1.25,0.25);
%		\draw[very thick,color={rgb:red,0.8500;green,0.3250;blue,0.0980}] (3.75,0.25)--(4.25,0.25);
%		\draw[very thick] (9.75,0.25)--(10.25,0.25);
%		\end{tikzpicture}
%		\vspace*{0.25cm}
%	\end{subfigure}
%	\begin{subfigure}[b]{.32\linewidth}
%		\includegraphics[width=\linewidth]{results/adv/advdg.eps}
%		\caption{DG}
%		\label{fig:ladvdg}
%	\end{subfigure}
%	\begin{subfigure}[b]{.32\linewidth}
%		\includegraphics[width=\linewidth]{results/adv/advlo.eps}
%		\caption{LO}
%		\label{fig:ladvlo}
%	\end{subfigure}
%	\begin{subfigure}[b]{.32\linewidth}
%		\includegraphics[width=\linewidth]{results/adv/advweno.eps}
%		\caption{WENO}
%		\label{fig:ladvweno}
%	\end{subfigure}
%	\caption{One-dimensional linear advection with constant velocity, numerical solutions at $t=100.0$ obtained using $E_h=128$ and $p\in\{1,2,3\}$.}
%	\label{fig:ladv}
%\end{figure}

\subsection{One-dimensional inviscid Burgers equation}
The first nonlinear problem we consider is the one-dimensional inviscid Burgers equation
\begin{equation}
\frac{\partial u}{\partial t}+\frac{\partial(u^2/2)}{\partial x} = 0 \quad \text{in } \Omega=(0,1).
\end{equation}
The initial condition
\begin{equation}
u_0(x)=\sin(2\pi x)
\end{equation}
remains smooth up to the critical time $t_c=\frac{1}{2\pi}$ of shock formation. Thus, we perform grid convergence studies at  the time $t=0.1$ for which the exact solution remains sufficiently smooth. The results are shown in Table \ref{tab:convburgers}. Similarly to the case of linear advection, both the DG and the WENO scheme yield optimal convergence rates. Once more, the WENO error converges to the Galerkin error. The numerical solutions obtained with the LO method maintain first-order accuracy in the nonlinear case.

\begin{table}[h!]
	\begin{subtable}{\textwidth}
		\centering
		\begin{tabular}{ccccccc}
			\hline
			&\multicolumn{2}{c}{DG}&\multicolumn{2}{c}{LO}&\multicolumn{2}{c}{WENO}\\
			\hline
			$E_h$ & $\|u_h-u_{\text{exact}}\|_{L^1}$ & EOC & $\|u_h-u_{\text{exact}}\|_{L^1}$ & EOC & $\|u_h-u_{\text{exact}}\|_{L^1}$ & EOC\\
			\hline
			16		 &   8.50e-03 & --       & 5.71e-02 & --       & 1.56e-02 &	--       \\
			32		 &   2.06e-03 &	2.05 & 2.97e-02 & 0.94 & 3.62e-03 & 2.11 \\
			64		 &   5.17e-04 & 1.99 & 1.54e-02 & 0.94 & 9.07e-04 & 2.00 \\
			128	   &   1.32e-04 & 1.97 & 7.98e-03 & 0.95 & 1.41e-04 & 2.68 \\
			256	   &   3.36e-05 & 1.98 & 4.09e-03 & 0.97 & 3.36e-05 & 2.07 \\
			512	   &   8.49e-06 & 1.98 & 2.08e-03 & 0.97 & 8.49e-06 & 1.98 \\
			1024 &   2.14e-06 & 1.99 & 1.05e-03 & 0.98 & 2.14e-06 & 1.99 \\
			2048 &   5.36e-07 & 1.99 & 5.29e-04 & 0.99 & 5.36e-07 & 1.99 \\
			\hline
		\end{tabular}
		\caption{$p=1$}
	\end{subtable}
	\begin{subtable}{\textwidth}
		\centering
		\begin{tabular}{ccccccc}
			\hline
			&\multicolumn{2}{c}{DG}&\multicolumn{2}{c}{LO}&\multicolumn{2}{c}{WENO}\\
			\hline
			$N_h$ & $\|u_h-u_{\text{exact}}\|_{L^1}$ & EOC & $\|u_h-u_{\text{exact}}\|_{L^1}$ & EOC & $\|u_h-u_{\text{exact}}\|_{L^1}$ & EOC\\
			\hline
			16      & 3.52e-04 & --       & 4.30e-02 & --       & 3.44e-04 & --       \\
			32      & 9.34e-05 & 1.92 & 2.24e-02 & 0.94 & 9.33e-05 & 1.88 \\
			64      & 1.21e-05 & 2.95 & 1.17e-02 & 0.94 & 1.21e-05 & 2.95 \\
			128   & 1.54e-06 & 2.97 & 6.11e-03 & 0.94 & 1.54e-06 & 2.97 \\
			256   & 1.95e-07 & 2.98 & 3.16e-03 & 0.95 & 1.95e-07 & 2.98 \\
			512   & 2.47e-08 & 2.98 & 1.62e-03 & 0.97 &	2.47e-08 & 2.98 \\
			1024 & 3.11e-09 & 2.99 & 8.21e-04 & 0.98 & 3.11e-09 & 2.99 \\
			2048 & 3.91e-10 & 2.99 & 4.14e-04 & 0.99 & 3.91e-10 & 2.99 \\
			\hline
		\end{tabular}
		\caption{$p=2$}
	\end{subtable}
	\begin{subtable}{\textwidth}
		\centering
		\begin{tabular}{ccccccc}
			\hline
			&\multicolumn{2}{c}{DG}&\multicolumn{2}{c}{LO}&\multicolumn{2}{c}{WENO}\\
			\hline
			$N_h$ & $\|u_h-u_{\text{exact}}\|_{L^1}$ & EOC & $\|u_h-u_{\text{exact}}\|_{L^1}$ & EOC & $\|u_h-u_{\text{exact}}\|_{L^1}$ & EOC\\
			\hline
			16    & 1.28e-04 & --       & 3.43e-02 & --       & 1.27e-04 & --       \\
			32    & 8.65e-06 & 3.89 & 1.76e-02 & 0.96 & 8.66e-06 & 3.87 \\
			64    & 5.10e-07 & 4.08 & 9.18e-03 & 0.94 & 5.11e-07 & 4.08 \\
			128 & 3.22e-08 & 3.99 & 4.76e-03 & 0.95 & 3.22e-08 & 3.99 \\
			256 & 2.04e-09 & 3.98 & 2.46e-03 & 0.95 & 2.04e-09 & 3.98 \\
			512 & 1.30e-10 & 3.98 & 1.26e-03 & 0.96 & 1.30e-10 & 3.98 \\
			\hline
		\end{tabular}
		\caption{$p=3$}
	\end{subtable}
	\caption{One-dimensional Burgers equation, grid convergence history for finite elements of degree $p\in\{1,2,3\}$.}
	\label{tab:convburgers}
\end{table}

\begin{figure}[!htb]
	\centering
	\begin{subfigure}[b]{\linewidth}
		\centering
		\begin{tikzpicture}
		\draw[rounded corners] (0, 0) rectangle (9.25, 0.5) node[pos=.5]{};
		\draw[very thick, color={rgb:red,0;green,0.4470;blue,0.7410}] (6.75,0.25)--(7.25,0.25);
		\node at (2,0.25) (a) {$p=1$};
		\node at (5,0.25) (a) {$p=2$};
		\node at (8,0.25) (a) {$p=3$};
		\draw[very thick,color={rgb:red,0.9290;green,0.6940;blue,0.1250}] (0.75,0.25)--(1.25,0.25);
		\draw[very thick,color={rgb:red,0.8500;green,0.3250;blue,0.0980}] (3.75,0.25)--(4.25,0.25);
		\end{tikzpicture}
		\vspace*{0.25cm}
	\end{subfigure}
	\begin{subfigure}[b]{.48\linewidth}
		\includegraphics[width=\linewidth]{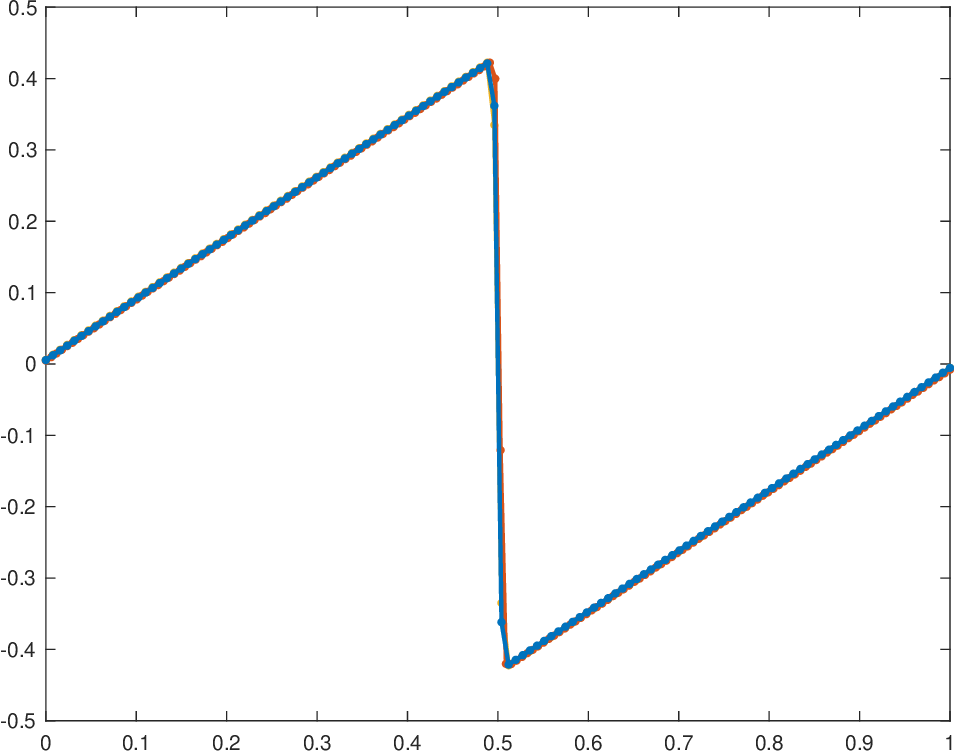}
		\caption{LO}
		\label{fig:burgerslo}
	\end{subfigure}
	\begin{subfigure}[b]{.48\linewidth}
		\includegraphics[width=\linewidth]{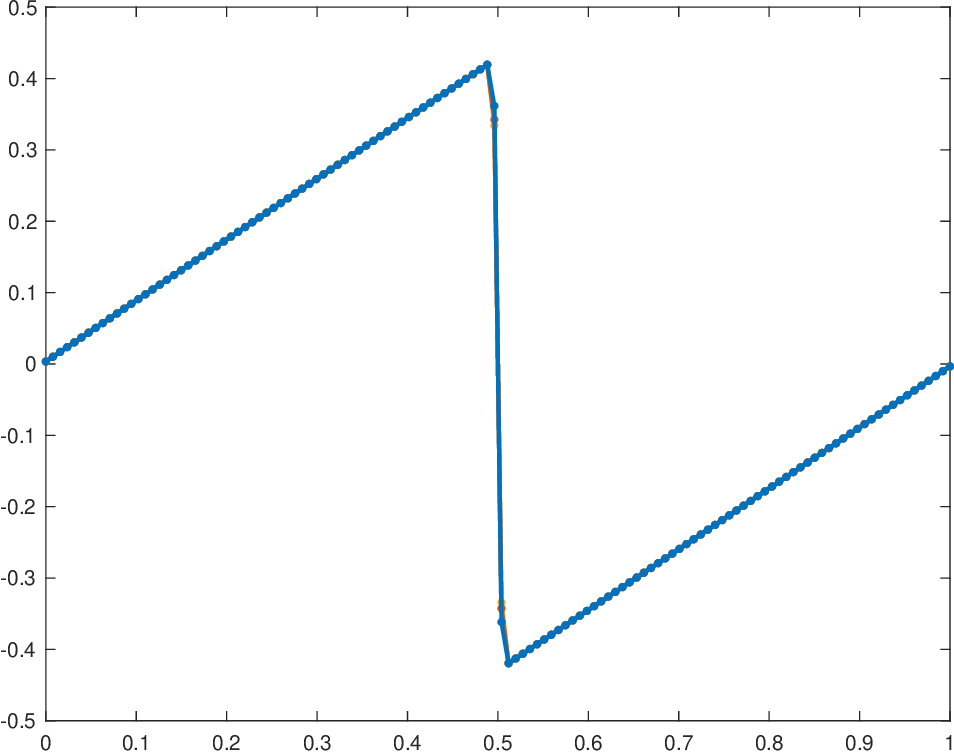}
		\caption{WENO}
		\label{fig:burgersweno}
	\end{subfigure}
	\caption{One-dimensional Burgers equation, numerical solutions at $t=1.0$ obtained using $E_h=128$ and $p\in\{1,2,3\}$.}
	\label{fig:burgers}
\end{figure}

To investigate the shock behavior for nonlinear problems, we extend the final time to $t=1.0$ after the shock has formed. For our numerical simulations, we employ finite elements of degree $p\in\{1,2,3\}$ on $E_h=128$ elements. The results obtained with the LO scheme and the WENO scheme are shown in Figs \ref{fig:burgerslo} and \ref{fig:burgersweno}, respectively. Due to the self-steepening nature of shocks, the LO scheme exhibits reduced diffusion. Interestingly, there are no observable differences between the solutions obtained using different polynomial degrees. This finding is consistent with the results reported in \cite{kuzmin2023a}.

% Burgers

\subsection{Two-dimensional solid body rotation}
Next, we consider LeVeque's solid body rotation benchmark \cite{leveque1996}. We solve the transport problem
\begin{equation}
\frac{\partial u}{\partial t}+\nabla\cdot(\mathbf{v}u) = 0 \quad \text{in } \Omega=(0,1)^2
\end{equation}
with the divergence-free velocity field ${\bf v}(x,y)=2\pi(0.5-y,x-0.5)$. In this two-dimensional test problem, a smooth hump, a sharp cone and a slotted cylinder are rotated around the center of the domain. After each revolution ($t=2\pi r$, $r\in \mathbb{N}$) the exact solution corresponds to the initial condition given by
\begin{equation}
u_0(x,y) = \begin{cases}
u_0^{\text{hump}}(x,y) & \text{if } \sqrt{(x-0.25)^2+(y-0.5)^2}\le 0.15, \\
u_0^{\text{cone}}(x,y) & \text{if } \sqrt{(x-0.5)^2+(y-0.25)^2}\le 0.15, \\
1 & \text{if }\begin{cases}
(\sqrt{(x-0.5)^2+(y-0.75)^2}\le0.15) \wedge \\
(|x-0.5|\ge0.025\vee y\ge0.85),
\end{cases}\\
0 & \text{otherwise},
\end{cases}
\end{equation}
where
\begin{align*}
u_0^{\text{hump}}(x,y) &= \frac{1}{4}+\frac{1}{4}\cos\bigg(\frac{\pi\sqrt{(x-0.25)^2+(y-0.5)^2}}{0.15}\bigg), \\
u_0^{\text{cone}}(x,y) &= 1-\frac{\sqrt{(x-0.5)^2+(y-0.25)^2}}{0.15}.
\end{align*}

\begin{figure}[!htb]
	\centering
	\begin{subfigure}[b]{.32\linewidth}
		\includegraphics[width=\linewidth]{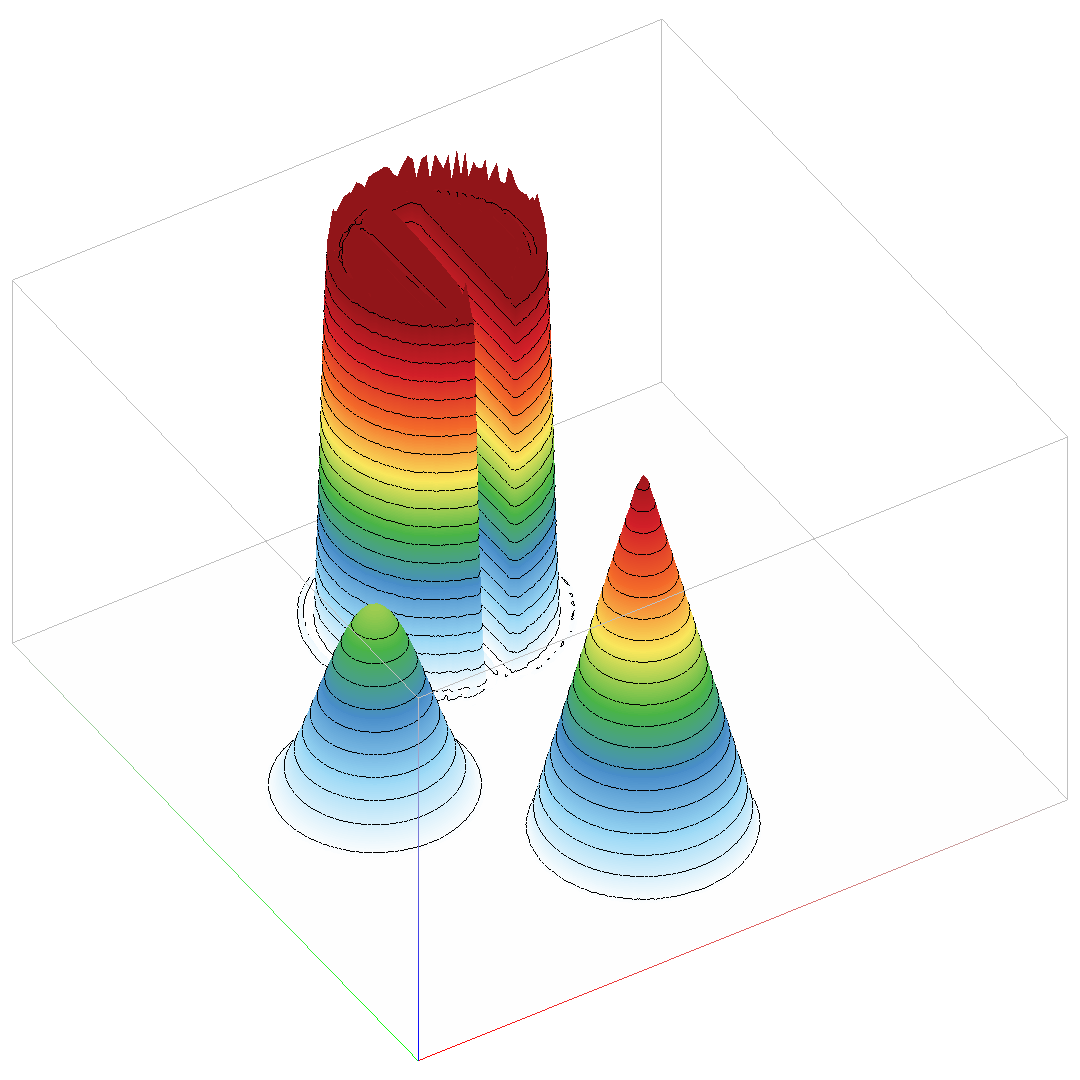}
	\end{subfigure}
	\begin{subfigure}[b]{.32\linewidth}
		\includegraphics[width=\linewidth]{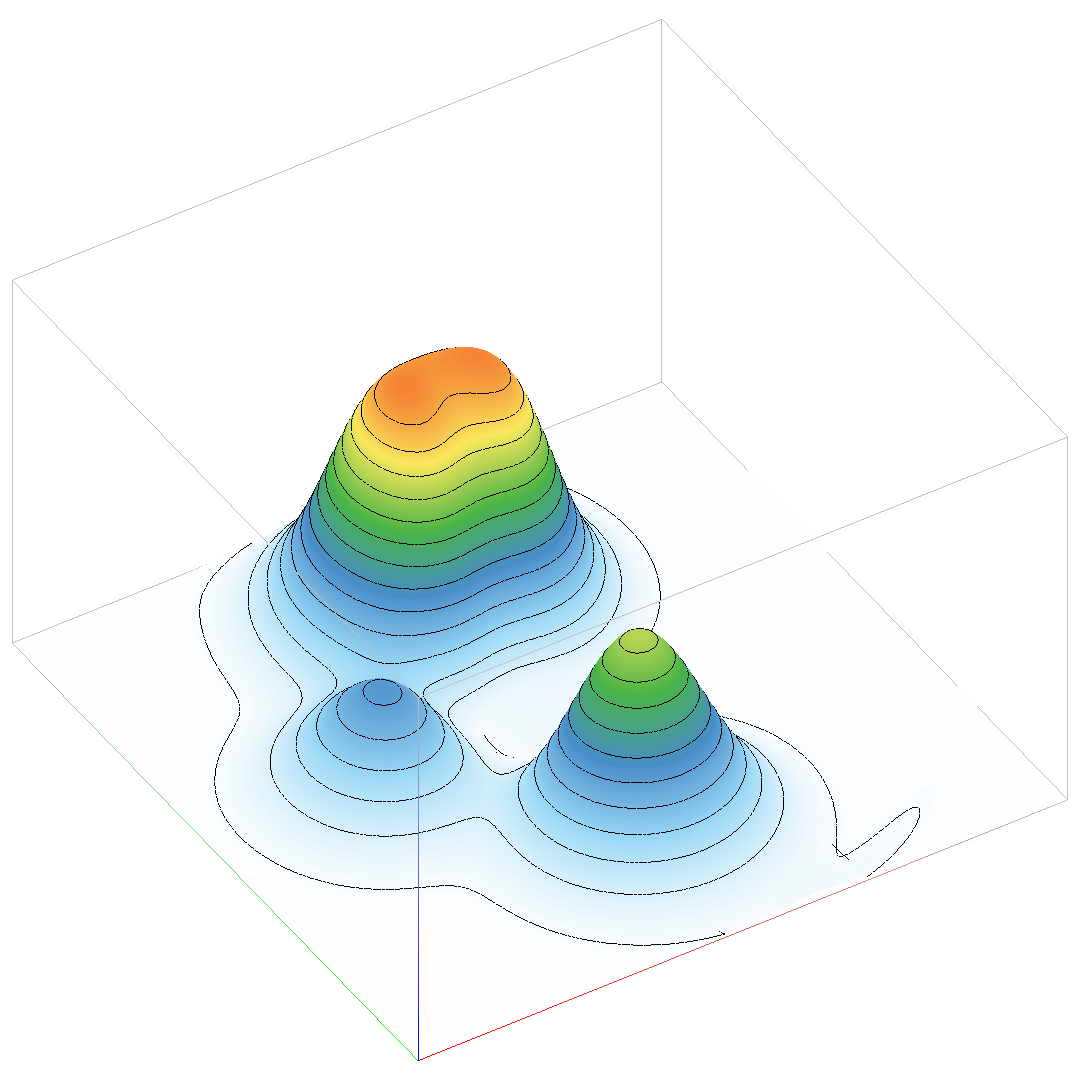}
	\end{subfigure}
	\begin{subfigure}[b]{.32\linewidth}
		\includegraphics[width=\linewidth]{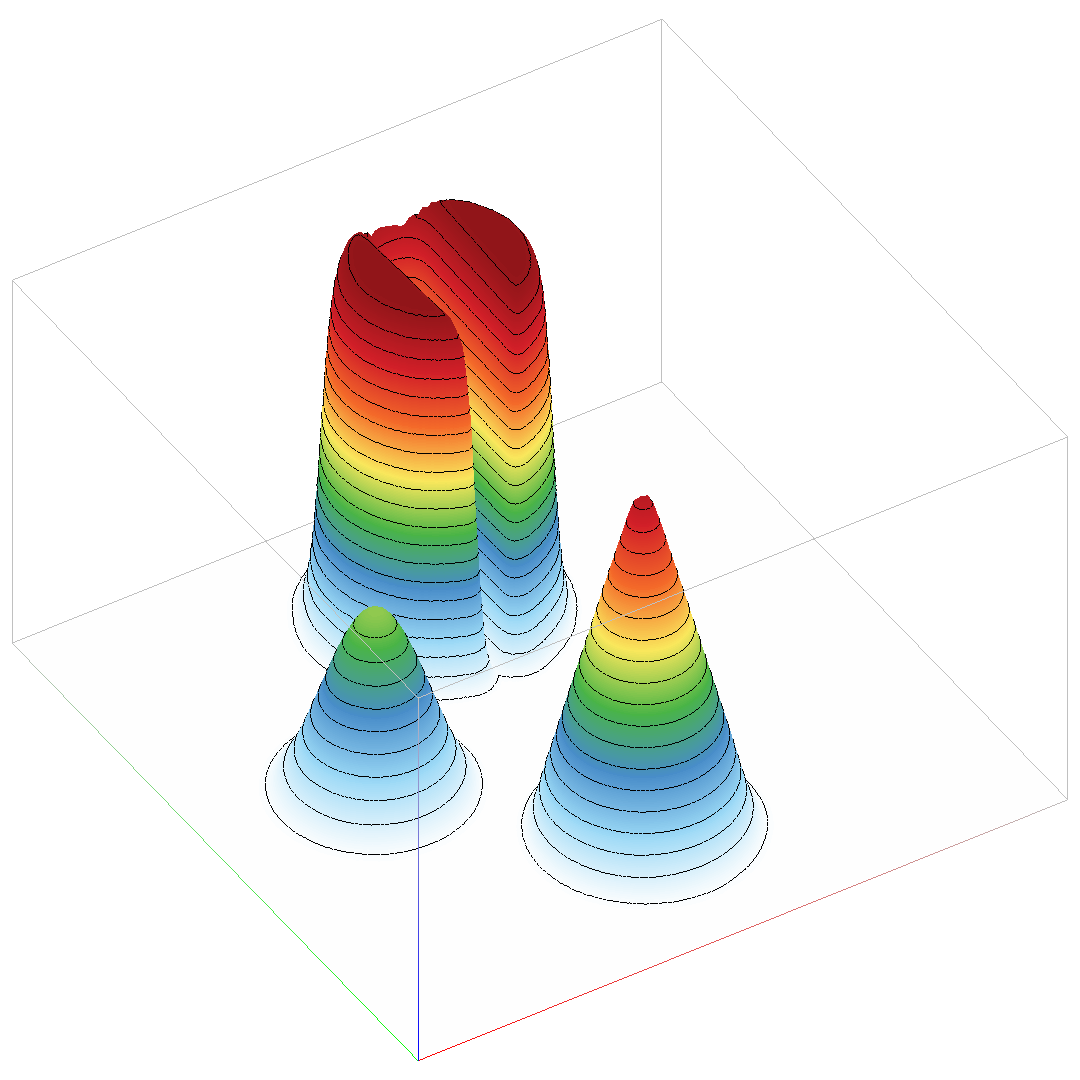}
	\end{subfigure}
	\begin{subfigure}[b]{.32\linewidth}
		\includegraphics[width=\linewidth]{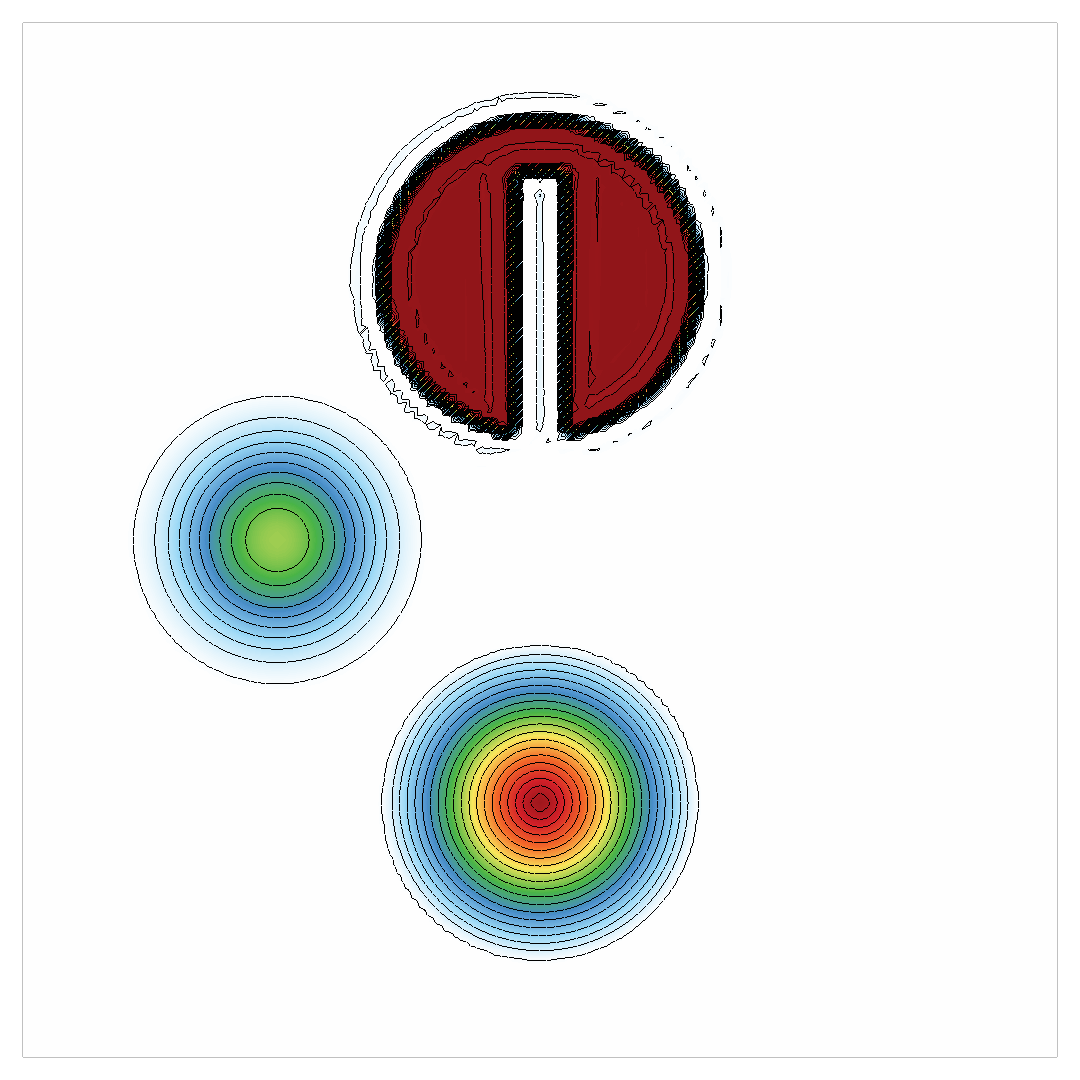}
		\caption{DG, $u_h \in [-0.107,1.120]$}
		\label{fig:sbrdg}
	\end{subfigure}
	\begin{subfigure}[b]{.32\linewidth}
		\includegraphics[width=\linewidth]{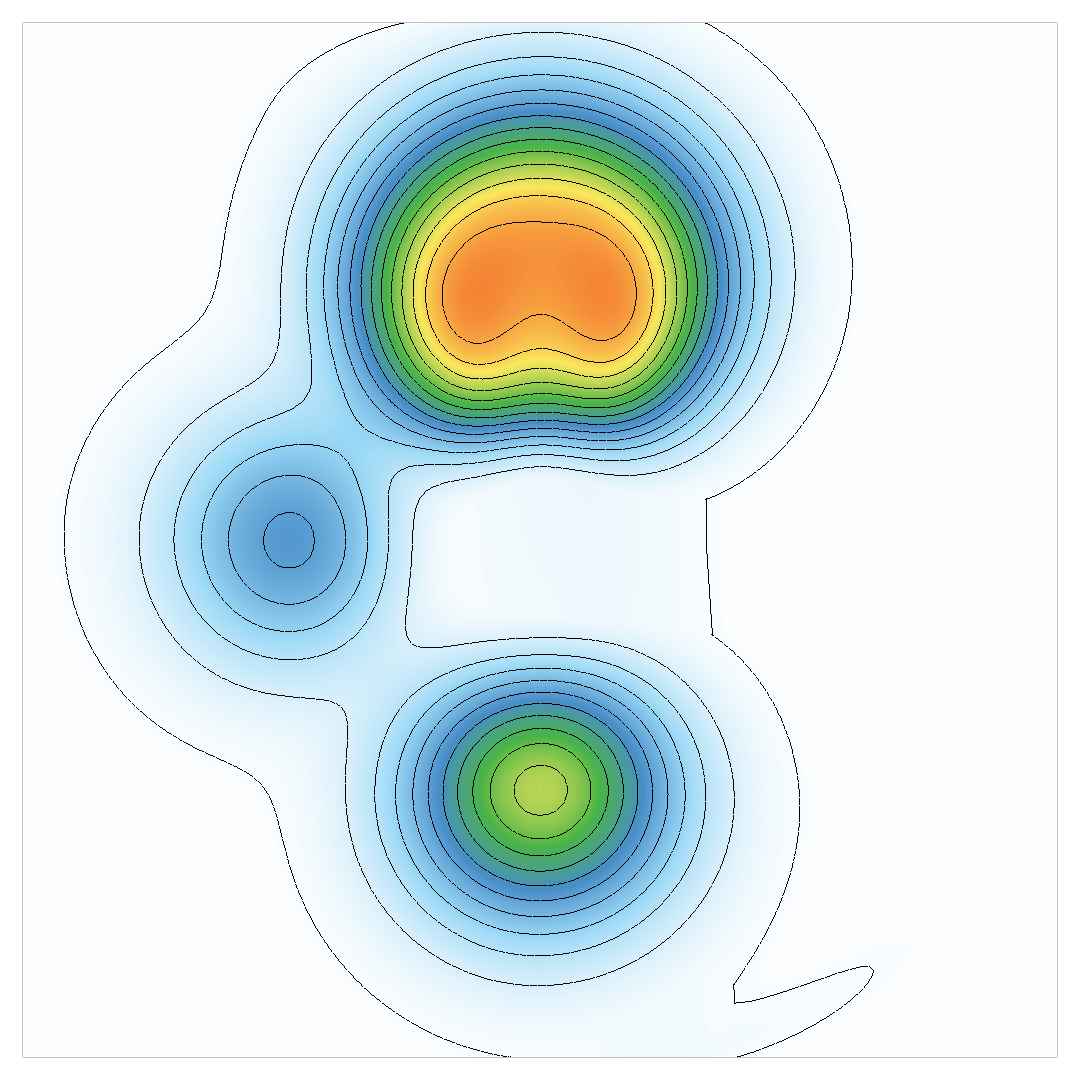}
		\caption{LO, $u_h \in [0.000,0.690]$}
		\label{fig:sbrlo}
	\end{subfigure}
	\begin{subfigure}[b]{.32\linewidth}
		\includegraphics[width=\linewidth]{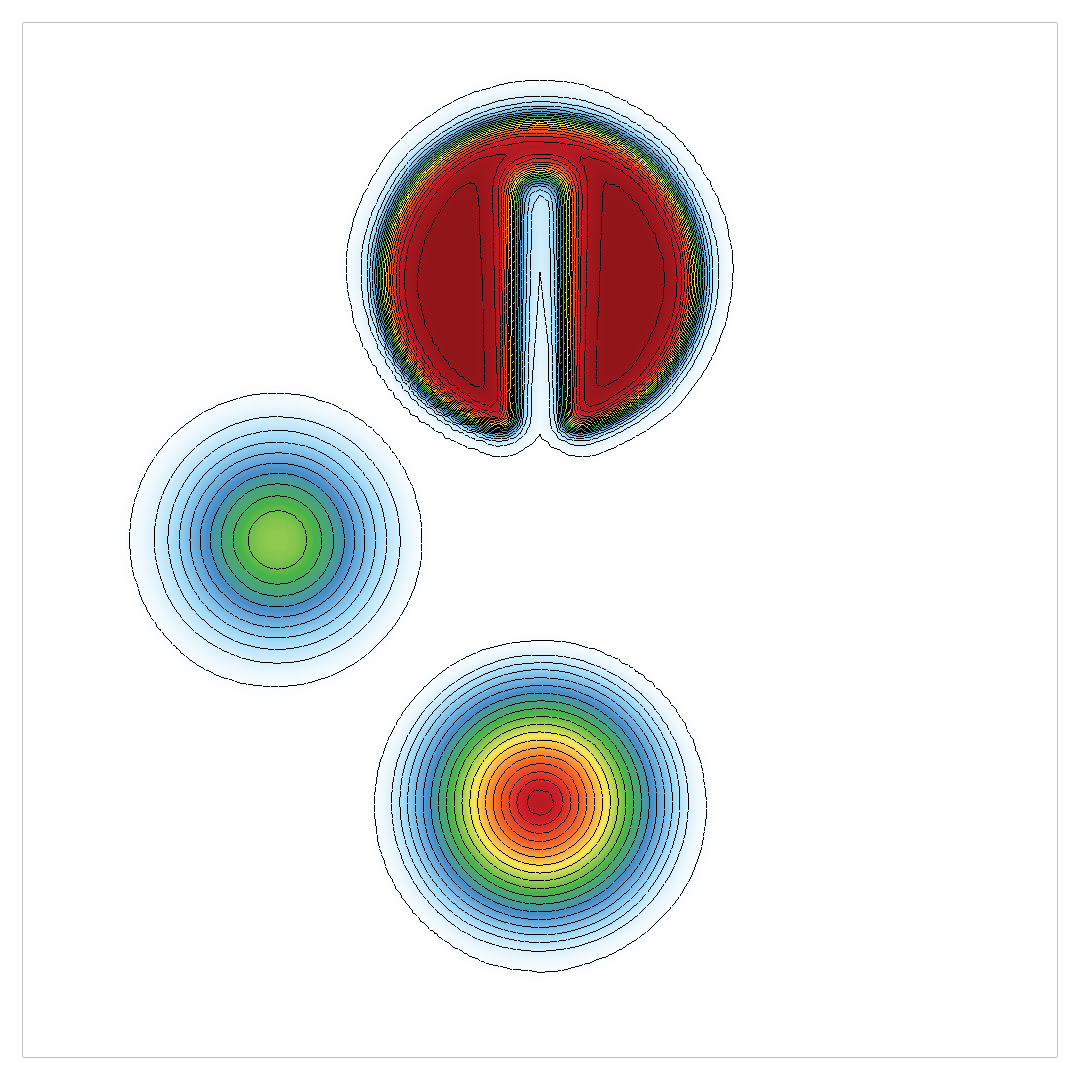}
		\caption{WENO, $u_h \in [0.000,1.000]$}
		\label{fig:sbrweno}
	\end{subfigure}
	\caption{Solid body rotation, numerical solutions at $t=1.0$ obtained using $E_h=128^2$ and $p=2$.}
	\label{fig:sbr}
\end{figure}

We perform numerical computations up to the final time $t=1.0$ on a uniform quadrilateral mesh using $E_h=128^2$ elements and $p=2$. Similarly to the one-dimensional scalar test problems, the DG solution, as presented in Fig. \ref{fig:sbrdg}, exhibits over- and undershoots near discontinuities, while solutions obtained with the LO scheme suffer from significant numerical dissipation; see Fig. \ref{fig:sbrlo}. In fact, the LO scheme fails to reproduce the geometric features of the initial condition accurately. Contrary to that, the WENO scheme effectively suppresses oscillations near discontinuities and accurately preserves the structure of all rotating objects; see Fig. \ref{fig:sbrweno}.

% Solid body rotation

\subsection{Two-dimensional Burgers equation}
Next, we consider the 2D inviscid Burgers equation 
\begin{align}
	\frac{\partial u}{\partial t}+\nabla\cdot \Big(\mathbf{v}\frac{u^2}{2}\Big)= 0 \quad \text{in } \Omega=(0,1)^2,
\end{align}
where $\mathbf{v}=(1,1)^\top$ is a constant vector. The initial condition is given by 
\begin{equation}
u_0(x,y)=\begin{cases}
-0.2 & \text{if }x<0.5 \wedge y>0.5,\\
-1.0 & \text{if }x>0.5 \wedge y>0.5,\\
\phantom{-}0.5 & \text{if }x<0.5 \wedge y<0.5,\\
\phantom{-}0.8 & \text{if }x>0.5 \wedge y<0.5.
\end{cases}
\end{equation}
The inflow boundary conditions are defined using the exact solution of the initial value problem in $\mathbb{R}^2$, as detailed in \cite{guermond2014a}. This solution stays in the invariant $\mathcal{G}=[-1.0,0.8]$.

\begin{figure}[!htb]
	\centering
	\begin{subfigure}[b]{.32\linewidth}
		\includegraphics[width=\linewidth]{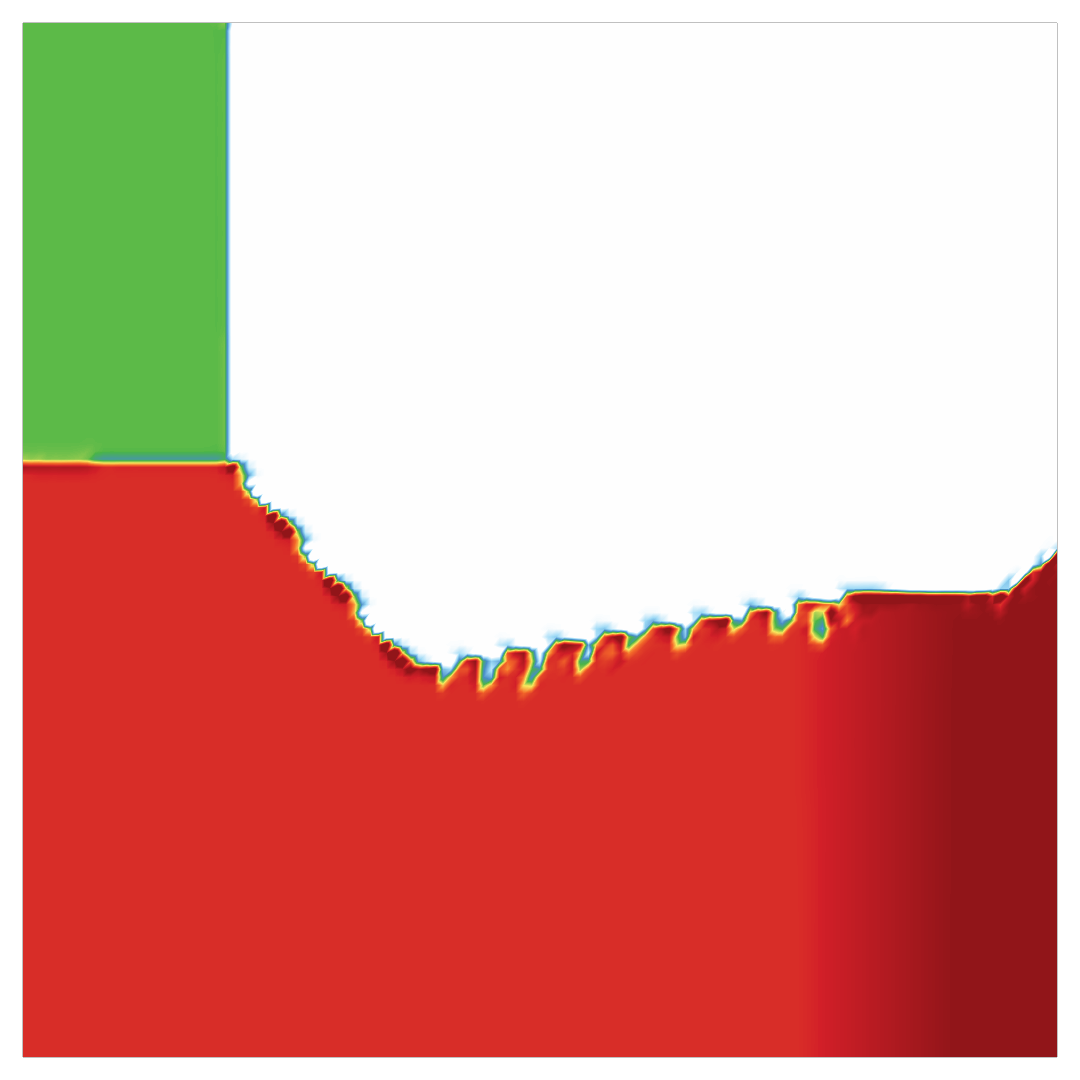}
		\caption{DG, $u_h \in [-2.091,1.691]$}
		\label{fig:burdg}
	\end{subfigure}
	\begin{subfigure}[b]{.32\linewidth}
		\includegraphics[width=\linewidth]{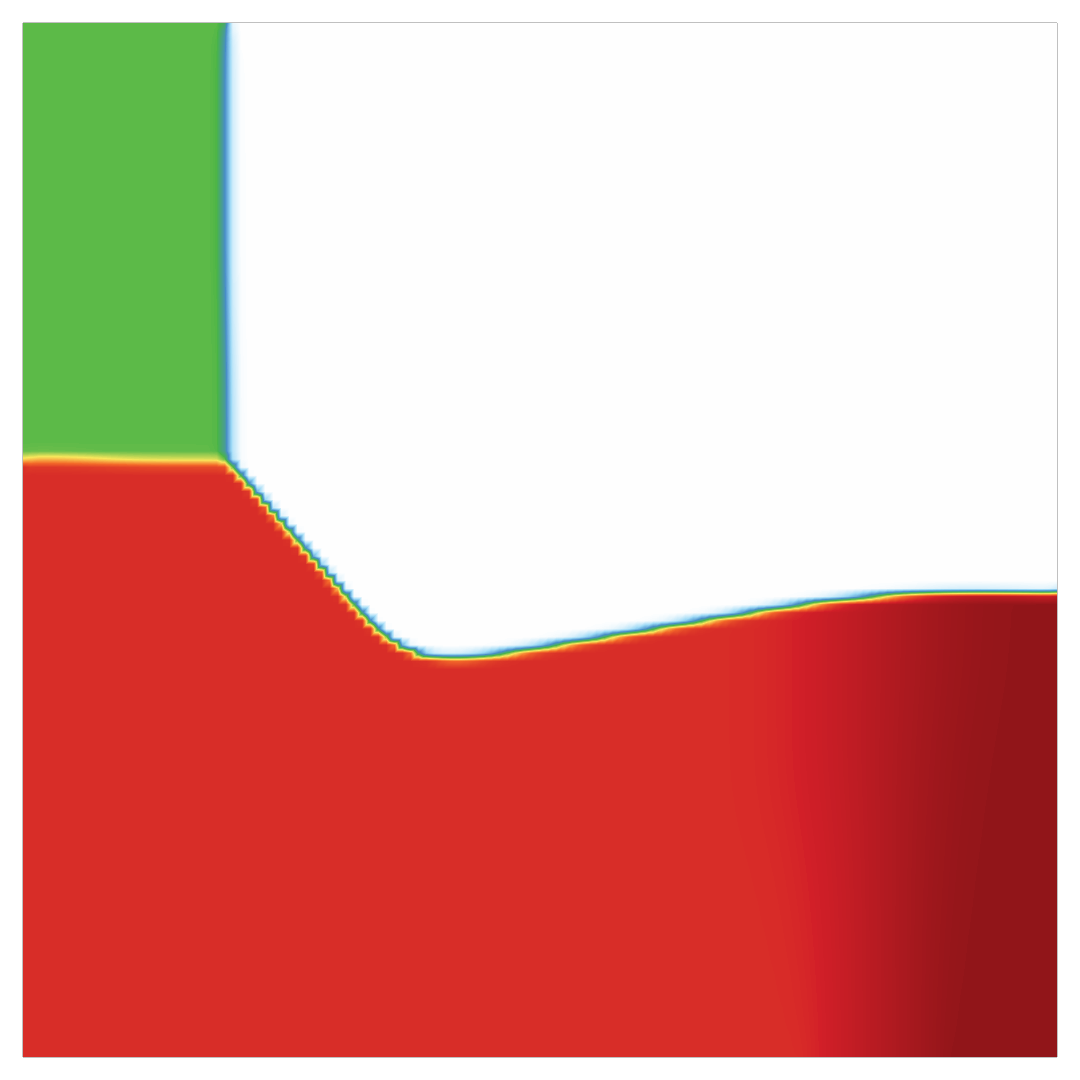}
		\caption{LO, $u_h \in [-1.000,0.800]$}
		\label{fig:burlo}
	\end{subfigure}
	\begin{subfigure}[b]{.32\linewidth}
		\includegraphics[width=\linewidth]{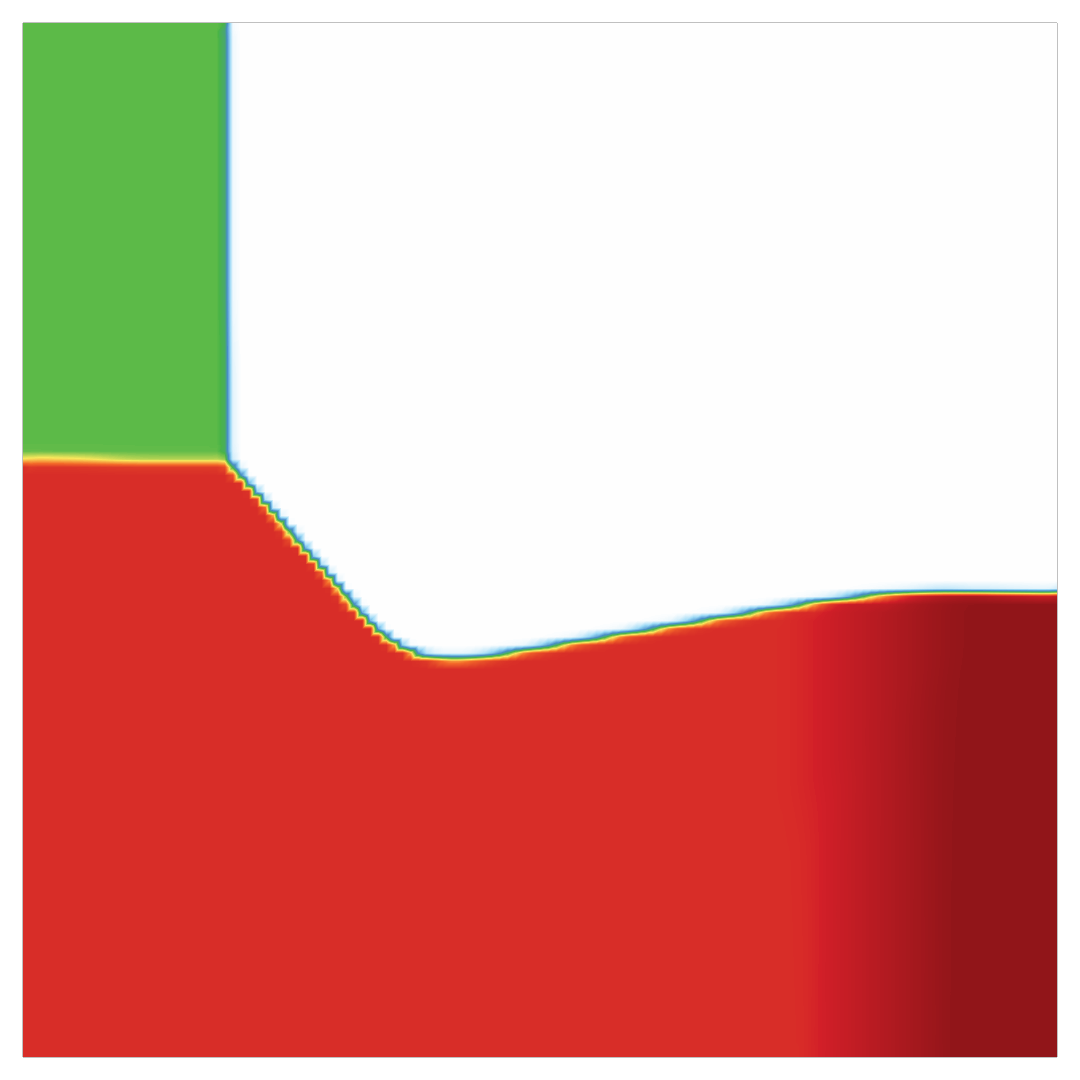}
		\caption{WENO,$u_h \in [-1.000,0.800]$}
		\label{fig:burweno}
	\end{subfigure}
	\caption{Two-dimensional Burgers equation, density profiles at $t=0.5$ obtained using $E_h=128^2$ and $p=2$.}
	\label{fig:bur}
\end{figure}

The numerical solutions obtained at the final time $t=0.5$ using $E_h=128^2$ elements and $p=2$ are shown in Fig. \ref{fig:bur}. As expected, the DG scheme produces highly oscillatory solutions near shocks. Similar to the one-dimensional benchmark, the LO scheme exhibits reduced diffusion. Solutions produced by both the LO and the WENO scheme remain within the set $\mathcal{G}$.

% 2D Burgers equation

\subsection{Two-dimensional KPP problem}
We investigate the entropy stability properties of our scheme using the two-dimensional KPP problem \cite{kurganov2007}. In this example, we solve
\begin{equation}
\frac{\partial u}{\partial t}+\nabla \cdot \mathbf{f}(u) = 0 ,
\end{equation}
where the nonconvex flux function is given by
\begin{equation}
\mathbf{f}(u)=(\sin(u),\cos(u)).
\end{equation}
The computational domain is $\Omega=(-2,2)\times(-2.5,1.5)$. The main challenge of this test problem lies in the potential convergence to incorrect weak solutions instead of the entropy solution, which exhibits a rotational wave structure. The initial condition is given by
\begin{equation}
u_0(x,y)=\begin{cases}
\frac{7\pi}{2} & \text{if }\sqrt{x^2+y^2}\le 1,\\
\frac{\pi}{4} & \text{otherwise}.
\end{cases}
\end{equation}

Once again, we perform numerical simulations up to the final time $t=1.0$ on a uniform quadrilateral mesh using $E_h=128^2$ elements and $p=2$.  The global upper bound for the maximum speed required to compute the viscosity parameter $\nu_e$ in \eqref{eq:visc} is $\lambda_e=1.0$. More accurate bounds can be found in \cite{guermond2017}. To ensure the LO component of the WENO scheme is sufficiently dissipative to suppress undershoots/ overshoots, we choose $\lambda_e=2.0$. Unlike the WENO-based stabilization approach applied to continuous finite elements in \cite{kuzmin2023a}, we do not modify the WENO parameters for this particular test problem. 

\begin{figure}[!htb]
	\centering
	\begin{subfigure}[b]{.32\linewidth}
		\includegraphics[width=\linewidth]{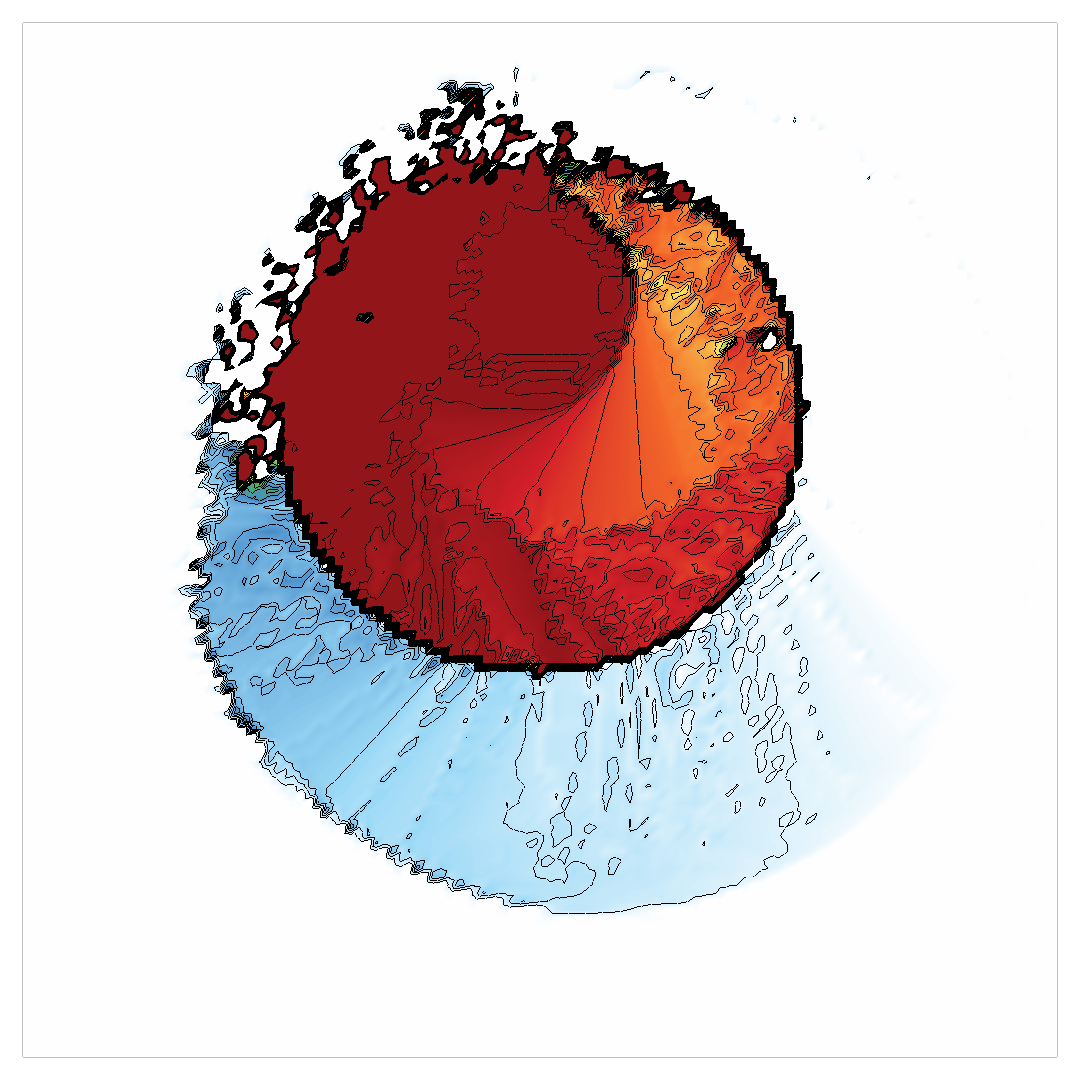}
		\caption{DG, $u_h \in [-44.944,43.235]$}
		\label{fig:kppdg}
	\end{subfigure}
	\begin{subfigure}[b]{.32\linewidth}
		\includegraphics[width=\linewidth]{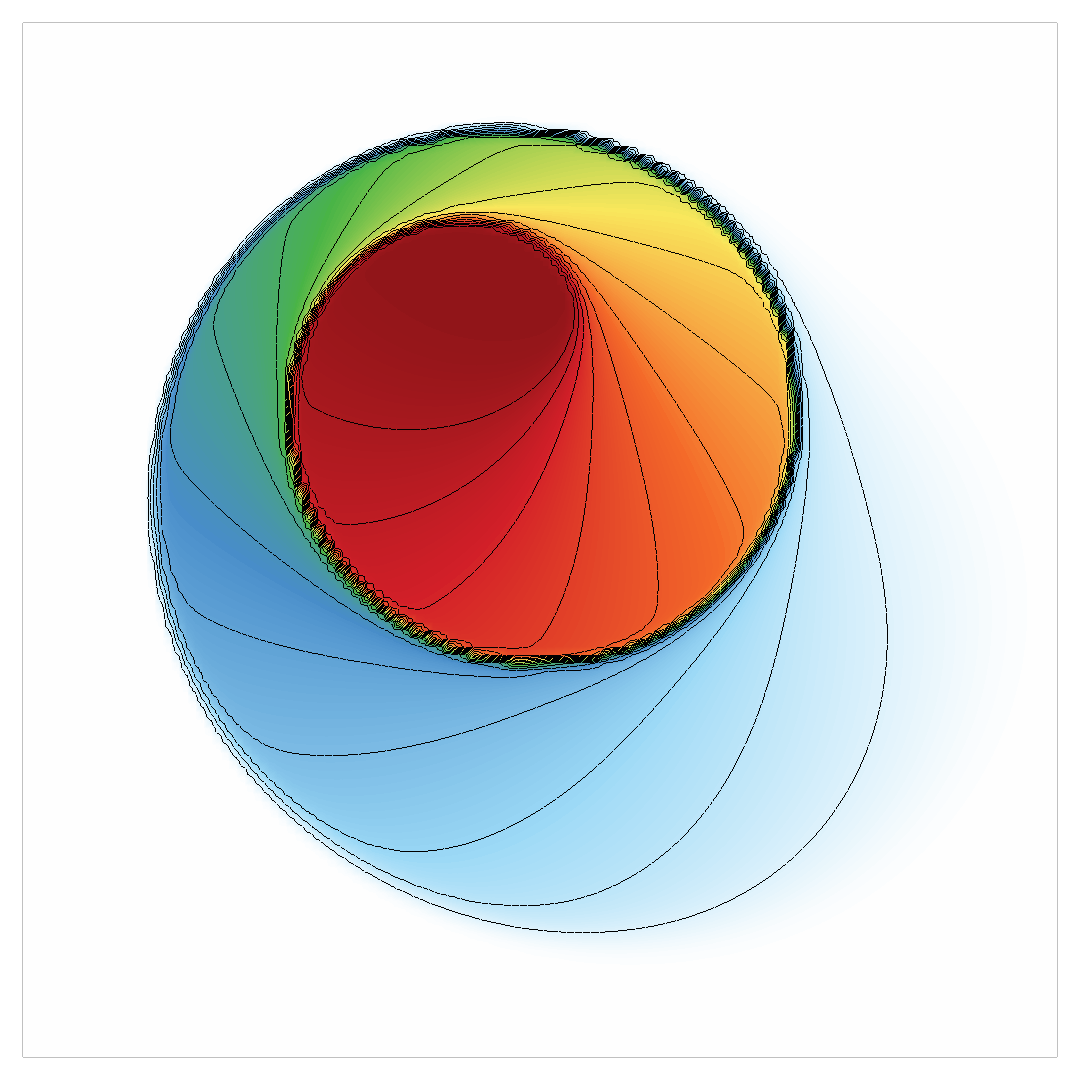}
		\caption{LO, $u_h \in [0.785,10.996]$}
		\label{fig:kpplo}
	\end{subfigure}
	\begin{subfigure}[b]{.32\linewidth}
		\includegraphics[width=\linewidth]{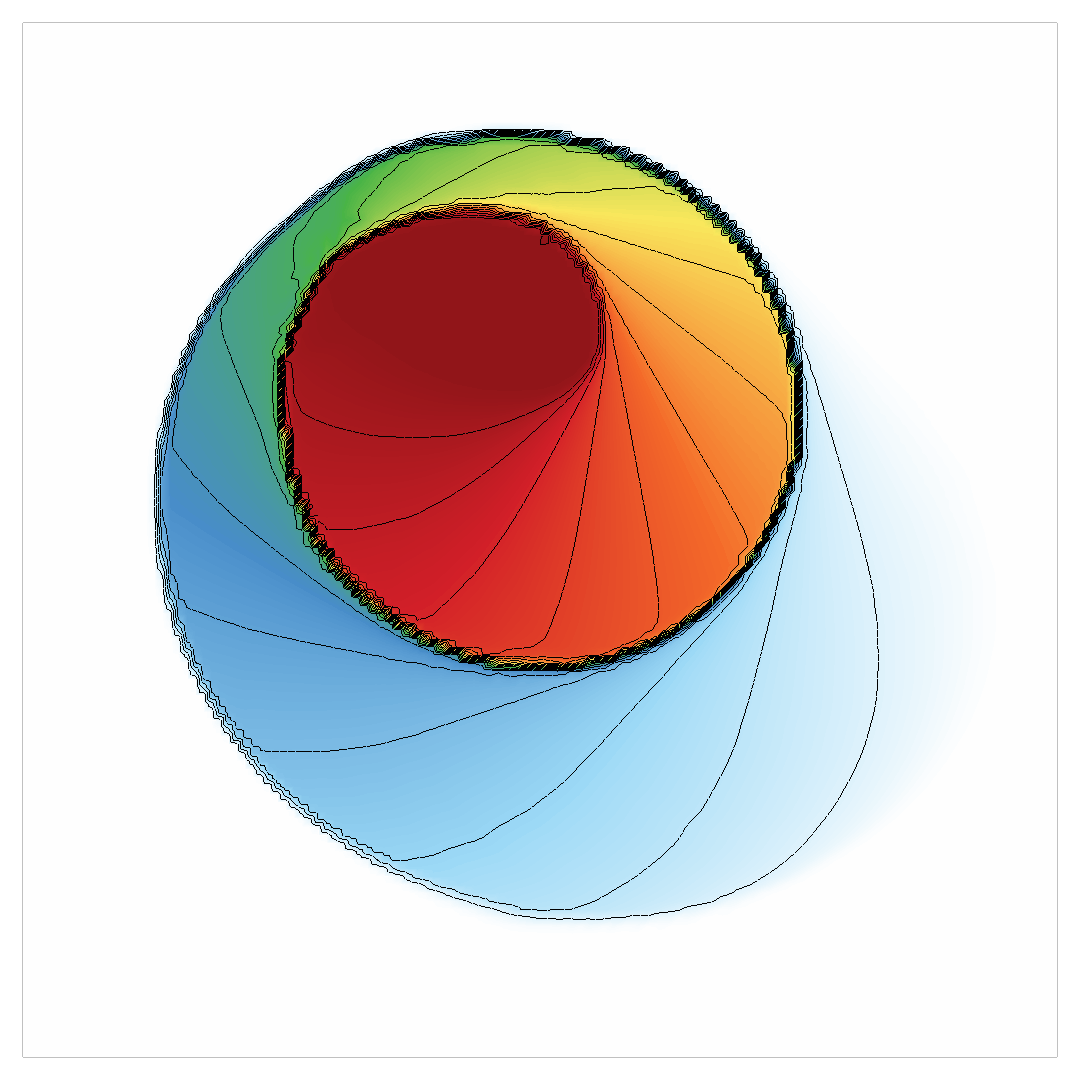}
		\caption{WENO, $u_h \in [0.785,11.066]$}
		\label{fig:kppweno}
	\end{subfigure}
	\caption{Two-dimensional KPP problem, numerical solutions at $t=1.0$ obtained using $E_h=128^2$ and $p=2$.}
	\label{fig:kpp}
\end{figure}

Fig. \ref{fig:kppdg} shows that the DG solution suffers from excessive oscillations and converges to an entropy-violating solution. Contrary to that, both the LO scheme and the WENO scheme converge to the correct entropy solution. The numerical solutions produced by these schemes are shown in Figs \ref{fig:kpplo} and \ref{fig:kppweno}, respectively. The WENO result is less diffusive, capturing shocks more sharply than the solution obtained using the LO scheme.

% KPP

\subsection{Euler equations of gas dynamics}
We consider the Euler equations of gas dynamics which represent the conservation of mass, momentum and total energy. The solution vector and flux matrix in \eqref{eq:pde} read
\begin{equation}
U = \begin{bmatrix}
\varrho \\
\varrho \mathbf{v}\\
\varrho E
\end{bmatrix}\in \mathbb{R}^{d+2}, \quad 
\mathbf{F(U)}=\begin{bmatrix}
\varrho \mathbf{v} \\
\varrho \mathbf{v} \bigotimes \mathbf{v}+p\mathbf{I}\\
(\varrho E + p)\mathbf{v}
\end{bmatrix}\in \mathbb{R}^{(d+2) \times d},
\end{equation}
where $\varrho$, $\mathbf{v}$, $E$ are the density, velocity and specific total energy, respectively. The identity matrix is denoted by $\mathbf{I}$, and the pressure $p$ is computed using the polytropic ideal gas assumption
\begin{equation}
p=\varrho e (\gamma-1).
\end{equation}
Here, $\varrho e$ and $\gamma$ are the internal energy and the heat capacity ratio, respectively. We use $\gamma=1.4$ in all of our numerical experiments.

Similarly to the approach presented in \cite{zhao2020}, our WENO-based shock sensor is designed to use information solely from the density field, resulting in significant computational savings \cite{pirozzoli2002}.
\subsubsection{Sod shock tube}

Sod's shock tube problem \cite{sod1978} serves as our first test for solving the Euler equations. This classical problem is commonly employed to assess the accuracy of high-order numerical methods. The computational domain $\Omega=(0,1)$ is delimited by reflecting walls and initially divided by a membrane into two distinct regions. Upon removing the membrane, a discontinuity at $x=0.5$ leads to the formation of a shock wave, rarefaction wave, and contact discontinuity. The initial condition is given by
\begin{equation}
\begin{bmatrix}
\varrho_L \\v_L\\p_L
\end{bmatrix}=
\begin{bmatrix}1.0\\0.0\\1.0
\end{bmatrix}, \quad 
\begin{bmatrix}
\varrho_R \\v_R\\p_R
\end{bmatrix}=
\begin{bmatrix}0.125\\0.0\\0.1
\end{bmatrix}.
\end{equation}

\begin{figure}[!htb]
	\centering
	\begin{subfigure}[b]{.32\linewidth}
		\includegraphics[width=\linewidth]{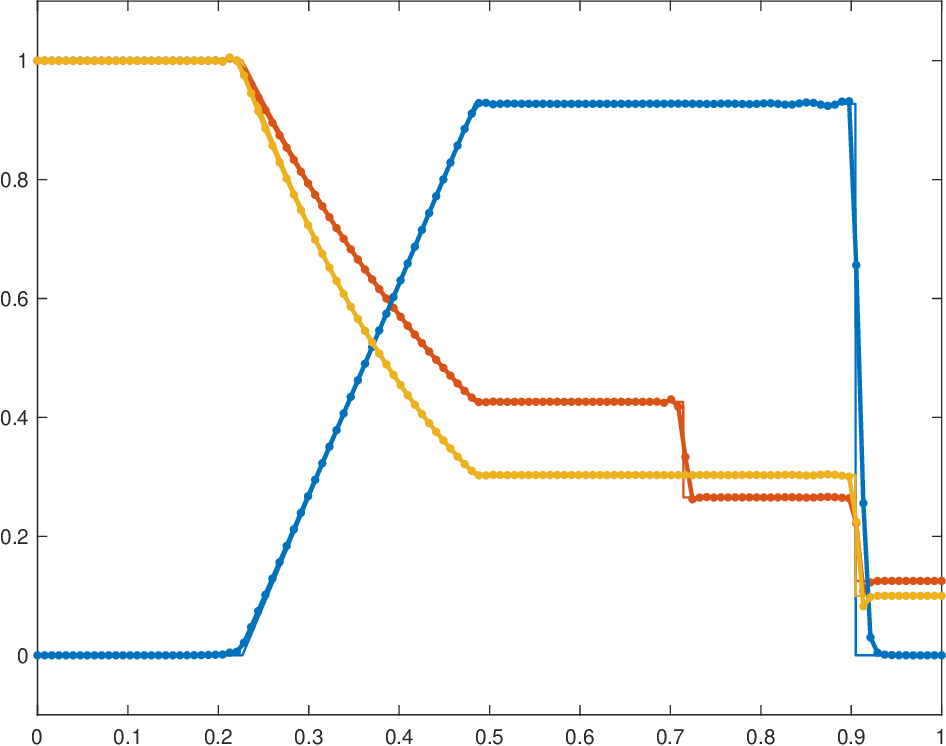}
		\caption{DG}
		\label{fig:soddg}
	\end{subfigure}
	\begin{subfigure}[b]{.32\linewidth}
		\includegraphics[width=\linewidth]{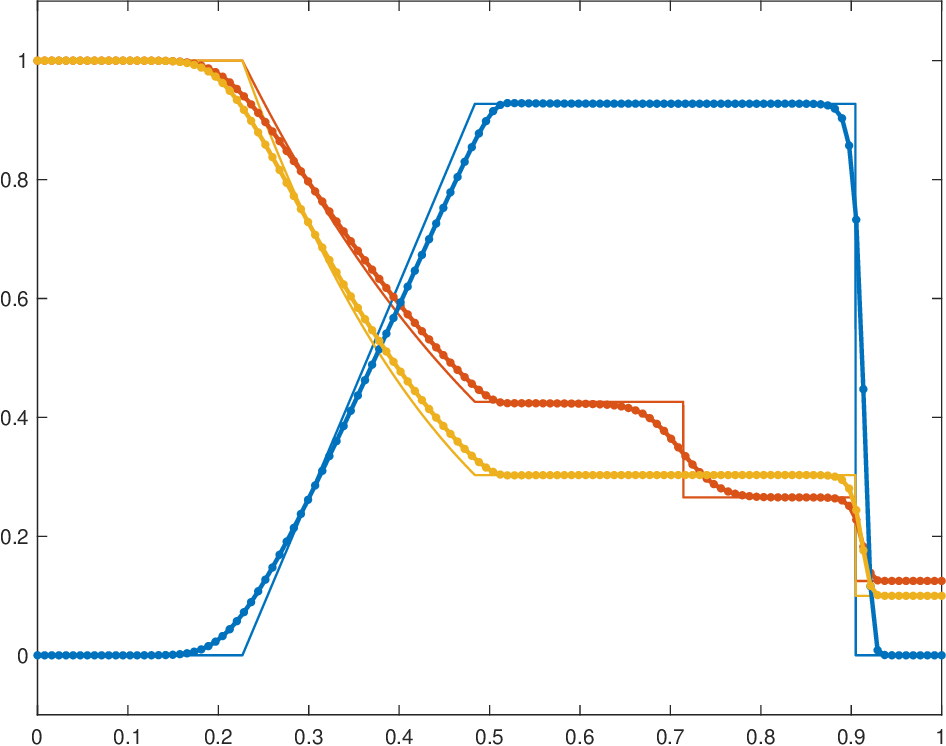}
		\caption{LO}
		\label{fig:sodlo}
	\end{subfigure}
	\begin{subfigure}[b]{.32\linewidth}
		\includegraphics[width=\linewidth]{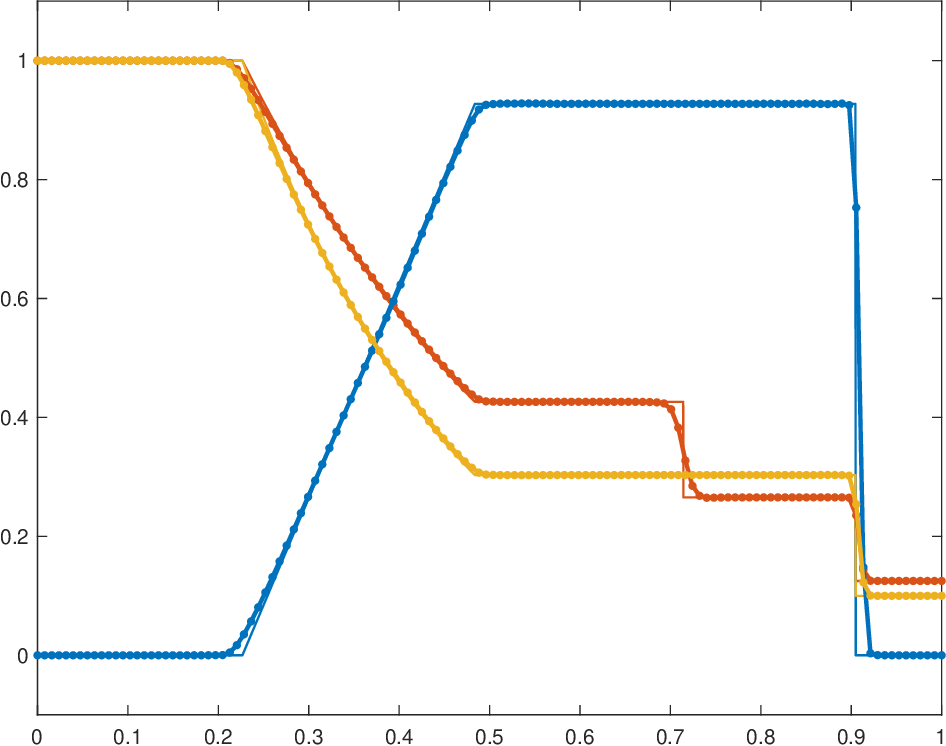}
		\caption{WENO}
		\label{fig:sodweno}
	\end{subfigure}
	\caption{Sod's shock tube, $(\varrho,v,E)$ at $t=0.231$ obtained using $E_h=128$ and $p=2$.}
	\label{fig:sod}
\end{figure}

We evolve numerical solutions up to the final time $t=0.231$ using $E_h=128$ elements and $p=2$. Fig. \ref{fig:soddg} illustrates that the DG solution exhibits over- and undershoots near discontinuities. The LO solution, as shown in Fig. \ref{fig:sodlo}, is free of spurious local extrema but has the drawback of being highly diffusive. The numerical solution obtained with the WENO scheme captures discontinuities sharply while preventing the occurrence of spurious oscillations; see Fig. \ref{fig:sodweno}.

% Sod

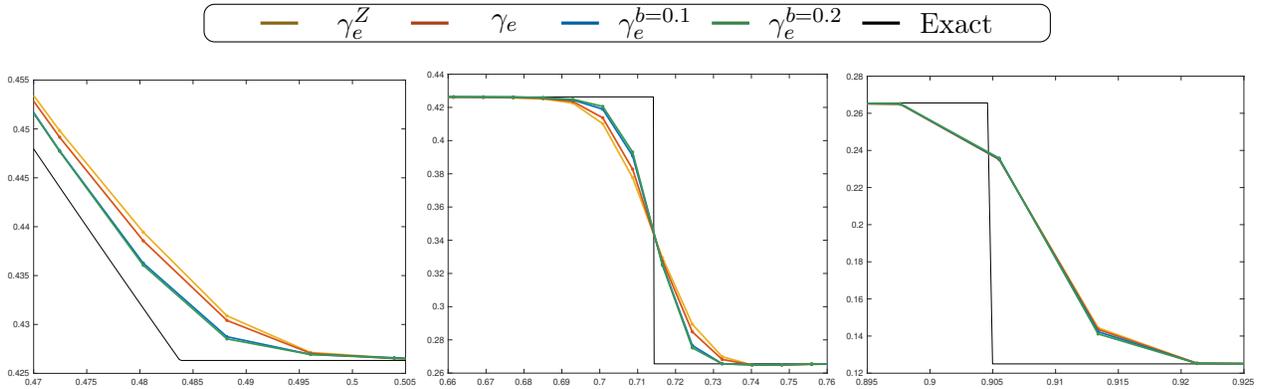
\begin{figure}[!htb]
	\centering
	\begin{subfigure}[b]{\linewidth}
		\centering
		\begin{tikzpicture}
		\draw[rounded corners] (0, 0) rectangle (11.25, 0.5) node[pos=.5]{};
		\draw[very thick, color={rgb:red,0.2431;green,0.5882;blue,0.3176}] (6.75,0.25)--(7.25,0.25);
		\draw[very thick, color=black] (8.75,0.25)--(9.25,0.25);
		\draw[very thick,color={rgb:red,0.9290;green,0.6940;blue,0.1250}] (0.75,0.25)--(1.25,0.25);
		\draw[very thick,color={rgb:red,0.8500;green,0.3250;blue,0.0980}] (2.75,0.25)--(3.25,0.25);
		\draw[very thick,color={rgb:red,0;green,0.4470;blue,0.7410}] (4.75,0.25)--(5.25,0.25);
		\node at (2.00,0.25) (a) {$\gamma_e^{Z}$};
		\node at (4.00,0.25) (a) {$\gamma_e$};
		\node at (6.00,0.25) (a) {$\gamma_e^{b=0.1}$};
		\node at (8,0.25) (a) {$\gamma_e^{b=0.2}$};
		\node at (10,0.25) (a) {Exact};
		\end{tikzpicture}
		\vspace*{0.25cm}
	\end{subfigure}
	\begin{subfigure}[b]{.32\linewidth}
		\includegraphics[width=\linewidth]{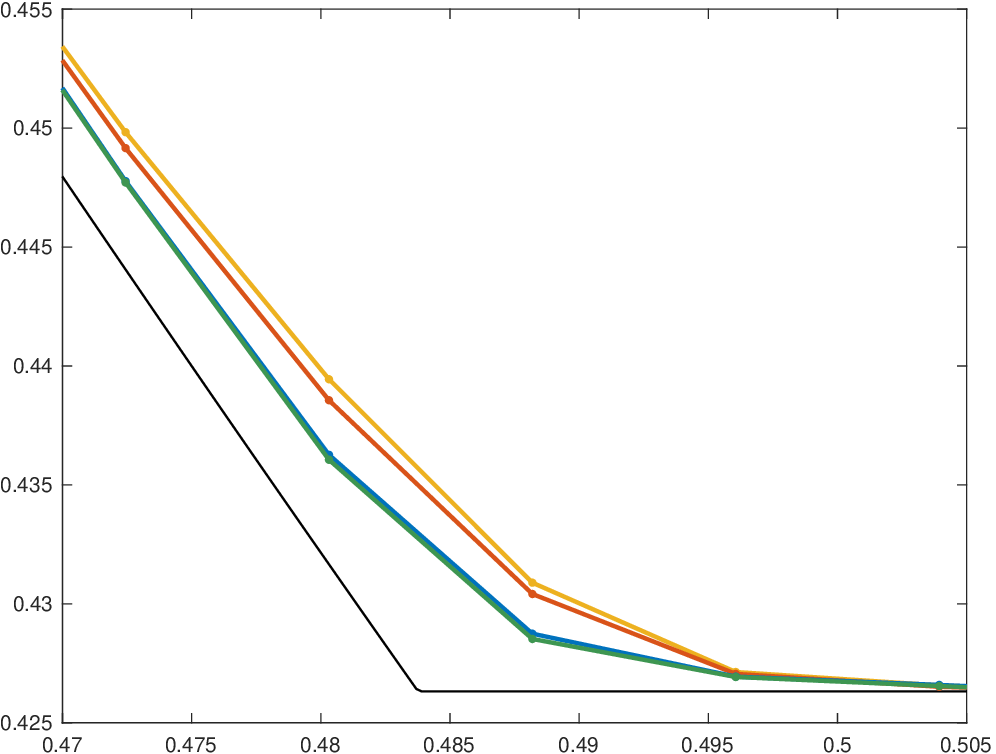}
		\label{fig:sodzoom1}
	\end{subfigure}
	\begin{subfigure}[b]{.32\linewidth}
		\includegraphics[width=\linewidth]{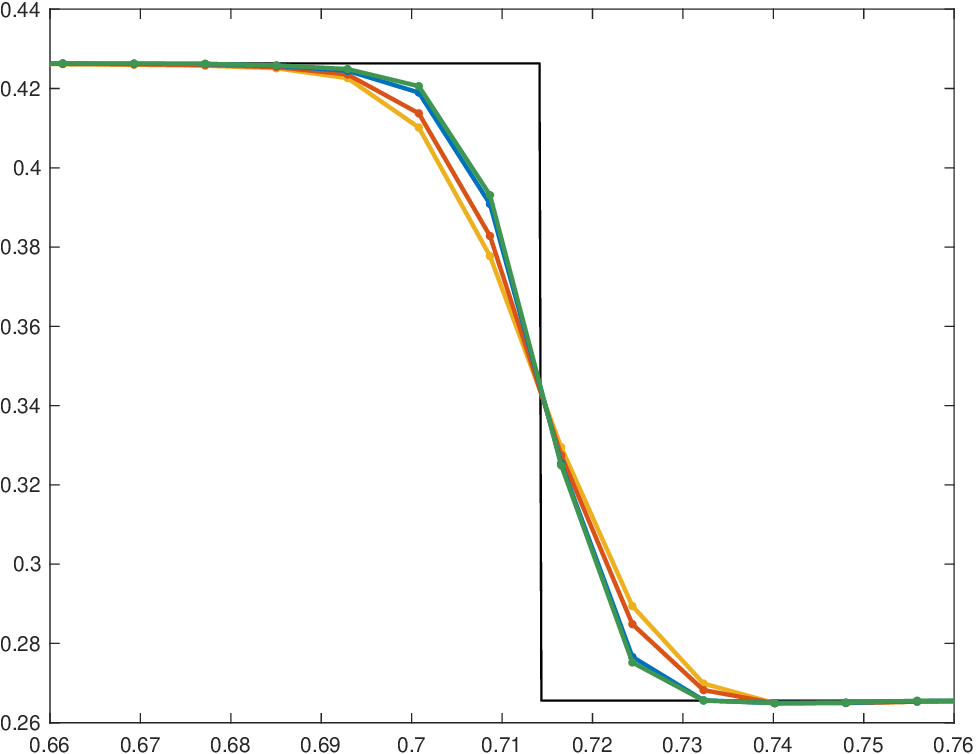}
		\label{fig:sodzoom2}
	\end{subfigure}
	\begin{subfigure}[b]{.32\linewidth}
		\includegraphics[width=\linewidth]{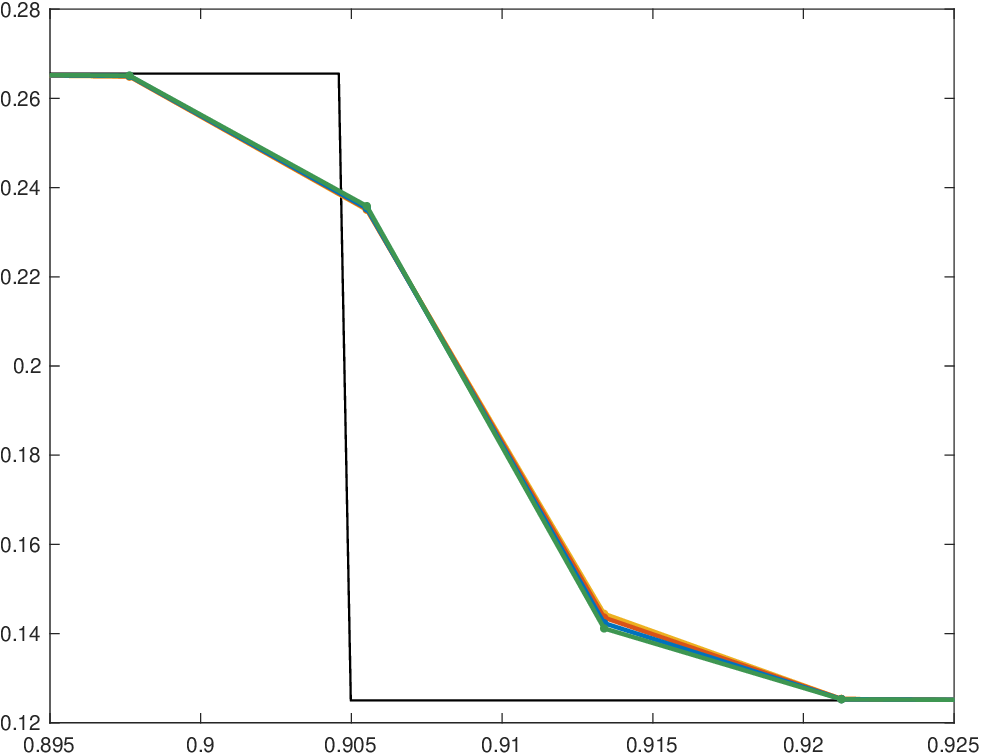}
		\label{fig:sodzoom3}
	\end{subfigure}
	\caption{Sod shock tube, density profiles $\varrho$ at $t=0.231$ obtained using $E_h=128$, $p=2$ and several shock sensors.}
	\label{fig:sodzoom}
\end{figure}

To further evaluate the accuracy of our scheme, we perform a comparison between the shock sensor $\gamma_e$ in \eqref{eq:sensor}, the modified version $\gamma_e^b$ displayed in \eqref{eq:gensensor}, and the shock sensor proposed by Zhao et al. \cite{zhao2020}, i.e.,
\begin{equation}
\gamma_e^Z=\frac{\sum_{l=0}^{m_e}|\frac{\omega_l}{\tilde{\omega}_l}-1|^{\theta}}{|\frac{1}{\min_l\tilde{\omega}_l}-1|^{\theta}+m_e}, \quad 0 \leq l \leq m_e.
\end{equation}
Here, deviations between nonlinear weights $\omega_l$, $l=0,\ldots,m_e$ and ideal linear weights $\tilde{\omega}_l$, $l=0,\ldots,m_e$ assigned to candidate polynomials are computed to assess the regularity of the solution. The free parameter $\theta$ controls the sensitivity to deviations from the ideal weights. We set $\theta = 1$. Continuing with the aforementioned setup, we present, in Fig. \ref{fig:sodzoom}, numerical solutions in different regions of the computational domain. Notably, our WENO-based shock sensor is less diffusive compared to the one introduced in \cite{zhao2020}. By employing the generalized smoothness sensor, the accuracy can be further improved. However, this approach may introduce instabilities for strong shocks.

% Sod-Zoom

\subsubsection{Modified Sod shock tube}
We investigate the entropy stability properties of our scheme for hyperbolic systems using the modified Sod shock tube problem \cite{toro2013}. Here, the initial condition is specified as
\begin{equation}
\begin{bmatrix}
\varrho_L \\v_L\\p_L
\end{bmatrix}=
\begin{bmatrix}1.0\\0.75\\1.0
\end{bmatrix}, \quad 
\begin{bmatrix}
\varrho_R \\v_R\\p_R
\end{bmatrix}=
\begin{bmatrix}0.125\\0.0\\0.1
\end{bmatrix}.
\end{equation}
For this particular problem, many high-order schemes encounter an entropy glitch due to the presence of a sonic point within the rarefaction wave. The computational domain is given by $\Omega=(0,1)$. Unlike the classical Sod shock tube, the left boundary functions as an inlet.

\begin{figure}[!htb]
	\centering
	\begin{subfigure}[b]{.48\linewidth}
		\includegraphics[width=\linewidth]{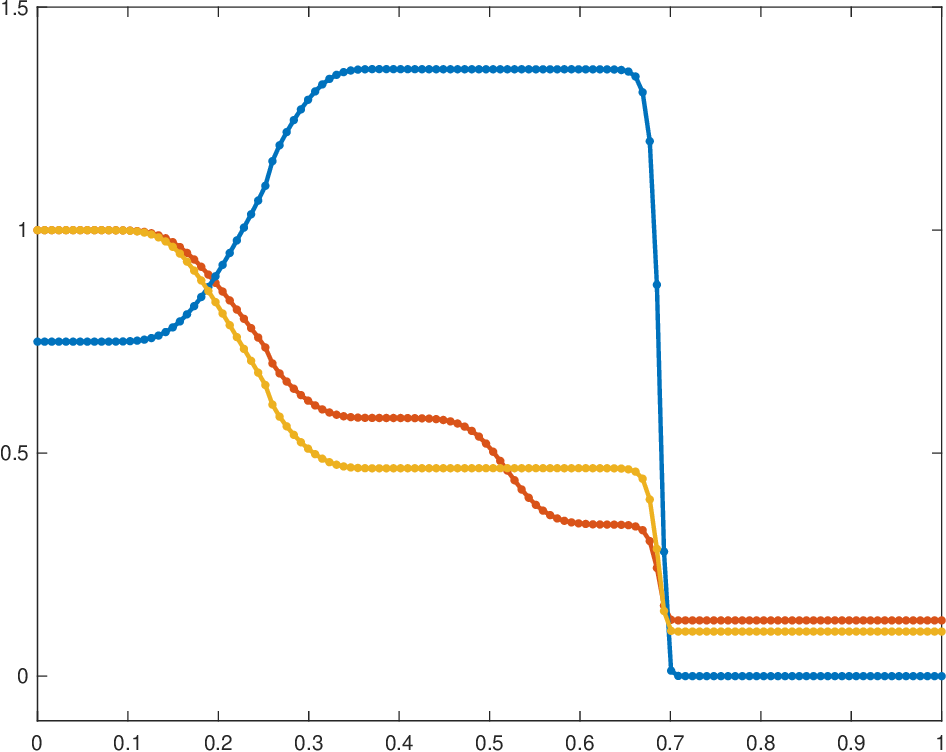}
		\caption{LO}
		\label{fig:msodlo}
	\end{subfigure}
	\begin{subfigure}[b]{.48\linewidth}
		\includegraphics[width=\linewidth]{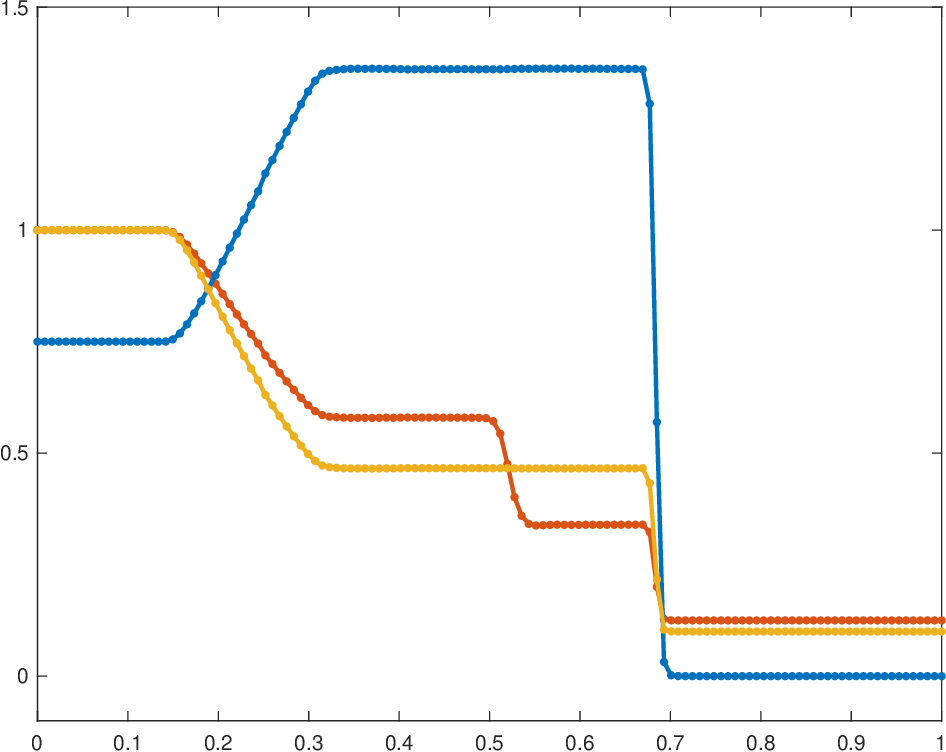}
		\caption{WENO}
		\label{fig:msodweno}
	\end{subfigure}
	\caption{Modified Sod shock tube, $(\varrho,v,E)$ at $t=0.2$ obtained using $E_h=128$ and $p=2$.}
	\label{fig:msod}
\end{figure}

Numerical simulations are carried out using $E_h=128$ elements and $p=2$. The results obtained with the LO method and the WENO method are depicted in Figs \ref{fig:msodlo} and \ref{fig:msodweno}, respectively. Interestingly, when employing the LO method, an entropy shock can be observed in all variables. In contrast, the numerical solutions obtained with the WENO scheme remain entropy stable, thereby eliminating the need for (semi-)discrete entropy fixes in our method.

% MSod

\subsubsection{Lax shock tube}
Next, we consider a 1D tube that is equipped with a diaphragm placed at its center, dividing the tube into two distinct regions with varying pressures. The rupture of the diaphragm at $t=0$ leads to the propagation of a rarefaction wave on the left side and the formation of a contact discontinuity and a shock wave on the right side. This test problem is known as the Lax shock tube problem \cite{lax1954}. The initial data
\begin{equation}
\begin{bmatrix}
\varrho_L \\v_L\\p_L
\end{bmatrix}=
\begin{bmatrix}0.445\\0.698\\3.528
\end{bmatrix}, \quad 
\begin{bmatrix}
\varrho_R \\v_R\\p_R
\end{bmatrix}=
\begin{bmatrix}0.5\\0.0\\0.571
\end{bmatrix}
\end{equation}
are prescribed in the computational domain $\Omega=(0,2)$.

\begin{figure}[!htb]
	\centering
	\begin{subfigure}[b]{\linewidth}
		\centering
		\begin{tikzpicture}
		\draw[rounded corners] (0, 0) rectangle (9.25, 0.5) node[pos=.5]{};
		\draw[very thick, color={rgb:red,0;green,0.4470;blue,0.7410}] (6.75,0.25)--(7.25,0.25);
		\node at (2,0.25) (a) {$p=1$};
		\node at (5,0.25) (a) {$p=2$};
		\node at (8,0.25) (a) {$p=3$};
		\draw[very thick,color={rgb:red,0.9290;green,0.6940;blue,0.1250}] (0.75,0.25)--(1.25,0.25);
		\draw[very thick,color={rgb:red,0.8500;green,0.3250;blue,0.0980}] (3.75,0.25)--(4.25,0.25);
		\end{tikzpicture}
		\vspace*{0.25cm}
	\end{subfigure}
	\begin{subfigure}[b]{.32\linewidth}
		\includegraphics[width=\linewidth]{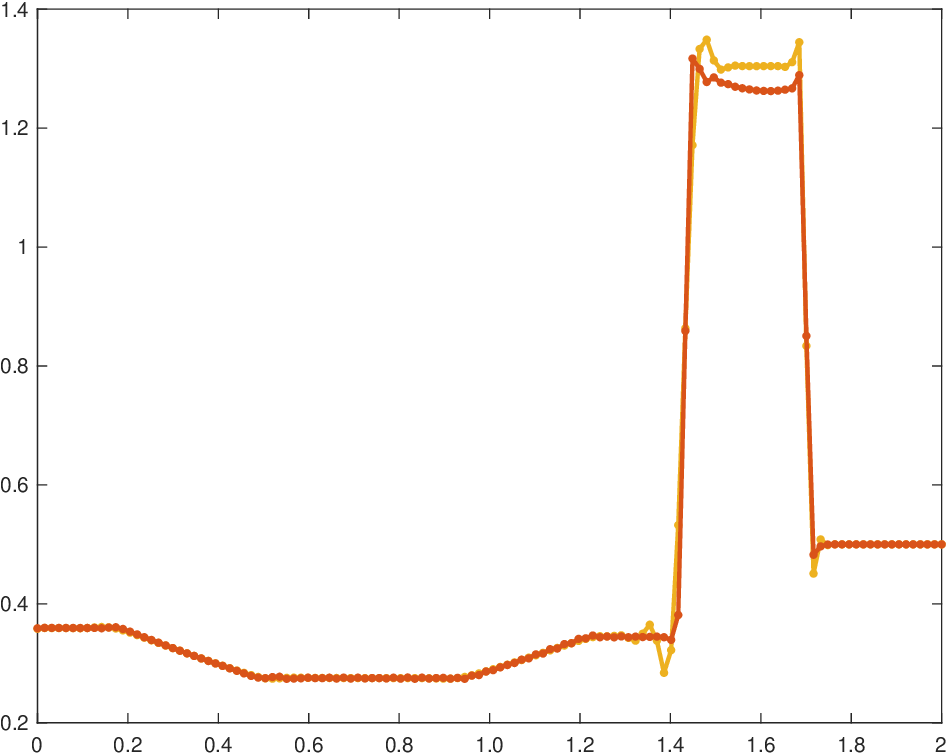}
		\caption{DG}
		\label{fig:laxdg}
	\end{subfigure}
	\begin{subfigure}[b]{.32\linewidth}
		\includegraphics[width=\linewidth]{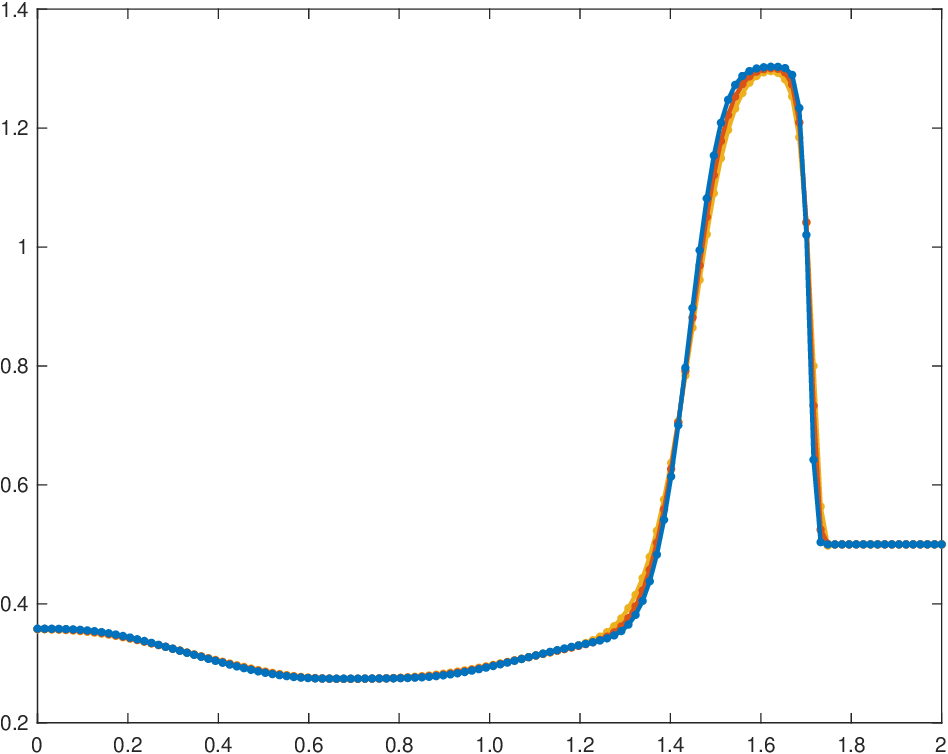}
		\caption{LO}
		\label{fig:laxlo}
	\end{subfigure}
	\begin{subfigure}[b]{.32\linewidth}
		\includegraphics[width=\linewidth]{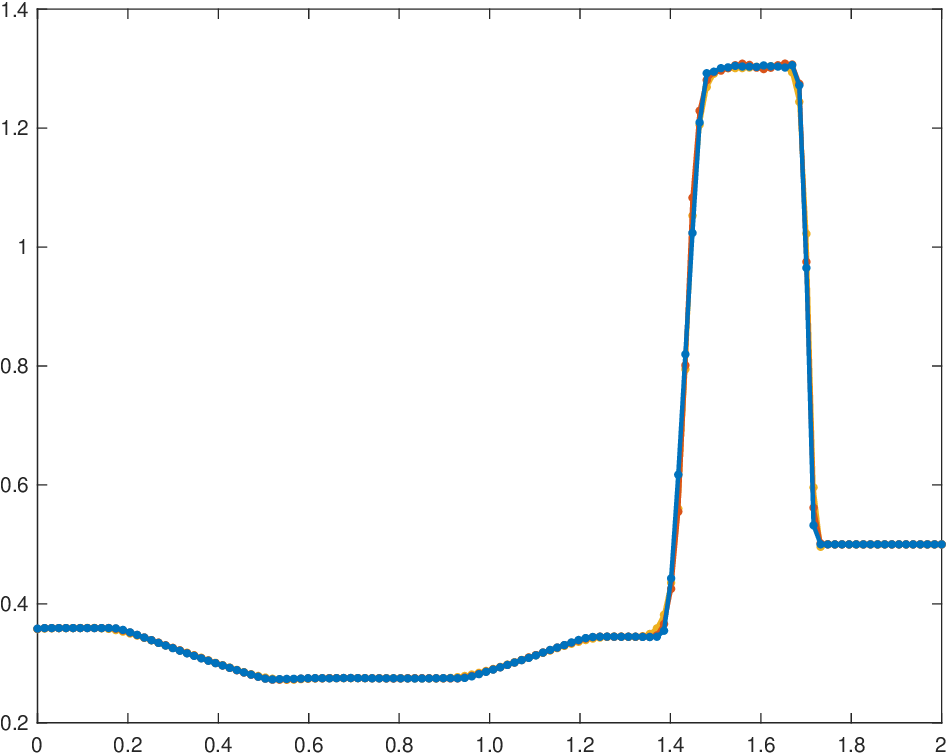}
		\caption{WENO}
		\label{fig:laxweno}
	\end{subfigure}
	\caption{Lax shock tube, density profiles $\varrho$ at $t=0.14$ obtained using $E_h=512$ and $p\in\{1,2,3\}$.}
	\label{fig:lax}
\end{figure}
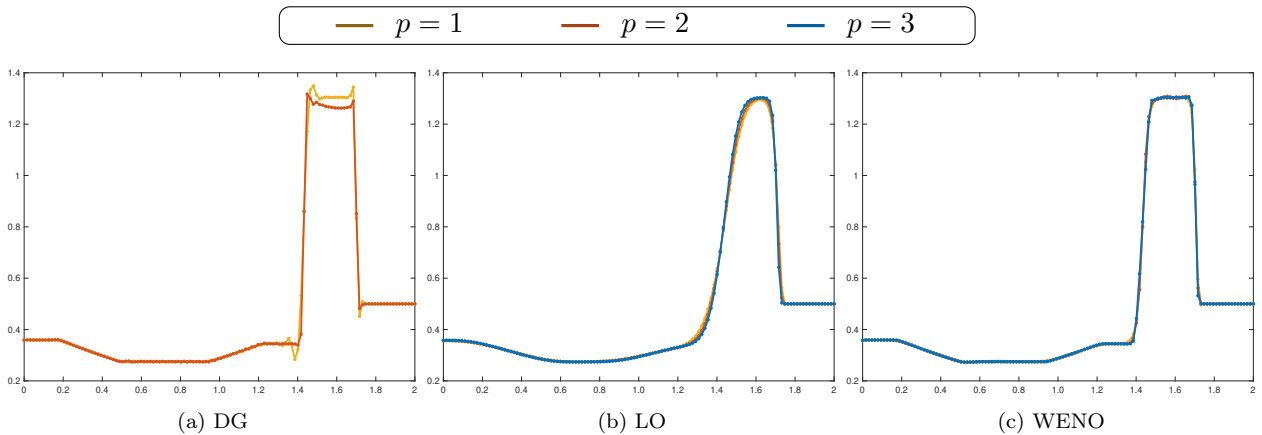

We perform numerical experiments up to the final time $t=0.14$ using $E_h=512$ elements and $p\in \{1,2,3\}$. The density profiles of the DG solutions, LO solutions and WENO solutions are displayed in Figs \ref{fig:laxdg}-\ref{fig:laxweno}, respectively.  Notably, the DG scheme equipped with cubic basis polynomials oscillates heavily, resulting in a distorted solution that is not shown here. Once again, the WENO method demonstrates sharp capturing of discontinuities, while the LO scheme suffers from excessive numerical dissipation.

% Lax

\subsubsection{Shu-Osher problem}
We investigate the sine-shock interaction problem, also known as the Shu-Osher problem \cite{shu1989}. The computational domain $\Omega=(-5, 5)$ is bounded by an inlet on the left boundary and a reflecting wall on the right boundary. This problem serves as a useful test to evaluate the resolution of higher frequency waves behind a shock. The problem is equipped with the initial conditions
\begin{equation}
\begin{bmatrix}
\varrho_L \\v_L\\p_L
\end{bmatrix}=
\begin{bmatrix}3.857143\\2.629369\\10.3333
\end{bmatrix}, \quad 
\begin{bmatrix}
\varrho_R \\v_R\\p_R
\end{bmatrix}=
\begin{bmatrix}1.0+0.2\sin(5x)\\0.0\\1.0
\end{bmatrix}.
\end{equation}

\begin{figure}[!htb]
	\centering
	\begin{subfigure}[b]{\linewidth}
		\centering
		\begin{tikzpicture}
		\draw[rounded corners] (0, 0) rectangle (9.25, 0.5) node[pos=.5]{};
		\draw[very thick, color={rgb:red,0;green,0.4470;blue,0.7410}] (6.75,0.25)--(7.25,0.25);
		\node at (2,0.25) (a) {$p=1$};
		\node at (5,0.25) (a) {$p=2$};
		\node at (8,0.25) (a) {$p=3$};
		\draw[very thick,color={rgb:red,0.9290;green,0.6940;blue,0.1250}] (0.75,0.25)--(1.25,0.25);
		\draw[very thick,color={rgb:red,0.8500;green,0.3250;blue,0.0980}] (3.75,0.25)--(4.25,0.25);
		\end{tikzpicture}
		\vspace*{0.25cm}
	\end{subfigure}
	\begin{subfigure}[b]{.48\linewidth}
		\includegraphics[width=\linewidth]{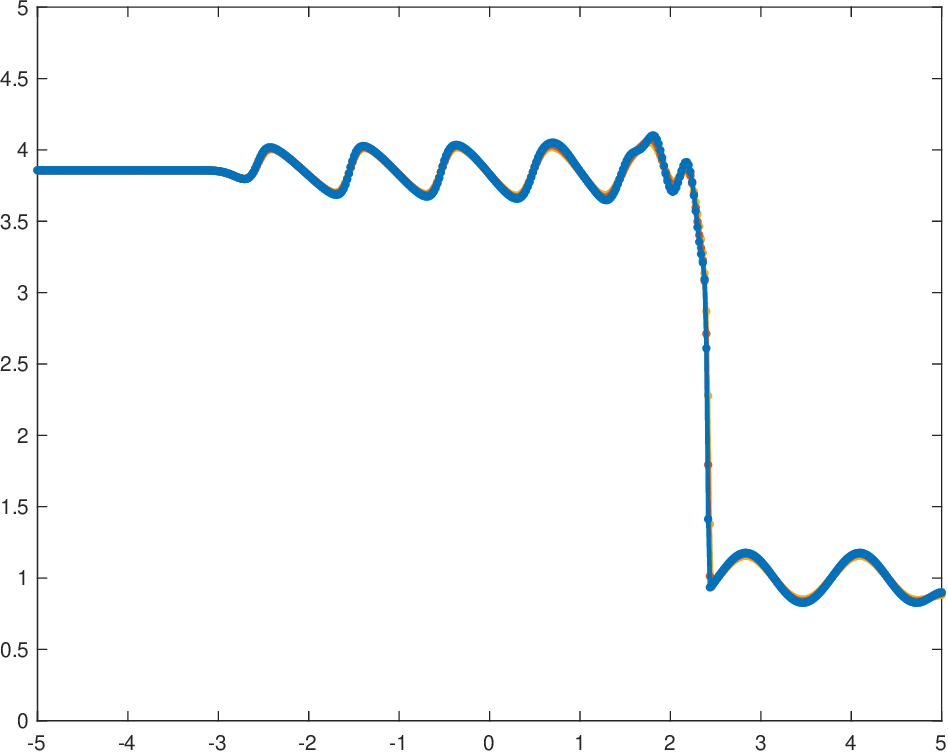}
		\caption{LO}
		\label{fig:solo}
	\end{subfigure}
	\begin{subfigure}[b]{.48\linewidth}
		\includegraphics[width=\linewidth]{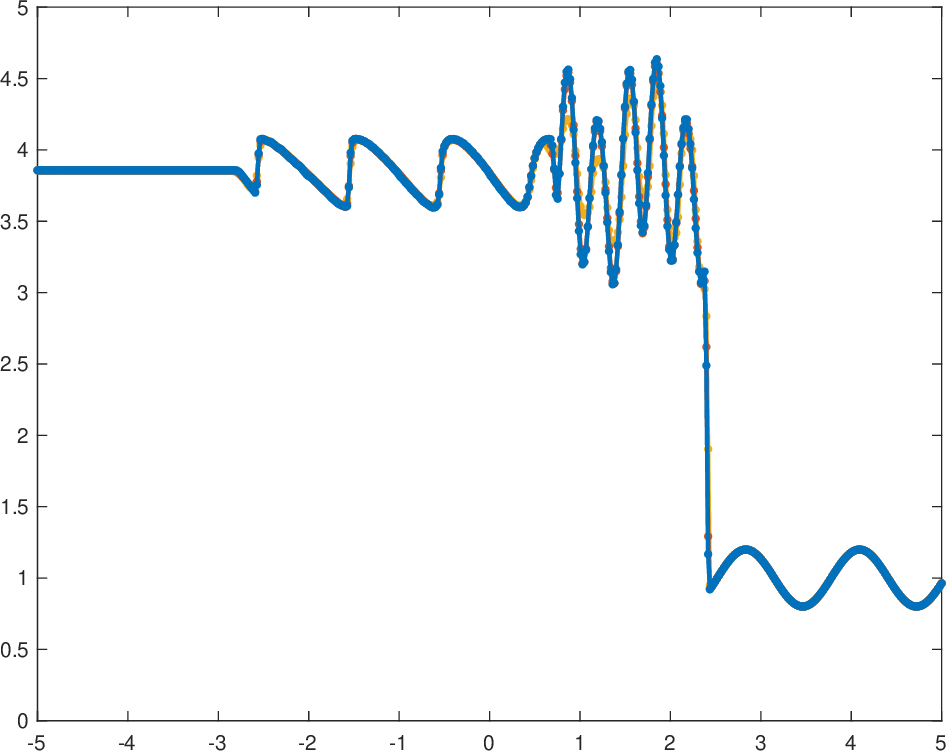}
		\caption{WENO}
		\label{fig:soweno}
	\end{subfigure}
	\caption{Shu-Osher problem, density profiles $\varrho$ at $t=1.8$ obtained using $E_h=512$ and $p\in\{1,2,3\}$.}
	\label{fig:so}
\end{figure}
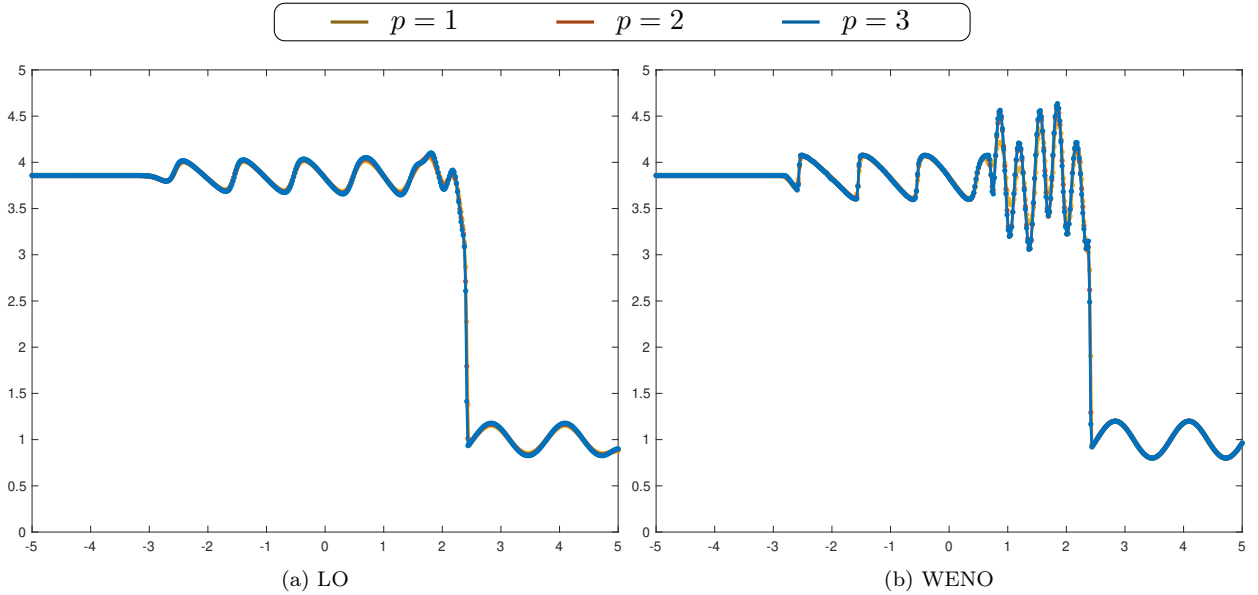

We perform numerical simulations up to the final time $t=1.8$ using $E_h=512$ elements and $p\in\{1,2,3\}$. Figs \ref{fig:solo} and \ref{fig:soweno} show the density profiles of the solutions obtained with the LO scheme and the WENO scheme, respectively. The LO scheme struggles to accurately capture the physical oscillations within the exact solution. Even the WENO scheme equipped with linear basis functions exhibits a significant amount of numerical dissipation. However, when the WENO scheme is equipped with high-order polynomials, it accurately captures all features of the exact solution.

% Shu-Osher

\subsubsection{Woodward-Colella blast wave problem}
We consider the Woodward-Colella blast wave problem \cite{woodward1984}, a challenging test for many high-order numerical schemes. Its solution involves multiple interactions of strong shock waves and rarefaction waves with each other and with contact waves. For a comprehensive understanding of all wave interactions, we refer the reader to \cite{woodward1984}. This problem serves as a rigorous test to assess the robustness of numerical schemes. The initial data
\begin{equation}
\begin{bmatrix}
\varrho_L \\v_L\\p_L
\end{bmatrix}=
\begin{bmatrix}1.0\\0.0\\1000.0
\end{bmatrix}, \quad 
\begin{bmatrix}
\varrho_M \\v_M\\p_M
\end{bmatrix}=
\begin{bmatrix}1.0\\0.0\\0.1
\end{bmatrix},\quad 
\begin{bmatrix}
\varrho_R \\v_R\\p_R
\end{bmatrix}=
\begin{bmatrix}1.0\\0.0\\100.0
\end{bmatrix}
\end{equation}
are prescribed in the computational domain $\Omega=(0,1)$, which is bounded by reflecting walls.

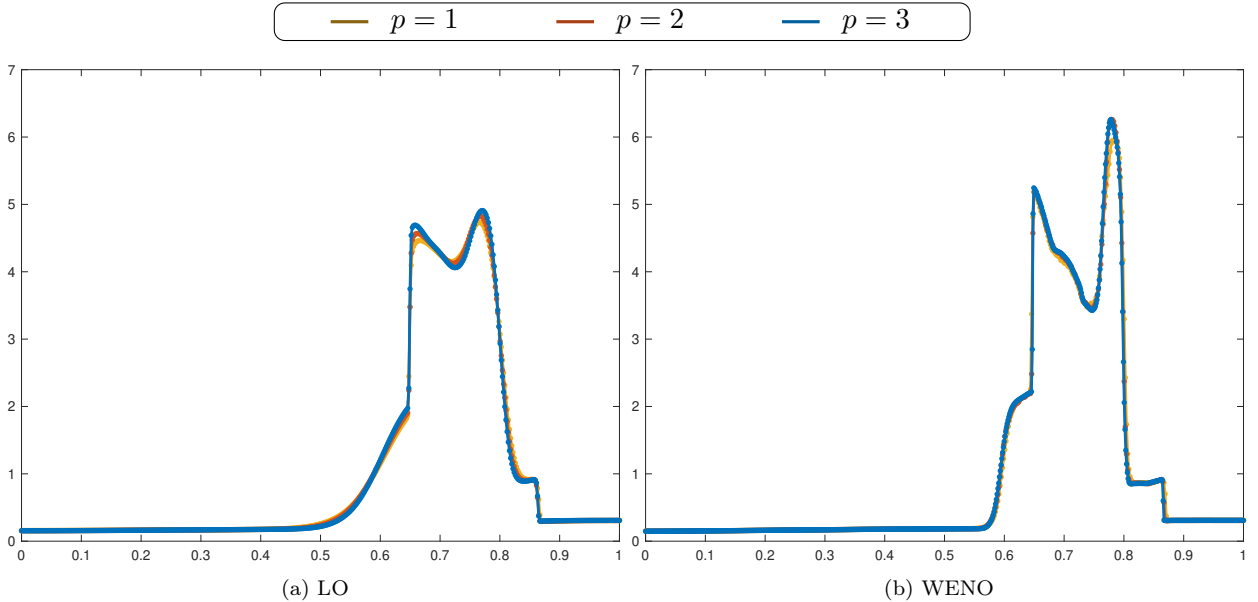
\begin{figure}[!htb]
	\centering
	\begin{subfigure}[b]{\linewidth}
		\centering
		\begin{tikzpicture}
		\draw[rounded corners] (0, 0) rectangle (9.25, 0.5) node[pos=.5]{};
		\draw[very thick, color={rgb:red,0;green,0.4470;blue,0.7410}] (6.75,0.25)--(7.25,0.25);
		\node at (2,0.25) (a) {$p=1$};
		\node at (5,0.25) (a) {$p=2$};
		\node at (8,0.25) (a) {$p=3$};
		\draw[very thick,color={rgb:red,0.9290;green,0.6940;blue,0.1250}] (0.75,0.25)--(1.25,0.25);
		\draw[very thick,color={rgb:red,0.8500;green,0.3250;blue,0.0980}] (3.75,0.25)--(4.25,0.25);
		\end{tikzpicture}
		\vspace*{0.25cm}
	\end{subfigure}
	\begin{subfigure}[b]{.48\linewidth}
		\includegraphics[width=\linewidth]{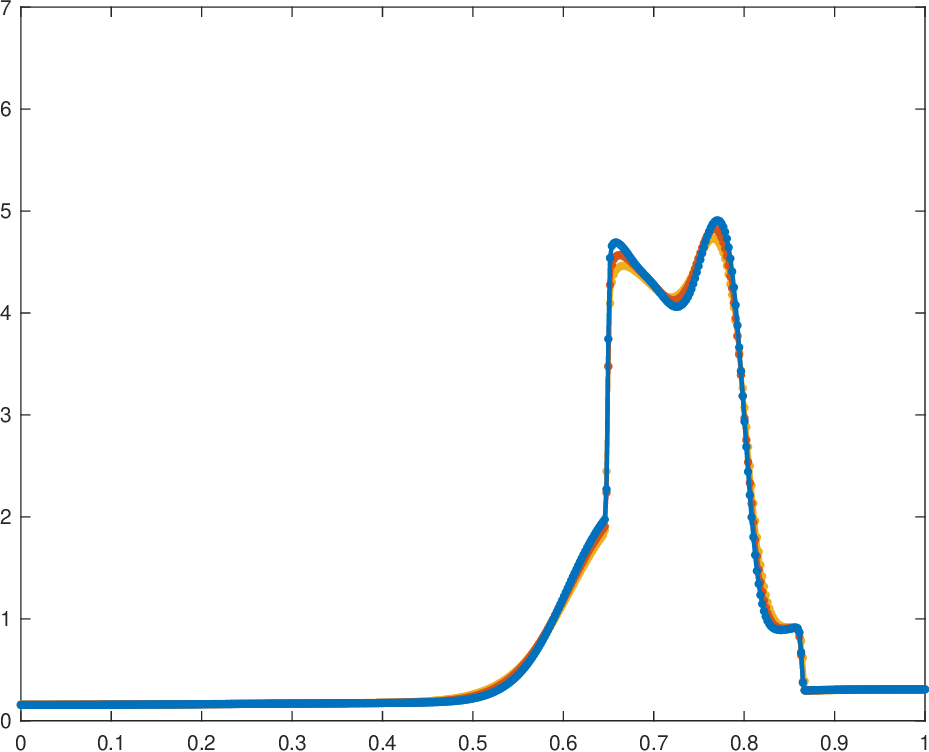}
		\caption{LO}
		\label{fig:wclo}
	\end{subfigure}
	\begin{subfigure}[b]{.48\linewidth}
		\includegraphics[width=\linewidth]{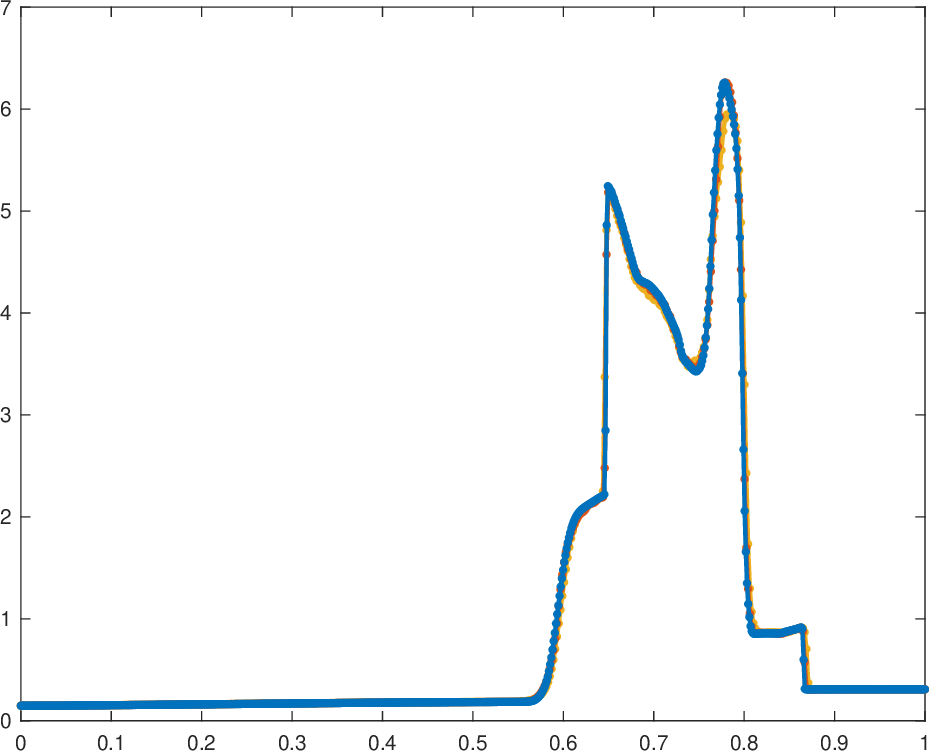}
		\caption{WENO}
		\label{fig:wcweno}
	\end{subfigure}
	\caption{Woodward-Colella blast wave problem, density profiles $\varrho$ at $t=0.038$ obtained using $E_h=512$ and $p\in\{1,2,3\}$.}
	\label{fig:wc}
\end{figure}

Fig. \ref{fig:wc} displays the density profiles of the LO solutions and the WENO solutions at the final time $t=0.038$. The numerical results were obtained using $E_h=512$ elements and $p\in\{1,2,3\}$. It is evident that the WENO scheme achieves higher accuracy compared to the dissipative LO scheme and  accurately captures all sharp features without any visible oscillations. Similarly to the Shu-Osher problem, the benefits of using high-order polynomials are clearly visible in this test.

% Woodward-Colella

\subsubsection{Kelvin-Helmholtz instability}

\begin{figure}[!htb]
	\centering
	\begin{subfigure}[b]{.48\linewidth}
		\includegraphics[width=\linewidth]{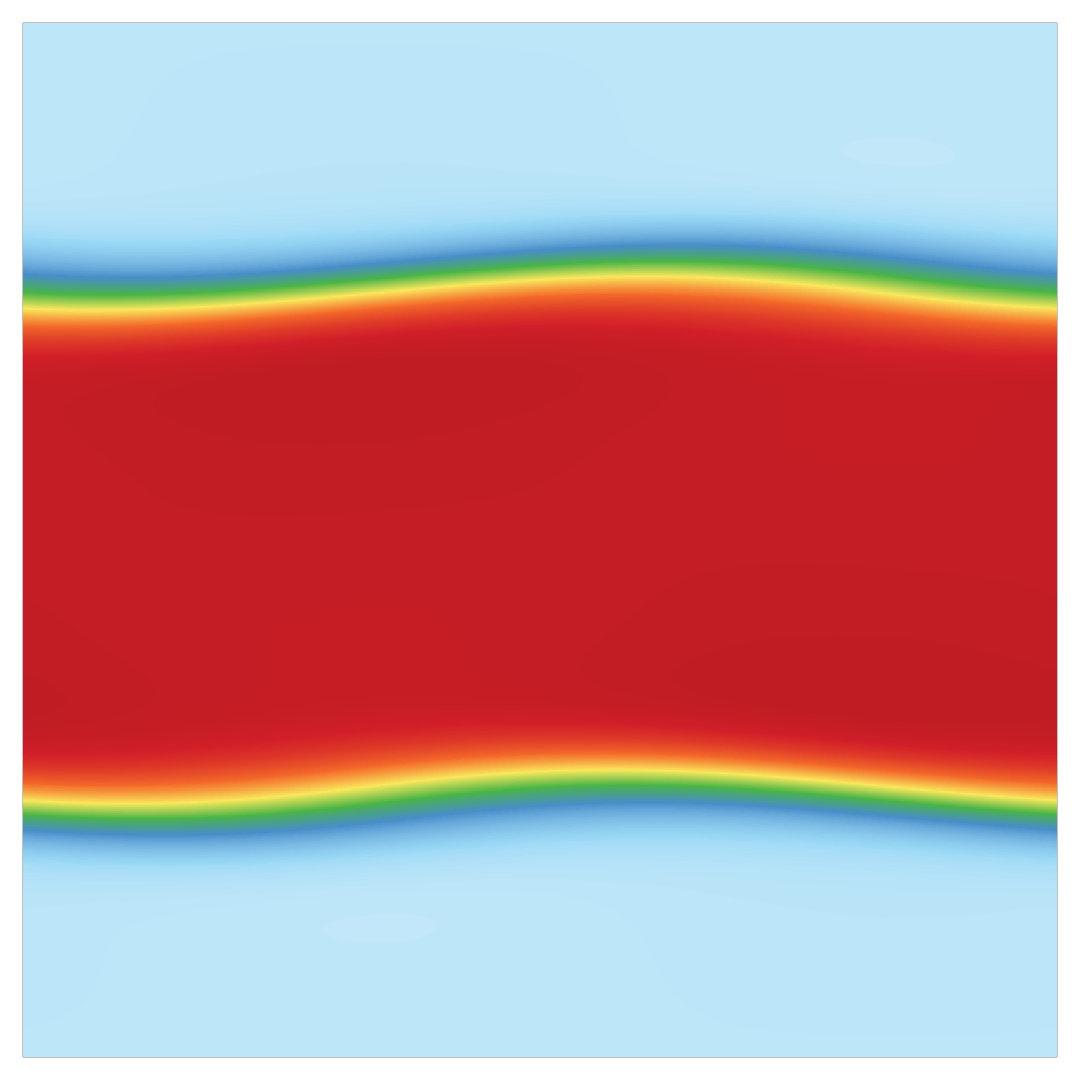}
		\caption{LO, $u_h \in [0.998,2.018]$}
		\label{fig:khlo}
	\end{subfigure}
	\begin{subfigure}[b]{.48\linewidth}
		\includegraphics[width=\linewidth]{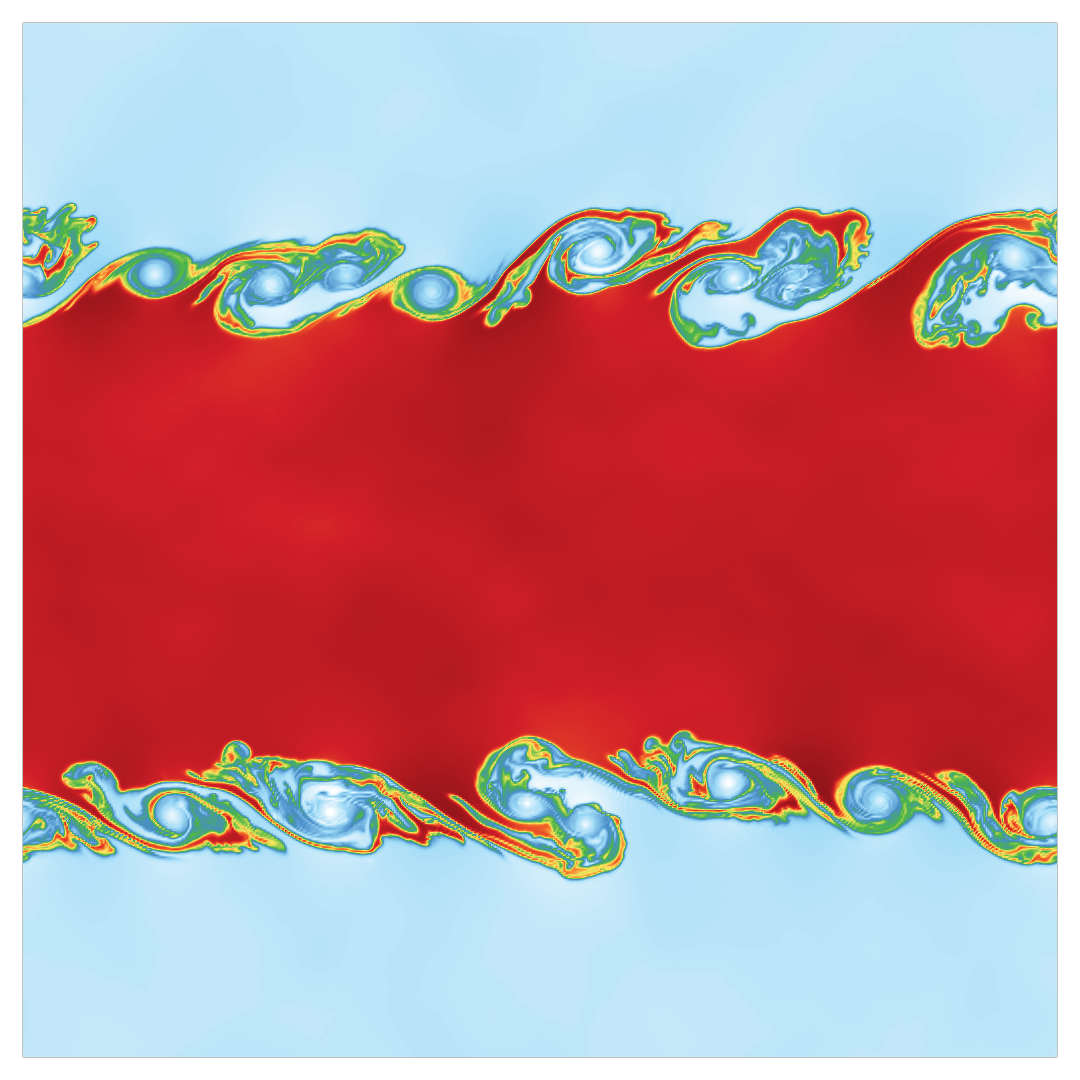}
		\caption{WENO, $u_h \in [0.817,2.179]$}
		\label{fig:khweno4}
	\end{subfigure}
	\caption{Kelvin-Helmholtz instability, density profiles $\varrho$ at $t=1.0$ obtained using $E_h=512^2$ and $p=2$.}
	\label{fig:kh}
\end{figure}

Next, we consider a Kelvin-Helmholtz instability problem \cite{ma2023}, which transforms a narrow shear layer into a complex pattern of vortices. Such problems are frequently employed as standard benchmarks to evaluate the accuracy and dissipation characteristics of numerical methods in converting a linear disturbance into a nonlinear state.
The initial data in the computational domain $\Omega=(0,1)\times(0,1)$ are given by
\begin{equation}
\begin{bmatrix}
\varrho_1 \\v_{x,1}\\v_{y,1}\\p_1
\end{bmatrix}=
\begin{bmatrix}\phantom{-}2.0\\-0.5\\0.01\sin(2\pi(x-0.5))\\\phantom{-}2.5
\end{bmatrix}, \quad 
\begin{bmatrix}
\varrho_2 \\v_{x,2}\\v_{y,2}\\p_2
\end{bmatrix}=
\begin{bmatrix}1.0\\0.5\\0.01\sin(2\pi(x-0.5))\\2.5
\end{bmatrix}.
\end{equation}
We run the simulation up to the final time $t=1.0$ on $E_h=512^2$ elements and $p=2$. Fig. \ref{fig:kh} displays the density profiles of both the LO and the WENO solution. It is evident that the LO scheme introduces a significant level of dissipation, preventing the resolution of vortical structures in the solution. In contrast, the solution obtained with the WENO scheme demonstrates improved accuracy and reveals fine-scale vortical structures.

% Kelvin-Helmholtz instability

\subsubsection{Rarefaction interaction}

We investigate the interaction of two-dimensional planar rarefaction waves. While one-dimensional rarefaction waves yield only continuous solutions, the behavior of two-dimensional rarefaction waves shifts from continuous flows to the presence of transonic shocks.  We consider four planar rarefaction waves defined by the initial conditions

\begin{equation}
\begin{bmatrix}
\varrho_1 \\v_{x,1}\\v_{y,1}\\p_1
\end{bmatrix}=
\begin{bmatrix}1.0\\0.6233\\0.6233\\1.5
\end{bmatrix}, \quad
\begin{bmatrix}
\varrho_2 \\v_{x,2}\\v_{y,2}\\p_2
\end{bmatrix}=
\begin{bmatrix}\phantom{-}0.389\phantom{0}\\-0.6233\\\phantom{-}0.6233\\\phantom{-}0.4\phantom{000}
\end{bmatrix}, \quad
\begin{bmatrix}
\varrho_3 \\v_{x,3}\\v_{y,3}\\p_3
\end{bmatrix}=
\begin{bmatrix}\phantom{-}1.0\phantom{000}\\-0.6233\\-0.6233\\\phantom{-}1.5\phantom{000}
\end{bmatrix}, \quad 
\begin{bmatrix}
\varrho_4 \\v_{x,4}\\v_{y,4}\\p_4
\end{bmatrix}=
\begin{bmatrix}\phantom{-}0.389\phantom{0}\\\phantom{-}0.6233\\-0.6233\\\phantom{-}0.4\phantom{000}
\end{bmatrix},
\label{eq:ic1}
\end{equation}

and 

\begin{equation}
\begin{bmatrix}
\varrho_1 \\v_{x,1}\\v_{y,1}\\p_1
\end{bmatrix}=
\begin{bmatrix}1.0\\0.0312\\0.0312\\0.5
\end{bmatrix}, \quad
\begin{bmatrix}
\varrho_2 \\v_{x,2}\\v_{y,2}\\p_2
\end{bmatrix}=
\begin{bmatrix}\phantom{-}0.927\phantom{0}\\-0.0312\\\phantom{-}0.0312\\\phantom{-}0.45\phantom{00}
\end{bmatrix}, \quad
\begin{bmatrix}
\varrho_3 \\v_{x,3}\\v_{y,3}\\p_3
\end{bmatrix}=
\begin{bmatrix}1.0\\-0.0312\\-0.0312\\0.5
\end{bmatrix}, \quad 
\begin{bmatrix}
\varrho_4 \\v_{x,4}\\v_{y,4}\\p_4
\end{bmatrix}=
\begin{bmatrix}\phantom{-}0.927\phantom{0}\\\phantom{-}0.0312\\-0.0312\\\phantom{-}0.45\phantom{00}
\end{bmatrix}.
\label{eq:ic2}
\end{equation}

\begin{figure}[!htb]
	\centering
	\begin{subfigure}[b]{.32\linewidth}
		\includegraphics[width=\linewidth]{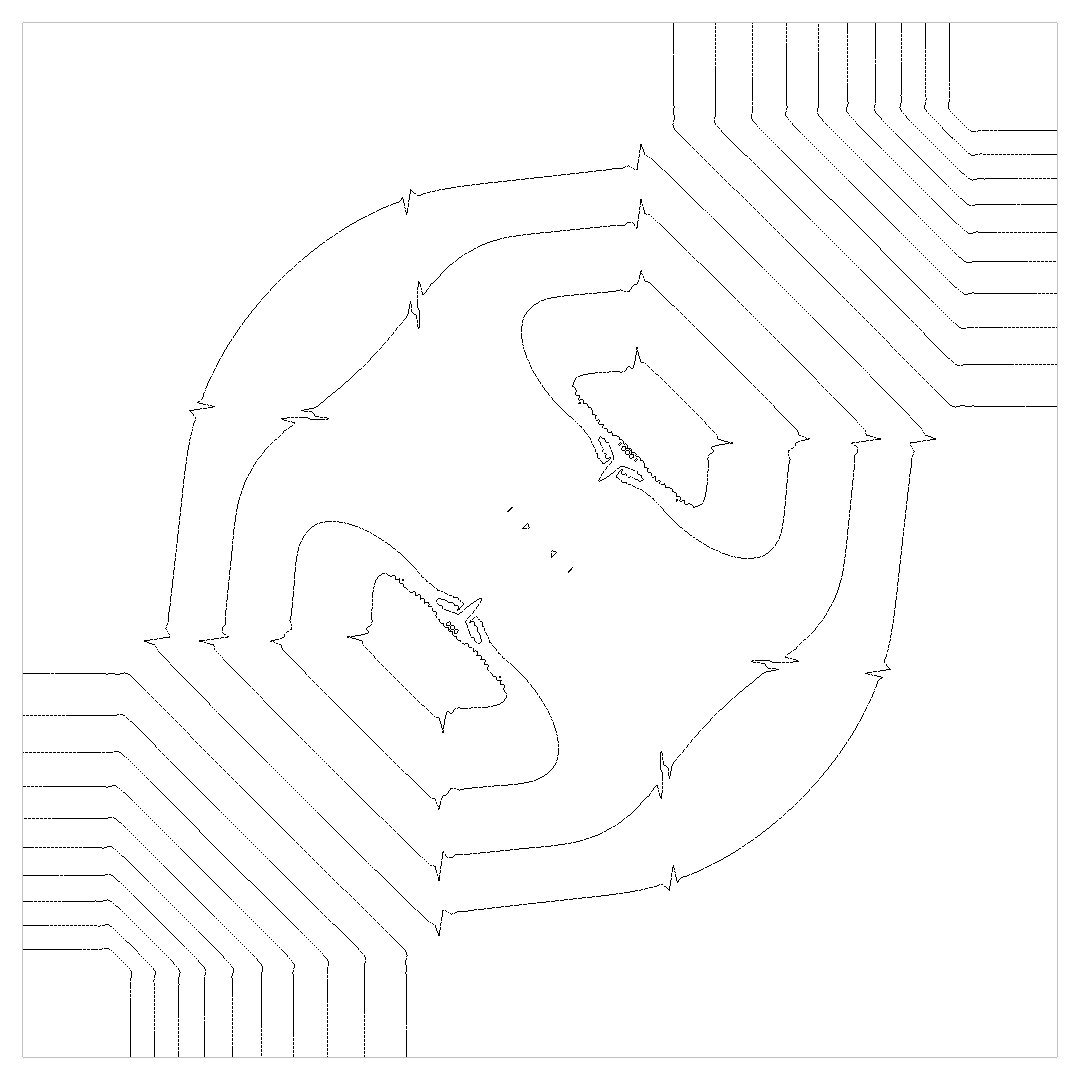}
		\caption{DG, $u_h \in [0.113,1.006]$}
		\label{fig:rriwdg}
	\end{subfigure}
	\begin{subfigure}[b]{.32\linewidth}
		\includegraphics[width=\linewidth]{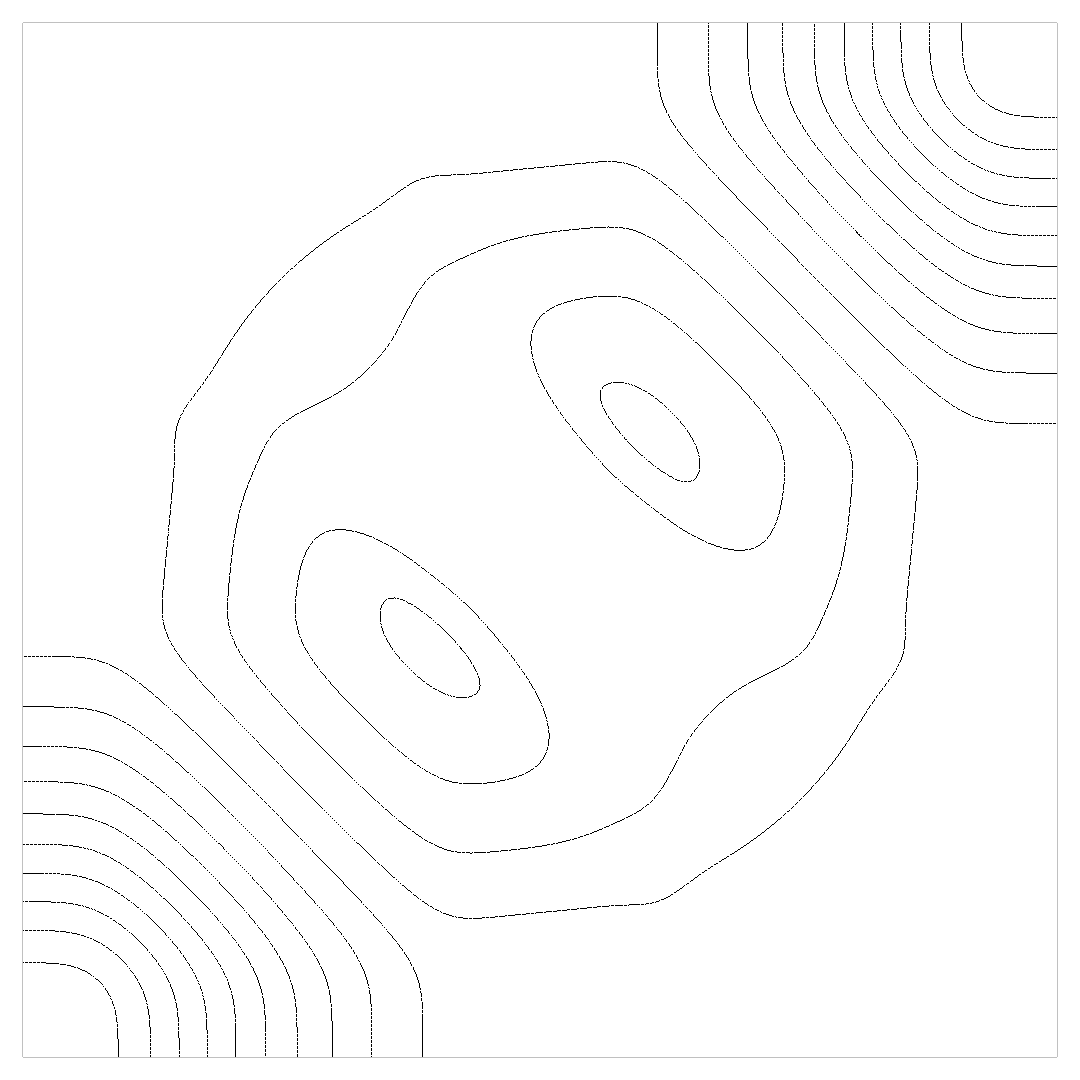}
		\caption{LO, $u_h \in [0.152,1.000]$}
		\label{fig:rriwlo}
	\end{subfigure}
	\begin{subfigure}[b]{.32\linewidth}
		\includegraphics[width=\linewidth]{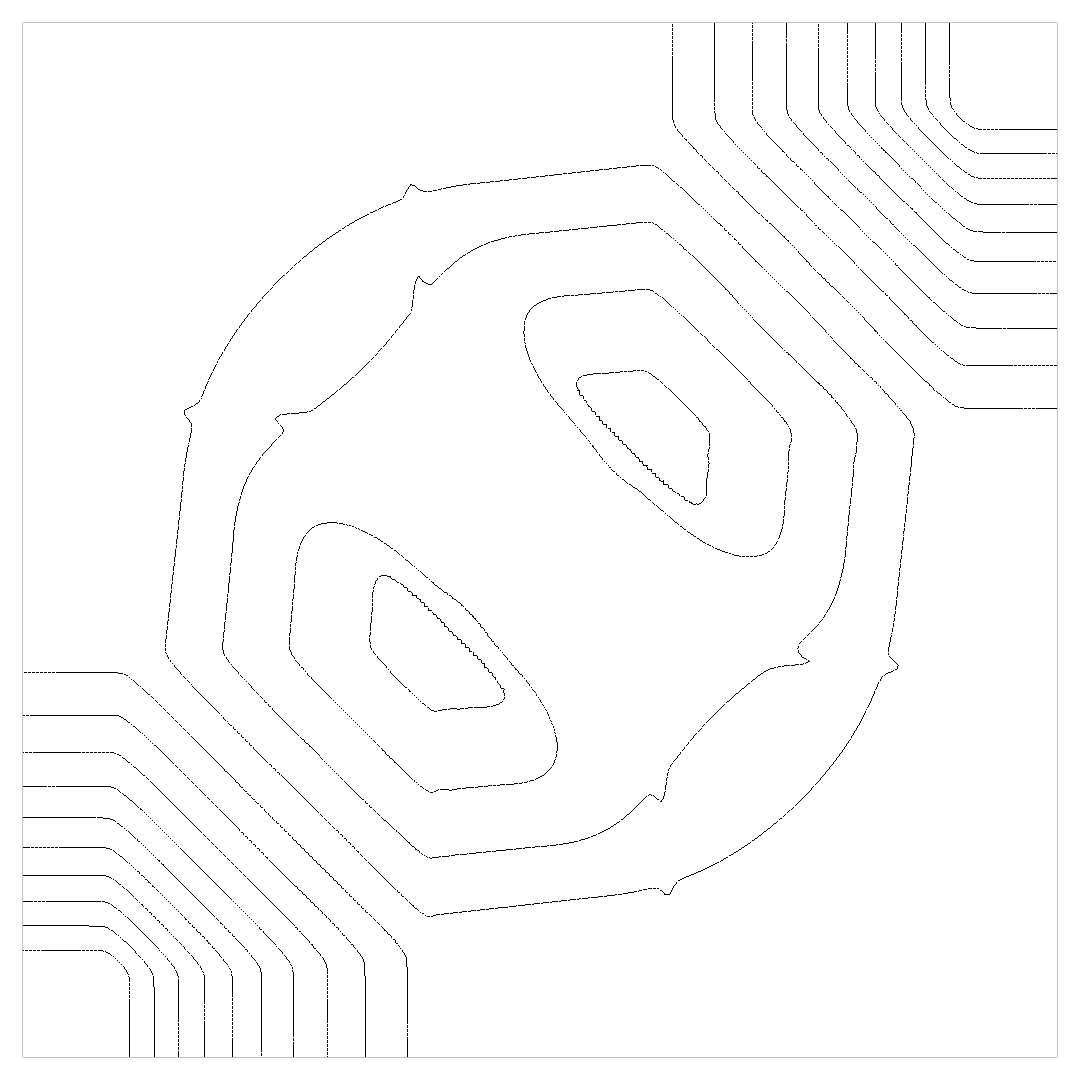}
		\caption{WENO, $u_h \in [0.134,1.000]$}
		\label{fig:rriwweno}
	\end{subfigure}
	\caption{Weak rarefaction interaction, density profiles at $t=0.2$ obtained using $E_h=256^2$ and $p=2$.}
	\label{fig:rriw}
\end{figure}

Numerical results are computed in the domain $\Omega=(0,1)\times(0,1)$ using outflow boundary conditions everywhere. The solutions in Fig. \ref{fig:rriw} are produced using the initial condition \eqref{eq:ic1}, where the rarefaction waves are weak, resulting in a continuous solution. Conversely, the application of initial condition \eqref{eq:ic2} induces strong rarefaction waves, leading to the emergence of shocks within the interior domain.  For more details, we refer the reader to \cite{glimm2008,li2010}. 

% Weak rarefaction interaction
% Strong rarefaction interaction
\begin{figure}[!htb]
	\centering
	\begin{subfigure}[b]{.32\linewidth}
		\includegraphics[width=\linewidth]{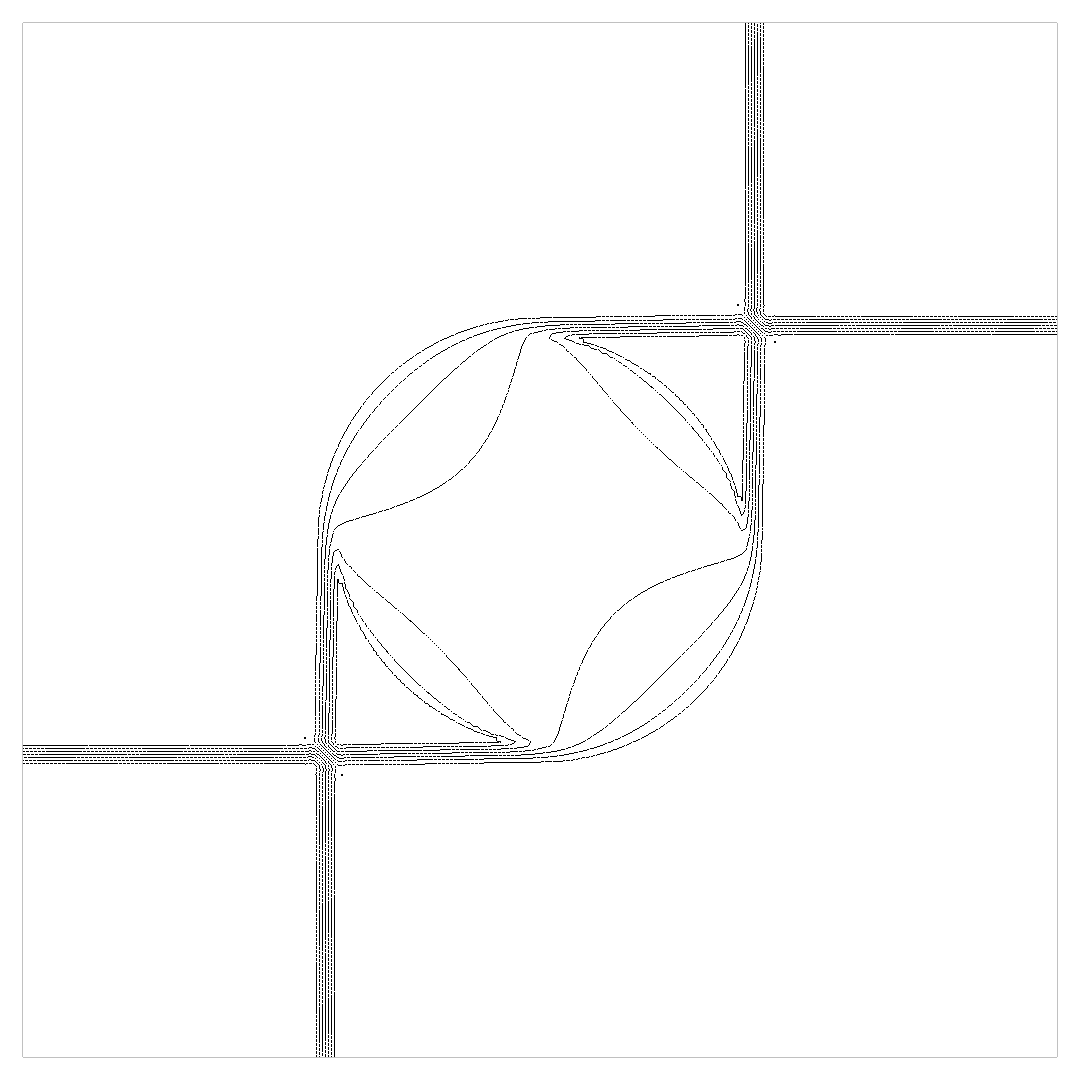}
		\caption{DG, $u_h \in [0.857,1.004]$}
		\label{fig:rrisdg}
	\end{subfigure}
	\begin{subfigure}[b]{.32\linewidth}
		\includegraphics[width=\linewidth]{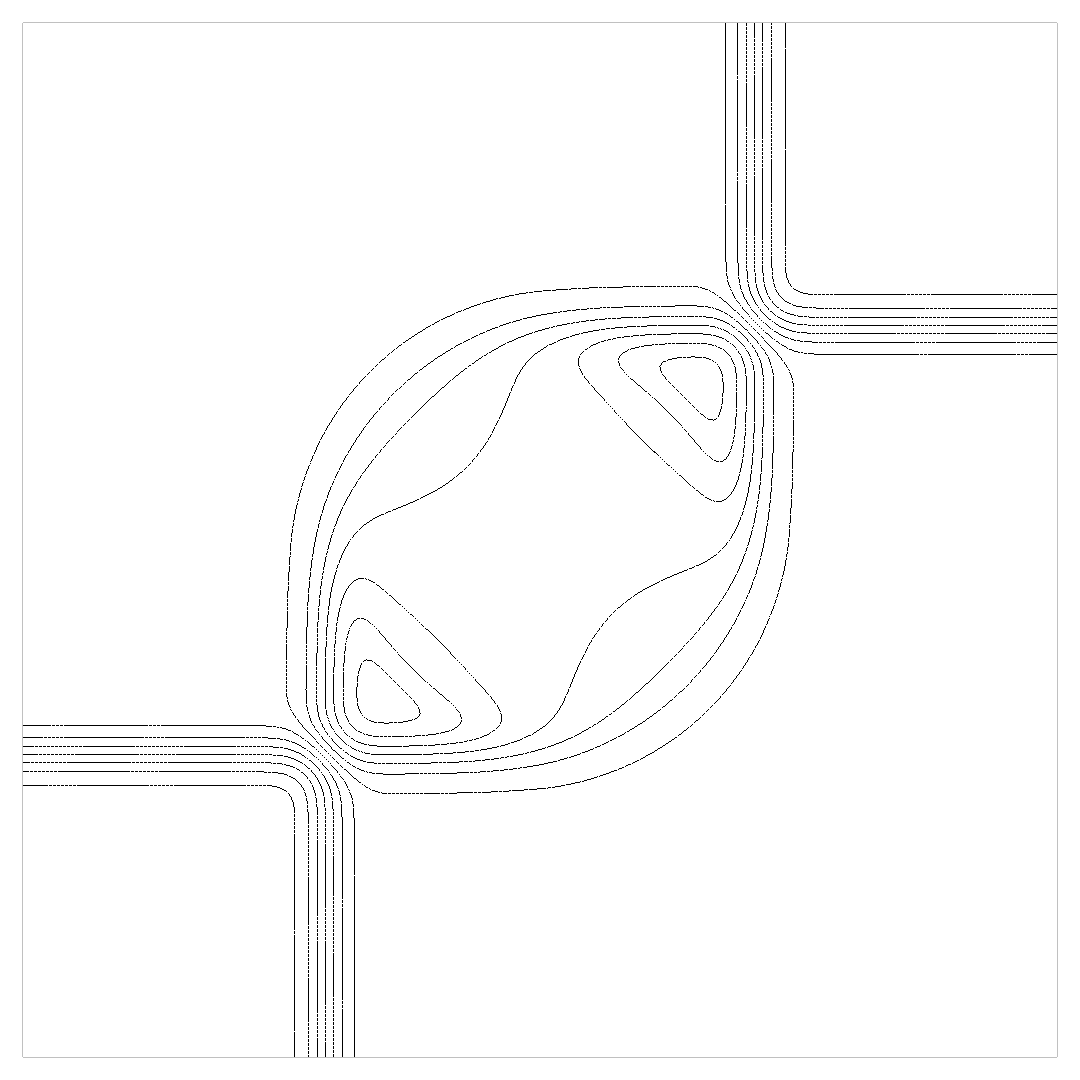}
		\caption{LO, $u_h \in [0.862,1.000]$}
		\label{fig:rrislo}
	\end{subfigure}
	\begin{subfigure}[b]{.32\linewidth}
		\includegraphics[width=\linewidth]{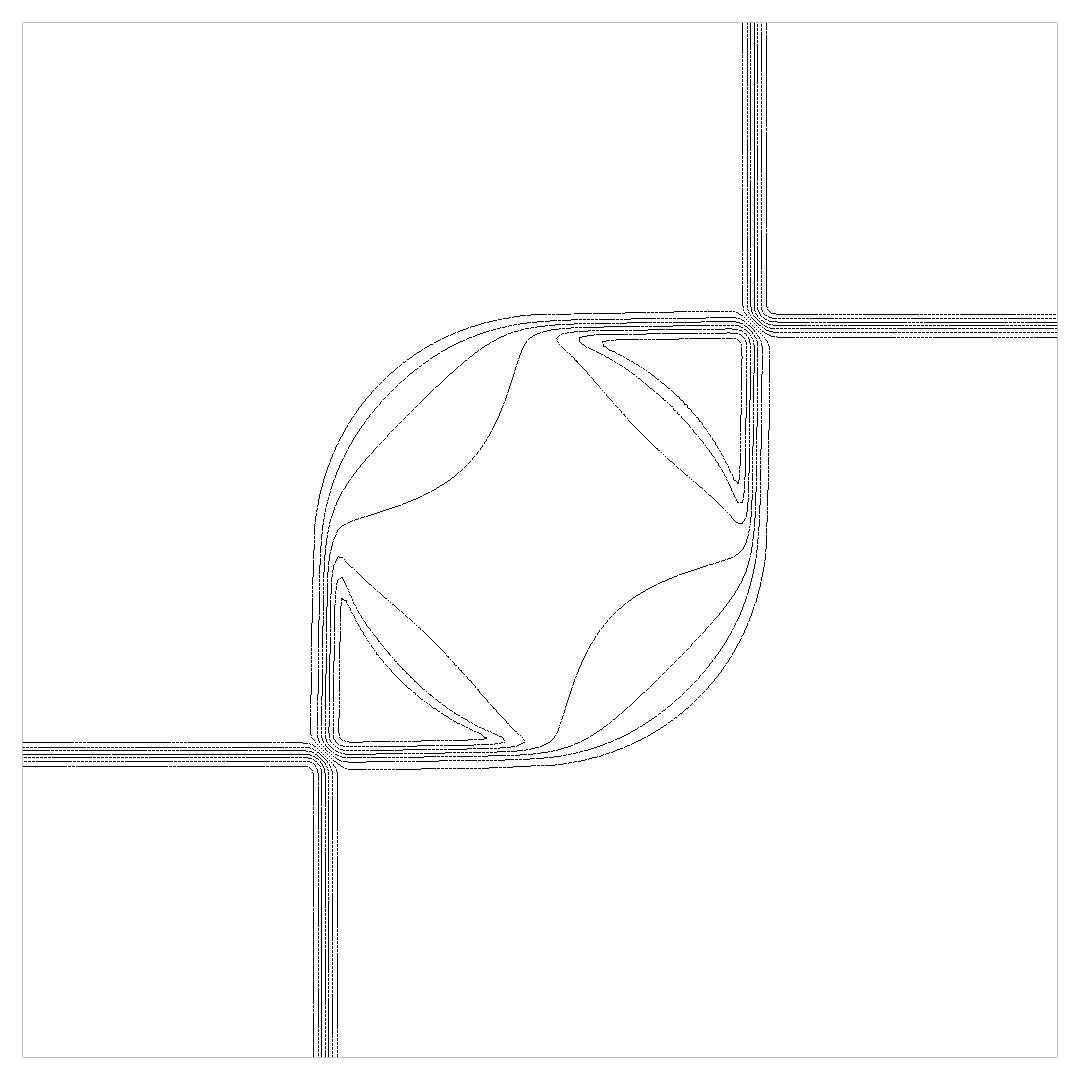}
		\caption{WENO, $u_h \in [0.860,1.000]$}
		\label{fig:rrisweno}
	\end{subfigure}
	\caption{Strong rarefaction interaction, density profiles at $t=0.25$ obtained using $E_h=256^2$ and $p=2$.}
	\label{fig:rris}
\end{figure}

\subsubsection{Double Mach reflection}
In our next numerical example, we consider the double Mach reflection problem of Woodward and Colella \cite{woodward1984}. The computational domain is the rectangle $\Omega=(0,4)\times(0,1)$. In this benchmark, the flow pattern features a Mach 10 shock in air which initially makes a $60^{\circ}$ angle with a reflecting wall. The following pre-shock and post-shock values of the flow variables are used
\begin{equation}
\begin{bmatrix}
\varrho_L \\v_{x,L}\\v_{y,L}\\p_L
\end{bmatrix}=
\begin{bmatrix}8.0\\\phantom{-}8.25\cos({30}^{\circ})\\-8.25\sin({30}^{\circ})\\116.5
\end{bmatrix}, \quad 
\begin{bmatrix}
\varrho_R \\v_{x,R}\\v_{y,R}\\p_R
\end{bmatrix}=
\begin{bmatrix}1.4\\0.0\\0.0\\1.0
\end{bmatrix}.
\end{equation}
Initially, the post-shock values (subscript L) are prescribed in the subdomain $\Omega_L=\{(x,y)\;|\;x<\frac{1}{6}+\frac{y}{\sqrt{3}}\}$ and the pre-shock values (subscript R) in $\Omega_R = \Omega \setminus \Omega_L$. The reflecting wall corresponds to $1/6 \leq x \leq 4$ and $y=0$. No boundary conditions are required along the line $x=4$. On the rest of the boundary, the post-shock conditions are prescribed for $x<\frac{1}{6}+\frac{1+20t}{\sqrt{3}}$ and the pre-shock conditions elsewhere. The so-defined values along the top boundary describe the exact motion of the initial Mach 10 shock.

\begin{figure}[!htb]
	\centering
	\begin{subfigure}[b]{.48\linewidth}
		\includegraphics[width=\linewidth]{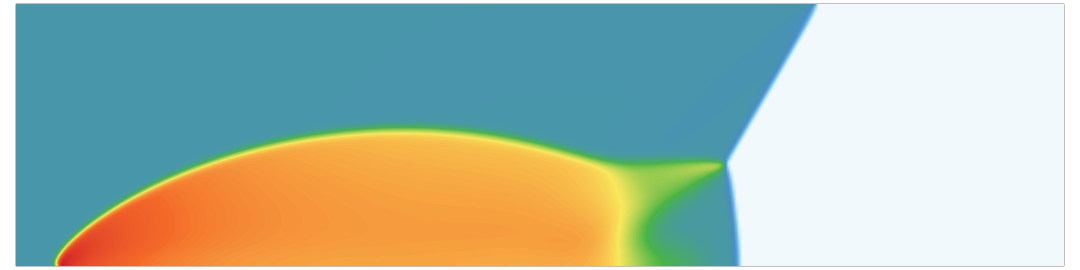}
		\caption{LO, $u_h \in [1.373,20.226]$, $E_h=384\cdot 96$}
		\label{fig:dmrlor4}
	\end{subfigure}
	\begin{subfigure}[b]{.48\linewidth}
		\includegraphics[width=\linewidth]{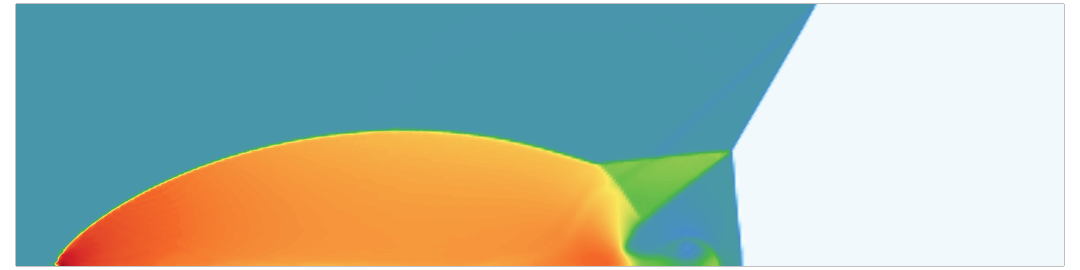}
		\caption{WENO, $u_h \in [1.013,23.613]$, $E_h=384\cdot 96$}
		\label{fig:dmrwenor4}
	\end{subfigure}
	\begin{subfigure}[b]{.48\linewidth}
		\includegraphics[width=\linewidth]{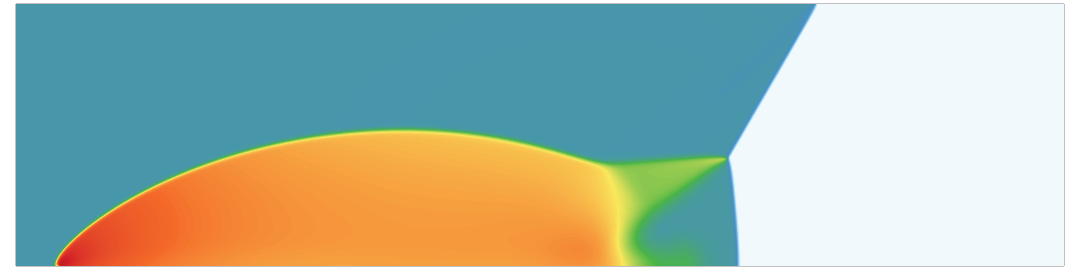}
		\caption{LO, $u_h \in [1.373,20.971]$, $E_h=768\cdot 192$}
		\label{fig:dmrlor5}
	\end{subfigure}
	\begin{subfigure}[b]{.48\linewidth}
		\includegraphics[width=\linewidth]{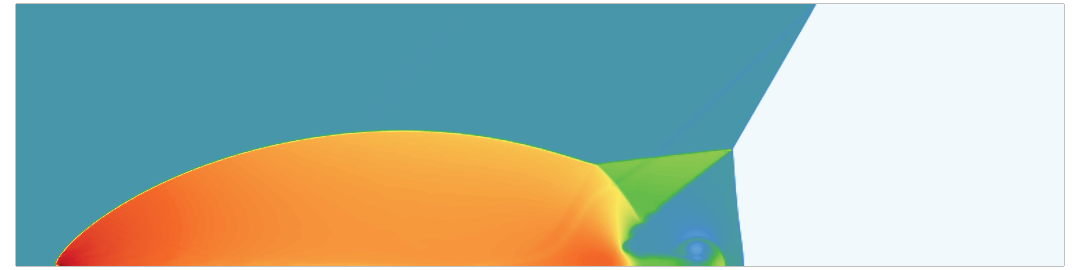}
		\caption{WENO, $u_h \in [0.957,23.201]$, $E_h=768\cdot 192$}
		\label{fig:dmrwenor5}
	\end{subfigure}
	\begin{subfigure}[b]{.48\linewidth}
		\includegraphics[width=\linewidth]{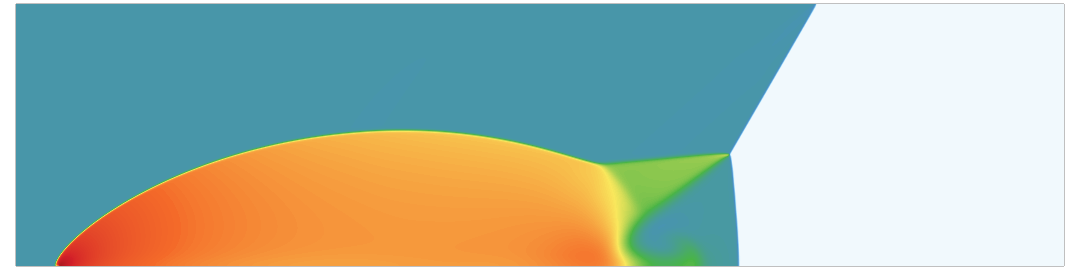}
		\caption{LO, $u_h \in [1.371,21.716]$, $E_h=1536\cdot384$}
		\label{fig:dmrlor6}
	\end{subfigure}
	\begin{subfigure}[b]{.48\linewidth}
		\includegraphics[width=\linewidth]{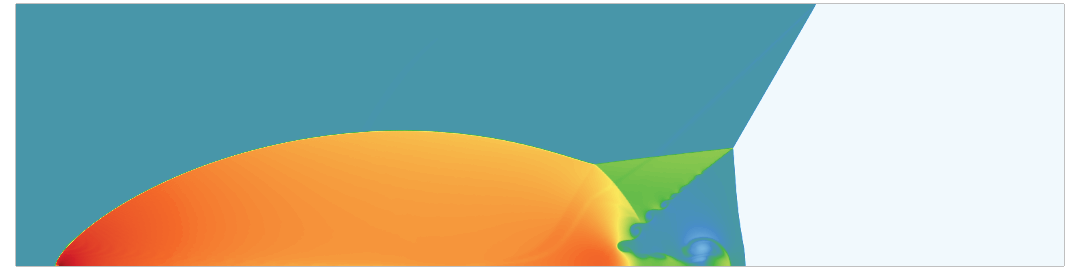}
		\caption{WENO, $u_h \in [0.892,23.148]$, $E_h=1536\cdot384$}
		\label{fig:dmrwenor6}
	\end{subfigure}
	\caption{Double Mach reflection, density profiles $\varrho$ at $t=0.2$ obtained using successively refined meshes and $p=2$.}
	\label{fig:dmr}
\end{figure}

In Fig. \ref{fig:dmr}, we present snapshots of the density distribution at the final time $t=0.2$ obtained using successively refined meshes and $p=2$. The LO approximation exhibits strong numerical diffusion, resulting in a poor resolution of the interacting shock waves. Contrarily, the WENO approximation captures all features sharply while introducing sufficient amounts of numerical dissipation to suppress spurious oscillations.

% Double Mach reflection

\subsubsection{Single-material triple point problem}
We investigate the single-material triple point problem \cite{kucharik2010} on the rectangle $\Omega = (0,7)\times (0,3)$, which is bounded by reflecting walls. The initial data read
\begin{equation}
\begin{bmatrix}
\varrho_1 \\v_{x,1}\\v_{y,1}\\p_1
\end{bmatrix}=
\begin{bmatrix}1.0\\0.0\\0.0\\1.0
\end{bmatrix}, \quad 
\begin{bmatrix}
\varrho_2 \\v_{x,2}\\v_{y,2}\\p_2
\end{bmatrix}=
\begin{bmatrix}0.125\\0.0\\0.0\\0.1
\end{bmatrix}, \quad 
\begin{bmatrix}
\varrho_3 \\v_{x,3}\\v_{y,3}\\p_3
\end{bmatrix}=
\begin{bmatrix}1.0\\0.0\\0.0\\0.1
\end{bmatrix},
\end{equation}
where the subscripts refer to the subdomains $\Omega_1=\{(x,y)\in[0.0,1.0]\times[0.0,3.0]\}$, $\Omega_2=\{(x,y)\in[1.0,7.0]\times[0.0,1.5]\}$ and $\Omega_3=\{(x,y)\in[1.0,7.0]\times[1.5,3.0]\}$, respectively. Initially forming a T-junction, the three states lead to the creation of a shock wave moving to the right, followed by the evolution of a vortex around the triple point.

\begin{figure}[!htb]
	\centering
	\begin{subfigure}[b]{.48\linewidth}
		\includegraphics[width=\linewidth]{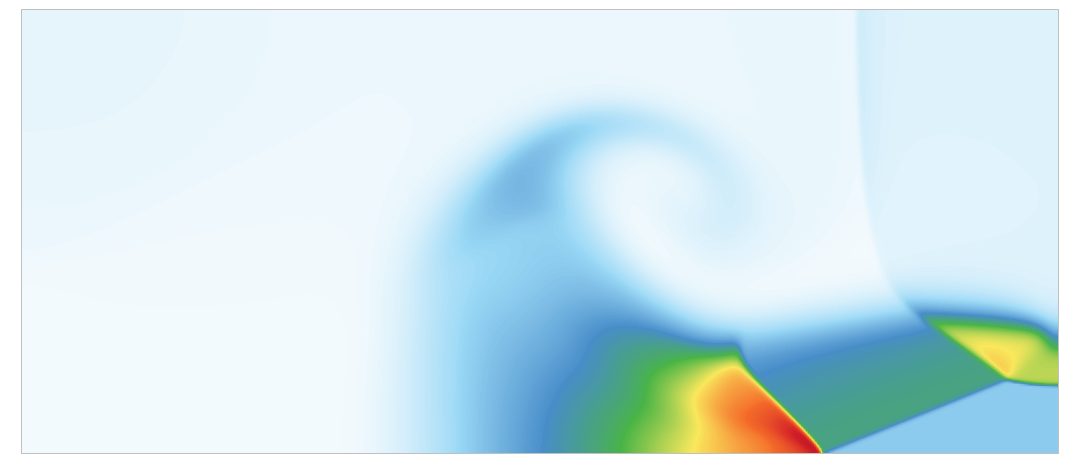}
		\caption{LO, $u_h \in [0.200,4.743]$}
		\label{fig:tpplo}
	\end{subfigure}
	\begin{subfigure}[b]{.48\linewidth}
		\includegraphics[width=\linewidth]{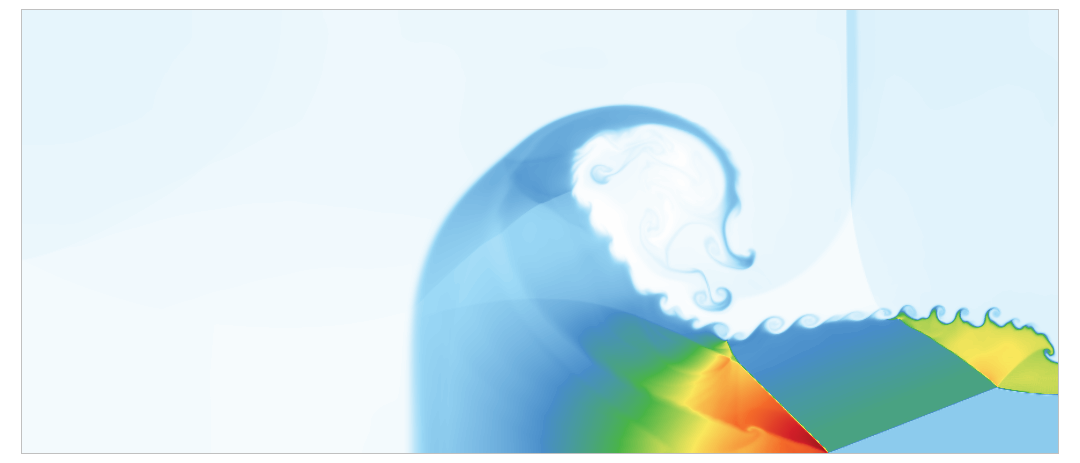}
		\caption{WENO, $u_h \in [0.116,5.271]$}
		\label{fig:tppweno}
	\end{subfigure}
	\caption{Single-material triple point problem, density profiles $\varrho$ at $t=4.0$ obtained using $E_h=896\cdot384$ and $p=2$.}
	\label{fig:tpp}
\end{figure}

In Fig. \ref{fig:tpp}, we present snapshots of the density distribution at the final time $t=4.0$ obtained using $E_h=896\cdot384$ elements and $p=2$. Both the LO and WENO schemes successfully capture the triple point. However, the WENO scheme demonstrates superior resolution of interfacial instability by incorporating smaller-scale structures.

% Single-material triple point problem

\subsubsection{Forward facing step problem}
In our last numerical example, we consider the forward facing step problem \cite{woodward1984}, also known as the Mach 3 wind tunnel test. The wind tunnel is 1 unit wide and 3 units long, featuring a step corner located at $(x,y)=(0.6,0.2)$. Initial conditions are set as $(\varrho, v_x,v_y,p) = (1.4, 3.0, 0.0, 1.0)$, representing a Mach 3 uniform flow impacting the step initially. Supersonic inflow and outflow conditions are applied on the left and right boundaries, respectively, while reflective wall conditions are imposed on the upper and lower boundaries. The corner of the step is treated as a singular point, akin to the approach in \cite{balsara2000,woodward1984}, assuming a nearly steady flow near the corner.

\begin{figure}[!htb]
	\centering
	\begin{subfigure}[b]{.48\linewidth}
		\includegraphics[width=\linewidth]{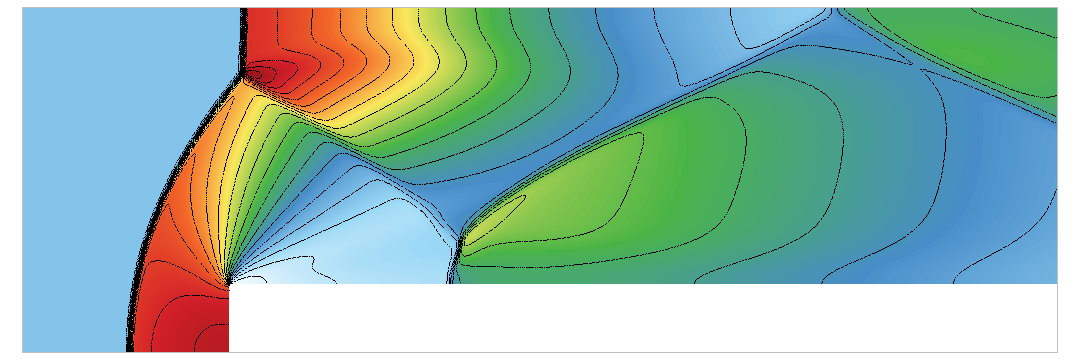}
		\caption{LO, $u_h \in [0.414,6.453]$}
		\label{fig:ffslo}
	\end{subfigure}
	\begin{subfigure}[b]{.48\linewidth}
		\includegraphics[width=\linewidth]{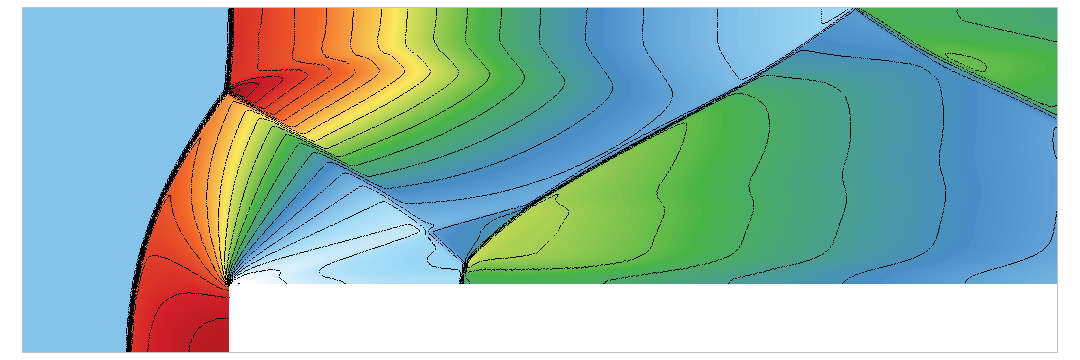}
		\caption{WENO, $u_h \in [0.263,6.508]$}
		\label{fig:ffsweno}
	\end{subfigure}
	\caption{Forward facing step problem, density profiles $\varrho$ at $t=4.0$ obtained using $E_h=480\cdot160$ and $p=2$.}
	\label{fig:ffs}
\end{figure}

We evolve numerical solutions up to the final time $t=4.0$ using $E_h=480\cdot160$ elements and $p=2$. Density profiles obtained from the LO and WENO scheme are shown in Fig. \ref{fig:ffs}. Both schemes successfully capture the shocks, with the WENO scheme providing a solution of higher resolution. 

%The vortex sheet originating from the Mach stem is adequately resolved with only a few zones across the sheet. Unfortunately, it is not.

% Forward facing step problem

%\subsubsection{Shock interaction}
%% Shock interaction
%\begin{figure}[!htb]
%	\centering
%	\begin{subfigure}[b]{.48\linewidth}
%		\includegraphics[width=\linewidth]{results/ssi/blue/ssilop2.png}
%		\caption{LO, $u_h \in [0.138,1.732]$}
%		\label{fig:ssilo}
%	\end{subfigure}
%	\begin{subfigure}[b]{.48\linewidth}
%		\includegraphics[width=\linewidth]{results/ssi/blue/ssiwenop2_dd.png}
%		\caption{DWENO, $u_h \in [0.121,1.758]$}
%		\label{fig:ssiwenodd}
%	\end{subfigure}
%	\caption{Shock interaction, density profiles at $t=0.8$ obtained using $E_h=512^2$ and $p=2$.}
%	\label{fig:ssi}
%\end{figure}

%%%
%%% Conclusions
%%%
\section{Conclusions}
\label{sec:concl}
We have shown that the methodology developed in \cite{kuzmin2023a} for stabilization of continuous Galerkin methods can be extended to and is ideally suited for RKDG discretizations of hyperbolic problems. By introducing low-order nonlinear numerical dissipation based on the relative differences between reconstructed candidate polynomials and the underlying DG approximation, we achieve high-order accuracy while suppressing spurious oscillations in the vicinity of discontinuities. Our algorithm exhibits a modular structure similar to the CG version presented in \cite{kuzmin2023a}. This modular design allows for easy customization of individual components, including the stabilization operator, smoothness sensor, and WENO reconstruction procedure. While our current smoothness sensor is piecewise constant, we aim to develop smoothness sensors inheriting a polynomial structure within each element. This modification will further improve the accuracy of our scheme, especially when dealing with coarse meshes. Theoretical studies of the WENO-based stabilization operator were performed in \cite{kuzmin2023a} in the CG context. We envisage that this preliminary analysis can be readily extended to the DG version.

\medskip

\begin{flushleft}
\textbf{Acknowledgments.} The development of the proposed methodology was sponsored by the German Research Association (DFG) under grant KU 1530/23-3. 
\end{flushleft}
%\section*{References}
\bibliographystyle{plain}
\bibliography{paper_vedral_dg_new}

\end{document}